\documentclass{amsart}
\usepackage{amsmath,amsthm,latexsym,amssymb}
\usepackage{eucal}
\usepackage{mathrsfs}
\usepackage[all]{xy}

\newtheorem{thm}{Theorem}[section]
\newtheorem{lem}[thm]{Lemma}
\newtheorem{prop}[thm]{Proposition}
\newtheorem{cor}[thm]{Corollary}
\newtheorem{schol}[thm]{Scholium}

\newtheorem*{DT}{Cobar Duality Theorem}
\newtheorem*{BDT}{Bar Duality Theorem}
\newtheorem*{AWT}{Alexander-Whitney Theorem}
\theoremstyle{definition}
\newtheorem{defn}[thm]{Definition}
\theoremstyle{remark}
\newtheorem{rmk}[thm]{Remark}
\newtheorem{ex}[thm]{Example}
\newtheorem{exs}[thm]{Examples}

\newtheorem{notn}[thm]{Conventions}

\DeclareMathAlphabet{\mathbfit}{OT1}{cmr}{bx}{it}

\newcommand{\cat}[1]{\mathbf{{#1}}}
\newcommand{\scat}[1]{\mathsf{#1}}
\newcommand{\op}[1]{\mathscr{#1}}

\newcommand{\N}{\mathbf{N}}

\newcommand{\Zmod}[1]{\mathbf{F}_{#1}}
\newcommand{\Hom}{F}
\newcommand{\CoHom}{\operatorname{CoHom}}
\newcommand{\End}{\operatorname{End}}
\newcommand{\Coend}{\operatorname{CoEnd}}
\newcommand{\boxprod}{\circ}                          

\newcommand{\alg}[1]{\op{{#1}}\text{-}\cat{Alg}}
\newcommand{\Alg}[2]{(\op{#1},#2)\text{-}\cat{Alg}}
\newcommand{\sAlg}[2]{(\op{#1},#2)\text{-}\scat{Alg}}
\newcommand{\coalg}[1]{\op{{#1}}\text{-}\cat{Coalg}}
\newcommand{\Coalg}[2]{(\op{#1},#2)\text{-}\cat{Coalg}}
\newcommand{\sCoalg}[2]{(\op{#1},{#2})\text{-}\scat{Coalg}}

\newcommand{\Ind}{\operatorname{Ind}}
\newcommand{\Lin}{\operatorname{Lin}}

\newcommand{\Ssigma}{\op{A}^\perp}
\newcommand{\diffract}{\Phi}
\newcommand{\csg}[1]{\cat{Comon}_{#1}}
\newcommand{\bimod}[2]{{}_{\op{#1}}\cat{Mod}_{\op{#2}}}
\newcommand{\ttma}{\cat{Mod}_{\op A}^{\text{tt}}}
\newcommand{\am}{{}_{\op A}\cat{Mod}}
\newcommand{\msym}{\cat {Ch}^\Sigma_{\text{mult}}}
\newcommand{\comsym}{\cat {Ch}^\Sigma_{\text{comult}}}
\newcommand{\acirc}{\underset {\op A}{\circ}}
\newcommand{\ind}{\operatorname{Ind}}
\newcommand{\lin}{\operatorname{Lin}}
\newcommand\Om{\Omega}
\newcommand{\bicat}{\mathbb}

\newcommand{\ob}{\operatorname{Ob}}
\newcommand{\mor}{\operatorname{Mor}}
\newcommand{\omdc}{\bicat{DC}^\Om}
\newcommand{\phidc}{\bicat{DC}^\Phi}
\newcommand{\del}{\partial}

\begin{document}

\title{Co-rings over operads characterize morphisms}

\author{Kathryn Hess}
\author{Paul-Eug\`ene Parent}
\author{Jonathan Scott}

\address{Institut de g\'eom\'etrie, alg\`ebre et topologie (IGAT) \\
    \'Ecole Polytechnique F\'ed\'erale de Lausanne \\
    CH-1015 Lausanne \\
    Switzerland}
    \email{kathryn.hess@epfl.ch}
\address{Department of Mathematics and Statistics \\
    University of Ottawa \\
    585 King Edward Avenue \\
    Ottawa, ON \\
    K1N 6N5 Canada}
    \email{pauleugene.parent@science.uottawa.ca}
\address{Institut de g\'eom\'etrie, alg\`ebre et topologie (IGAT) \\
    \'Ecole Polytechnique F\'ed\'erale de Lausanne \\
    CH-1015 Lausanne \\
    Switzerland}
    \email{jonathan.scott@epfl.ch}
    \thanks{K.H. thanks the Institut Mittag-Leffler (Djursholm,
Sweden) for its hospitality during a crucial phase of this research.}


\begin{abstract}
Let $\cat M$ be a bicomplete, closed symmetric monoidal category.  Let $\op{P}$ be an operad in $\cat M$, i.e., a monoid in the category of symmetric sequences of objects in $\cat M$, with its composition monoidal structure.  Let $\op R$ be a $\op P$-co-ring, i.e., a comonoid in the category of $\op P$-bimodules.  The co-ring $\op{R}$  induces a natural ``fattening'' of the category
of $\op{P}$-(co)algebras, expanding the morphism sets while
leaving the objects fixed.   Co-rings over operads are thus ``relative operads,'' parametrizing morphisms as operads parametrize (co)algebras.

Let $\op{A}$ denote the associative operad in the category of
chain complexes.  We define a ``diffracting'' functor $\diffract$
that produces $\op{A}$-co-rings from symmetric sequences of chain
coalgebras, leading to a multitude of ``fattened'' categories of
(co)associative  chain (co)algebras.  In particular, we obtain a
purely operadic description of the categories $\cat{DASH}$ and
$\cat{DCSH}$ first defined by Gugenheim and Munkholm, via an
$\op{A}$-co-ring $\op{F}$, which has the two-sided Koszul resolution of $\op A$ as its underlying $\op A$-bimodule.

The diffracting functor plays a crucial role in enabling us to prove existence  of higher, ``up to homotopy'' structure of morphisms via acyclic models methods.  It has already been successfully applied in this sense in \cite {hpst:04}, \cite {hl}, \cite {hess}, \cite {hps2} and \cite {hr}.
\end{abstract}

\maketitle

\section{Introduction}\label{sec:introduction}

Over the past 30 years, operads have proven to be an excellent
tool for encoding the multi-layered structure of certain classes
of algebraic  \emph{objects}, e.g., the coherent systems of higher
homotopies describing $A_{\infty}$ or $E_{\infty}$ structure for
the multiplication map of an algebra.   In this article we study
co-rings over operads, which can be seen as ``relative operads,''  as they serve to parametrize the deep algebraic structure of \emph{morphisms}.

The structure of the categories $\cat{DASH}$ and $\cat{DCSH}$,
first defined by Gugenheim and Munkholm in the early
1970's~\cite{gugenheim-munkholm:74}, motivates our work. The
objects of these categories have a relatively simple algebraic
description, while that of the morphisms is rich and complex.  The
objects of $\cat{DASH}$ are augmented, connected associative chain
algebras, and a morphism from $A$ to $A'$ is a morphism of chain
coalgebras $B(A)\rightarrow B(A')$, where $B$ denotes the bar
construction. Dually, the objects of $\cat{DCSH}$ are counital, connected
coassociative chain coalgebras, and a morphism from $C$ to $C'$ is
a morphism of chain algebras $\Omega C\rightarrow \Omega C'$, where
$\Omega$ denotes the cobar construction.  The research presented here stemmed from our desire to find conditions guaranteeing existence of a morphism of chain coalgebras $B(A)\to B(A')$ extending  a given chain map $A\to A'$ between two chain algebras or, dually, of constructing a morphism of chain algebras $\Omega C\to \Omega C'$ extending a given chain map $C\to C'$ between two chain coalgebras.

The categories $\cat{DASH}$ and $\cat{DCSH}$ play an important
role in topology. Let $C_*$ denote the normalized chains functor
from simplicial sets to chain complexes.  Let $K$ be any $1$-reduced
simplicial set, and let $GK$ be the Kan loop group on $K$. The
usual coproduct on $C_*K$  is a morphism in $\cat{DCSH}$.
Moreover, as  shown in~\cite{hpst:04} using tools developed in this article, if $K$ is $1$-reduced, there is a natural
coproduct on $\Omega C_*K$ with respect to which the natural
equivalence of chain algebras $\Omega C_*K\rightarrow C_*GK$
defined by Szczarba~\cite{szczarba:61} is also a morphism in
$\cat{DCSH}$.  Furthermore,
Bousfield and Gugenheim showed that the Stokes weak equivalence $A_{PL}(X)
\xrightarrow{\simeq} C^{*}(X,\mathbf{Q})$, from the commutative cochain algebra of piecewise-linear rational forms on a space $X$ to its rational singular cochains, is a
$\cat{DASH}$ morphism~\cite{bousfield-gugenheim}. In addition,
Stasheff and Halperin have exploited $\cat{DASH}$-morphisms to
study the collapse of the Eilenberg-Moore spectral
sequence~\cite{stasheff-halperin}.

As an application of the general theory we develop here, we obtain
a purely operadic characterization of the categories $\cat{DASH}$
and $\cat{DCSH}$, in terms of an explicit co-ring $\op{F}$ over
the associative operad $\op{A}$ in the category of chain
complexes. The $\op A$-bimodule underlying $\op{F}$ is in fact the
two-sided Koszul resolution of $\op{A}$.  Since $\op F$ is a semifree $\op A$-co-ring, it is now possible to establish the existence of DASH- and DCSH-structures by acyclic models methods.  This new possibility has already proved invaluable in \cite {hpst:04}, \cite {hps2}, \cite {hl} and \cite {hess}.

We obtain $\op F$  as a special case of a natural ``diffraction'' construction associating a semifree $\op A$-co-ring $\Phi(\op X)$  to any symmetric sequence of chain coalgebras  $\op X$.  We can therefore also use acyclic models methods to prove the existence of parametrized families of higher homotopies that generalize DASH- and DCSH-structures.  In particular, considering the case $\op X=\op F$, we have shown in \cite{hps2} that, if $EK$ is the simplicial suspension of a reduced simplicial set $K$, then there is a natural coproduct on $\Omega^2C_{*}(EK)$ such that the natural weak equivalence of chain algebras $\Omega ^2C_{*}(EK)\to C_{*}(G^2EK)$ is the linear part of a morphism in $\cat {DCSH}\big(\Omega ^2C_{*}(EK), C_{*}(G^2EK)\big)$.

We now provide a more detailed summary of our approach to characterizing morphisms in terms of co-rings over operads.  Let $(\cat{M},\otimes
,I)$ be a bicomplete, closed, symmetric monoidal category.   Let
$\Sigma _n$ denote the symmetric group on $n$ letters. Consider the
category $\cat{M}^{\Sigma}$ of symmetric sequences in $\cat M$, i.e., sequences $\op{X}=\{ \op{X}(n)\mid n\geq 0\}$
of objects in $\cat{M}$  such that $\op{X}(n)$ admits a right
action of $\Sigma_n$ for all $n$.

The category $\cat{M}^{\Sigma}$ admits three distinct monoidal
structures:  the \emph{level} structure $(\cat{M}^{\Sigma},
\otimes, \op{C})$, the \emph{graded} structure $(\cat{M}^{\Sigma},
\odot, \op{U})$, and the \emph{composition} structure
$(\cat{M}^{\Sigma}, \circ, \op J)$;  see
section~\ref{sec:bimodules}.  The level and graded monoidal
structures are both symmetric and closed.  The composition
structure, however, is not symmetric and is only right closed. The
name of the composition structure is justified by the fact that
there is a monoidal functor from $\cat{M}^\Sigma$ to the category
of endofunctors on $\cat{M}$ with monoidal product given by composition;
see~\cite{rezk:96}.

Given the composition monoidal structure, it is easy to define
operads, their (bi)modules and their (co)algebras.  An
\emph{operad} is a unital composition monoid $(\op{P}, \gamma, \eta)$,
i.e., a symmetric sequence $\op{P}$ endowed with an associative
multiplication $\gamma : \op{P}\circ \op{P}\rightarrow \op{P}$
admitting a unit $\eta: \op{J}\rightarrow \op{P}$.  A \emph{left
$\op{P}$-module} consists of a symmetric sequence $\op{M}$ endowed
with a left action $\lambda : \op{P}\circ \op{M}\rightarrow
\op{M}$ satisfying the usual associativity and unit axioms.
\emph{Right $\op{P}$-modules} $(\op{M}, \rho)$ are defined
analogously, and \emph{$(\op{P},\op{Q})$-bimodules} $(\op{R},
\lambda, \rho)$ are symmetric sequences endowed with commuting
left and right actions of $\op{P}$ and $\op{Q}$. We denote the
categories of right and left $\op{P}$-modules and of morphisms of symmetric sequences respecting their module structure by
$\cat{Mod}_{\op{P}}$ and ${}_{\op{P}}\cat{Mod}$.  The category of $(\op P,\op Q)$-bimodules is denoted ${}_{\op{P}}\cat{Mod}_{\op Q}$.  In the bimodule category ${}_{\op{P}}\cat{Mod}_{\op P}$, we can define a monoidal product $\underset {\op P}\circ$ via the usual coequalizer construction.

Let $\op A$ denote the (reduced) \emph{associative operad} in $\cat M$, i.e., for all $n>0$, $\op A(n)=I[\Sigma_{n}]$, the free $\Sigma _{n}$-object generated by the unit object $I$ in $\cat M$, and $\op A(0)=O$, the initial object in $\cat M$.  The composition product $\gamma :\op A\circ\op A\to \op A$ is given by block permutation. In section \ref {ssec:modassocop} of this article, we explain that left $\op A$-modules are graded monoids in $\cat M$ and that right $\op A$-modules are precosimplicial objects in $\cat M$.   In particular, we show that $\op A$-bimodules generalize the notion of  ``operads with multiplication'' due to Gerstenhaber and Voronov~\cite{gerstenhaber-voronov:95}.

The definition of algebras and coalgebras over a given operad
$\op{P}$ is somewhat less concise.  A \emph{$\op{P}$-algebra} is
an object $A$ of $\cat{M}$, together with a set of equivariant
morphisms in $\cat{M}$
\[
    \{ \theta_n : \op{P}(n) \otimes A^{\otimes n} \rightarrow  A
        \mid n \geq 0 \},
\]
where $\Sigma_n$ acts on $A^{\otimes n}$ by permuting factors, and
commuting appropriately with the composition product $\gamma$ on
$\op{P}$. Dually, a \emph{$\op{P}$-coalgebra} consists of an
object $C$ of $\cat{M}$, together with a set of equivariant
morphisms in $\cat{M}$
\[
    \{\theta_n : C \otimes \op{P}(n) \rightarrow
    C^{\otimes n} \mid n\geq 0 \}
\]
commuting appropriately with the composition product $\gamma$. The
categories of $\op{P}$-algebras and $\op{P}$-coalgebras  and of the morphisms in $\cat M$ respecting their (co)algebra structure are
denoted $\alg{P}$ and $\coalg{P}$.

Our characterization of morphisms depends on the observation
that the categories of $\op{P}$-algebras and of
$\op{P}$-coalgebras embed in the categories of left
$\op{P}$-modules and of right $\op{P}$-modules. Questions
concerning \emph{duality} of algebras and coalgebras can therefore
be viewed as questions of \emph{chirality} of right and left
modules, thanks to the asymmetry of the composition product.

The embeddings are defined as follows.  Let $O $ be the
initial object in $\cat M$.   We define two embeddings of $\op{P}$-algebras as
left $\op{P}$-modules.  Given a $\op{P}$-algebra $A$, the
\emph{constant symmetric sequence} $c(A)$ has $c(A)(n)=A$ for all
$n>0$ and $c(A)(0)=O$.  The \emph{trivial symmetric sequence}
$z(A)$ has $z(A)(0)=A$ and $z(A)(n)=O$ for $n \geq 1$.  Both $c(A)$
and $z(A)$ are naturally left $\op{P}$-modules.  On the other
hand, if $C$ is a $\op{P}$-coalgebra, then the \emph{free
graded monoid} $\op{T}(C)$, with $\op{T}(C)(n)=C^{\otimes n}$
for all $n>0$ and $\op{T}(C)(0)=I$, admits a natural right
$\op{P}$-module structure. The algebra embedding $z$ is well-known
(cf. Kapranov-Manin,~\cite{kapranov-manin:01}), and, while we
suspect that the algebra embedding $c$ and the coalgebra embedding
are a part of operad folklore, we were not able to find them in
print.

There are numerous other sources of (bi)modules over operads.  For
example, a morphism of operads $\op{P}\rightarrow \op{Q}$ endows
$\op{Q}$ with the structure of a $\op{P}$-bimodule.  Moreover, if
$\op{X}$ is any symmetric sequence, then $\op{P}\circ \op{X}\circ
\op{Q}$ is naturally a $(\op{P}, \op{Q})$-bimodule.  Finally, as
Ching~\cite{ching:05} and McCarthy and
Minasian~\cite{mccarthy-minasian:04} recently showed, the functor
calculus is a rich source of (bi)modules over operads. Ching
proved  that the derivatives of the identity functor on based
spaces form an operad $\partial_{*}I$ and that any based space
naturally gives rise to a right $\partial_{*}I$-module.  On the
other hand, McCarthy and Minasian explained how to construct an
operad $\mathbf{a}_F$ from any monad $(F,\mu, \eta)$ on
the category of $S$-modules and showed that if the Goodwillie
tower of $F(X)$ splits for every spectrum $X$, then $F(X)$ is an
$\mathbf{a}_F$-algebra, i.e., $c\bigl(F(X)\bigr)$ is a left
$\mathbf{a}_F$-module.

In this article, we need to work with (usually  noncounital) comonoids in the category of $\op P$-bimodules, which are called \emph{$\op P$-co-rings}.  The category of $\op P$-co-rings and of $\op P$-bimodule morphisms respecting the comonoidal structure is denoted $\cat{CoRing}_{\op{P}}$.  The full subcategory of counital co-rings is denoted $\cat{CoRing}_{\op P, *}$. Co-rings over unitary rings were first defined by Sweedler in 1975.  Almost 25 years later, Takeuchi observed that entwined modules, such as Doi-Hopf modules, could be seen as comodules over co-rings, sparking renewed interest in co-rings among algebraists.  Co-rings are now understood to provide a unifying language for various types of Galois descent theory, related to (Hopf-)Galois extensions of commutative rings and to coalgebra Galois extensions.

A $\op{P}$-co-ring $\op{R}$ can be used to ``fatten up'' the categories of left and
right $\op{P}$-modules, as well as those of $\op{P}$-algebras and
$\op{P}$-coalgebras, leaving the objects fixed but expanding the
morphism sets.  As we explain below, the categories $\cat{DASH}$
and $\cat{DCSH}$ are ``fattened" versions of $\alg{A}$ and
$\coalg{A}$, where $\cat{M}=\cat {Ch}$ is the category of chain complexes over a commutative ring $R$.

Let $\psi: \op{R}\rightarrow \op{R} \underset{\op{P}}{\circ}
\op{R}$ denote the coproduct on $\op{R}$, which is coassociative.
Assume that it is counital with respect to a $\op{P}$-bimodule morphism
$\varepsilon: \op{R} \rightarrow \op{J}$.  We then define the
fattened categories ${}_{(\op{P},\op R)}\cat{Mod}$,
$\cat{Mod}_{(\op{P}, \op R)}$, $\Alg{P}{\op R}$, and
$\Coalg{P}{\op R}$ to have the same objects as their thinner
counterparts, but to have morphisms given by
\[
    {}_{(\op{P},\op R)}\cat{Mod}(\op{M}, \op{N})
    :={}_{\op{P}}\cat{Mod}(\op{R}\underset{\op P}{\circ} \op M, \op{N});
\]
and
\[
    \Alg{P}{\op R}(A,A') := {}_{(\op{P},\psi)}\cat{Mod}\left(c(A), c(A')\right)={}_{\op{P}}\cat{Mod}\left(\op{R}\underset{\op P}{\circ} c(A), c(A')\right);
\]
similarly for right $\op{P}$-modules and $\op{P}$-coalgebras.

The composition of morphisms in these categories is defined  in
terms of $\psi$.  Given $\theta \in
{}_{(\op{P},\op R)}\cat{Mod}(\op{M}, \op{M}')$ and $\theta'\in
{}_{(\op{P},\op R)}\cat{Mod}(\op{M}', \op{M}'')$, their composite
$\theta'\theta\in{}_{(\op{P},\op R)}\cat{Mod}(\op{M}, \op{M}'')$ is
given by composing the following sequence of (strict) morphisms of
left $\op{P}$-modules.
\[
    \op{R} \underset{\op{P}}{\circ} \op{M}
        \xrightarrow{\psi\underset{\op{P}}{\circ} Id_{\op{M}}}
        \op{R} \underset{\op{P}}{\circ} \op{R} \underset{\op{P}}{\circ} \op{M}
        \xrightarrow{Id_{\op{R}}\underset{\op{P}}{\circ} \theta}
        \op{R} \underset{\op{P}}{\circ} \op{M}'
        \xrightarrow{\theta'}
        \op{M}''
\]
For all right $\op P$-modules $\op M$, the identity morphism in ${}_{(\op{P},\psi)}\cat{Mod}(\op{M}, \op{M})$ is the strict morphism of $\op P$-modules
$$\op R\underset {\op P}\circ \op M\xrightarrow{\varepsilon\underset {\op P}\circ id_{\op M}}\op P\underset {\op P}\circ \op M\cong \op M.$$
 Composition and identities in $\cat{Mod}_{(\op{P},\op R)}$ are defined similarly,
while composition in $\Alg{P}{\op R}$ and in $\Coalg{P}{\op R}$ is
obtained by restriction from  ${}_{(\op{P},\op R)}\cat{Mod}$ and
$\cat{Mod}_{(\op{P},\op R)}$. We call
${}_{(\op{P},\op R)}\cat{Mod}$, $\cat{Mod}_{(\op{P}, \op R)}$,
$\Alg{P}{\op R}$, and $\Coalg{P}{\op R}$ the
\emph{$\op{R}$-governed} versions of their strict
counterparts.  Morphisms in an $\op R$-governed category are denoted $\underset {\op R}\to$.

Note that ${}_{\op{P}}\cat{Mod}(\op{M}, \op{N})$ embeds naturally
in ${}_{(\op{P},\op R)}\cat{Mod}(\op{M}, \op{N})$, by sending
$\varphi$ to
\[
    \varepsilon \underset{\op{P}}{\circ} \varphi
        : \op{R} \underset{\op{P}}{\circ} \op{M} \rightarrow
        \op{P} \underset{\op{P}}{\circ} \op{N}\cong \op{N}.
\]
Similar embeddings exist of the ``strict'' right module, algebra
and coalgebra categories into their ``fattened'' versions.

The reader familiar with category theory will have recognized that
the $\op{R}$-governed categories of left and right
$\op{P}$-modules are precisely the Kleisli categories associated
to the comonads  $\op{R}\underset{\op{P}}{\circ}-$ and
$-\underset{\op{P}}{\circ} \op{R}$, respectively.

If the co-ring $\op R$ is not counital, then there are no canonical candidates for identity morphisms, though the associative composition of $\op R$-governed morphisms still makes sense.  We obtain therefore $\op R$-governed \emph{semicategories} of modules $ {}_{(\op{P},\op R)}\scat{Mod}$ and $\scat{Mod}_{(\op{P},\op R)}$ and corresponding semicategories of (co)algebras.

The definitions made above are clearly canonical, i.e., if $\cat M$ is a small category and $\op P$ is an operad in $\cat M$, there are functors
$$L:\cat{CoRing}_{\op P}^{op}\to \cat {SemiCat}: \op R\mapsto {}_{(\op{P},\op R)}\scat{Mod}$$
$$R:\cat{CoRing}_{\op P}^{op}\to \cat {SemiCat}: \op R\mapsto \scat{Mod}_{(\op{P},\op R)}$$
and
$$L:\cat{CoRing}_{\op P,*}^{op}\to \cat {Cat}: \op R\mapsto {}_{(\op{P},\op R)}\cat{Mod}$$
$$R:\cat{CoRing}_{\op P,*}^{op}\to \cat {Cat}: \op R\mapsto \cat{Mod}_{(\op{P},\op R)},$$
where $\cat {SemiCat}$ is the category of small semicategories and $\cat {Cat}$ is the category of small categories.

A plentiful source of composition comonoids is thus essential to
producing ``fattened" (semi)categories of modules over an operad.  In
this article we develop a tool for constructing composition
comonoids over the associative operad $\op{A}$, when $\cat{M}=\cat {Ch}$ is
the category of connective (i.e., bounded-below) chain complexes over a commutative ring $R$.  Let
$\cat{Comon}_{\otimes}$ denote the category of level comonoids in
$\cat{M}^{\Sigma}$, i.e., of symmetric sequences $\op{M}$ endowed
with a coassociative level comultiplication $\Delta: \op{M}
\rightarrow \op{M}\otimes \op{M}$ that is not necessarily
counital. We define a functor, called the \emph{diffracting
functor},
\[
    \diffract : \cat {Comon}_{\otimes} \rightarrow
       \cat{CoRing}_{\op{A}}.
\]
As a consequence, there exist
$\diffract(\op{M})$-governed semicategories of left
and right $\op{A}$-modules and of $\op{A}$-algebras and
$\op{A}$-coalgebras.

There is a
natural transformation, which we call the \emph{generalized
Milgram map},
\[
    q : \diffract(-\otimes -) \Longrightarrow
        \diffract(-)\otimes \diffract(-)
\]
of functors from $\cat {Comon}_{\otimes}^{\times 2}$ to
${}_{\op{A}}\cat{Mod}_{\op{A}}$.  The transformation $q$ is a sort
of generalized Alexander-Whitney map, along the same lines as the
original Milgram map $\Omega (C\otimes C')\rightarrow  \Omega
C\otimes \Omega C'$ for chain coalgebras $C$ and
$C'$~\cite{milgram:66}. Consequently, if the level
comultiplication $\Delta$ on $\op{M}$ is itself a morphism of
level comonoids, then $\diffract(\op{M})$ is naturally endowed with a level comultiplication
$q\diffract(\Delta)$, which is not necessarily coassociative.

The diffracting functor satisfies a beautiful duality relation,
with respect to both the cobar construction and the bar
construction.  Before stating the duality results, we fix some basic terminology and notation.

Let $\op{M}$ be a level comonoid.  If $\op{Y}$
and $\op{Z}$ are level (co)monoids, then a morphism $\op{Y}
\circ \op{M} \rightarrow \op{Z}$ is called
\emph{(co)multiplicative} if it is compatible with the
(co)multiplications in $\op{Y}$ and $\op{Z}$ and with the
comultiplication in $\op{M}$; see
section~\ref{sec:induction-linearization} for the precise
definitions. Let $\cat M^\Sigma_{\text{mult}}(\op Y\circ \op M, \op Z)$ denote the set of multiplicative morphisms of symmetric sequences, when $\op Y$ and $\op Z$ level monoids.  Let $\comsym (\op Y\circ\op M,  \op Z)$ denote the set of comultiplicative morphisms of symmetric sequence, when $\op Y$ and $\op Z$ are level comonoids.   We show in section \ref {ssec:bimodules} that if $A$ is a monoid, then $\op T(A)$ is a level monoid, while if $B$ is a comonoid, then $c(B)$ is a level comonoid.

In the statement below, $\cat C^{\cat D}$ denotes the category of functors and natural transformations from a small category $\cat D$ to a category $\cat C$.

\begin{BDT}
Let $\cat D$ be any small category.
There are mutually inverse, natural isomorphisms
$$\ind: \am\big(\Phi(-)\acirc c(-),c(-)\big)\longrightarrow \comsym \big(-\circ \,c(B-),c(B-)\big),$$
and
$$\lin :\comsym \big(-\circ\, c(B-),c(B-)\big)\longrightarrow\am\big(\Phi(-)\acirc c(-),c(-)\big)$$
of  functors from $(\alg{A})^{\cat D}\times\cat {Comon}_{\otimes}\times (\alg{A})^{\cat D}$ to  $\cat{Set}^{\cat D}.$
\end{BDT}

 The natural isomorphisms $\ind$ and $\lin$ are called \emph{induction} and \emph{linearization}, respectively.

Let $A$ and $A'$ be any two associative chain algebras.  By the Bar Duality Theorem, an $\op{M}$-governed morphism $BA
\xrightarrow[\op{M}]{} BA'$ can be obtained by suspending and then
comultiplicatively lifting a family
\[
    \{ \diffract(\op{M})(m)\otimes A^{\otimes m}
        \rightarrow A' \mid m\geq 1\}
\]
of appropriately equivariant morphisms of chain complexes.   On
the other hand, by restriction and desuspension, an
$\op{M}$-governed morphism $BA \xrightarrow[\op{M}]{} BA'$ gives
rise to a $\Phi(\op{M})$-governed morphism of chain algebras $A
\xrightarrow[\Phi(\op{M})]{} A'$.

As an application of Bar Duality, we can prove the existence of
natural $\op{X}$-governed morphisms between bar constructions by
acyclic models methods, since $\diffract(\op{X})$ is a semifree $\op A$-co-ring.  This fundamental and extremely useful existence result accounts for our interest in the diffracting functor.

\begin{thm}[Existence of $\op{M}$-governed morphisms
I]\label{thm:exist-1} Let $X,Y: \cat D\to\alg A$ be functors, where $\cat D$ is a category admitting a set of models $\mathfrak M$ with respect to which  $X$ is free and globally connective and $Y$ is acyclic.  Let $(\op M, \Delta)$ be a level comonoid under $\op J$.  Let $\tau:UX\to UY$ be a natural transformation, where $U $ is the forgetful functor down to $\cat {Ch}$.  Then there is a natural, comultiplicative  transformation
$$\theta:\op T(B X)\circ \op M\to \op T(B Y)$$
lifting $s\tau$,  i.e., the following composite is equal to $\tau$.
$$X\xrightarrow{s}sX\hookrightarrow \op c(B X)\circ \op J\to c(BX)\circ \op M\xrightarrow{\theta}c(BY)\xrightarrow {\text{proj.}}sY\xrightarrow {s^{-1}}Y$$
\end{thm}

In other words, for each $d\in \ob \cat D$, the natural chain map $\tau(d):F(d)\to G(d)$ extends naturally to an $\op M$-governed coalgebra morphism $\theta_{\op M}(d):BF(d)\underset {\op X}\to BG(d)$.

We prove as well an Eckmann-Hilton dual to the Bar Duality Theorem, which is absolutely essential to the work in \cite {hpst:04}, \cite {hl}, \cite {hess} and \cite {hps2}. Here, $\cat{Mod}_{\op A}^{\text{tt}}\big(\op T(X)\acirc\Phi( \op X) ,\op T(Y)\big)$ denotes the set of transposed tensor morphisms of right $\op A$-modules (cf., Definition \ref {def:tt-morph}), where $\op T(X)$ and $\op T(Y)$ are both right $\op A$-modules.

\begin{DT} Let $\cat D$ be any small category.
There are mutually inverse, natural isomorphisms
$$\ind: \ttma\big(\op T(-)\acirc \Phi(-), \op T(-)\big)\longrightarrow \msym \big(\op T(\Om-)\circ -,\op T(\Om-)\big)$$
and
$$\lin :\msym \big(\op T(\Om-)\circ -,\op T(\Om-)\big)\longrightarrow\ttma\big(\op T(-)\acirc \Phi(-),\op T(-)\big)$$
of  functors from $(\coalg{A})^{\cat D}\times\cat {Comon}_{\otimes}\times (\coalg{A})^{\cat D}$ to  $\cat{Set}^{\cat D}.$
\end{DT}

 The natural isomorphisms $\ind$ and $\lin$ are again called \emph{induction} and \emph{linearization}, respectively.

Let $C$ and $C'$ be any two coassociative chain coalgebras.  The
Cobar Duality Theorem implies that an $\op{M}$-governed morphism
$\Omega C \xrightarrow[\op{M}]{} \Omega C$ can be obtained by
desuspending and then multiplicatively extending a family
\[
    \{ C \otimes \diffract(\op{M})(m) \rightarrow
    (C')^{\otimes m} \mid m \geq 1\}
\]
of appropriately equivariant morphisms of chain complexes. On the
other hand, by restriction and suspension,  an $\op{M}$-governed
morphism $\Omega C \xrightarrow[\op{M}]{} \Omega C'$ gives rise to
a $\diffract(\op{M})$-governed morphism of chain coalgebras $C
\xrightarrow[\diffract(\op{M})]{} C'$.

Acyclic models methods again permit us to establish the existence
of natural $\op{M}$-governed morphisms, now between cobar
constructions.

\begin{thm}[Existence of $\op{M}$-governed morphisms
II]\label{thm:exist-2} Let $X,Y: \cat D\to\coalg A$ be functors, where $\cat D$ is a category admitting a set of models $\frak M$ with respect to which  $X$ is free and globally connective and $Y$ is acyclic.  Let $(\op M, \Delta)$ be a level comonoid under $\op J$.  Let $\tau:UX\to UY$ be a natural transformation, where $U $ is the forgetful functor down to $\cat {Ch}$.  Then there is a natural, multiplicative  transformation
$$\theta:\op T(\Om X)\circ \op M\to \op T(\Om Y)$$
extending $s^{-1}\tau$,  i.e., the following composite is equal to $\tau$.
$$X\xrightarrow{s^{-1}}s^{-1}X\hookrightarrow \op T(\Om X)\circ \op J\to \op T(\Om X)\circ \op M\xrightarrow{\theta}\op T(\Om Y)\xrightarrow {\text{proj.}}s^{-1}Y\xrightarrow {s}Y$$
\end{thm}

In other words, for each $d\in \ob \cat D$, the natural chain map $\tau(d):F(d)\to G(d)$ extends naturally to an $\op M$-governed algebra morphism $\theta_{\op M}(d):\Om F(d)\underset {\op M}\to \Om G(d)$.

We prove an Enriched Cobar Duality Theorem in section \ref {ssec:bicatbndl-coalg}, which we then use to show that the induction isomorphism $\ind$ can be made to transport extra structure, in the following sense. Let  $\theta: \op T(C)\acirc \Phi(\op M)\to \op T(C')$ be a transposed tensor morphism, where $C$ and $C'$ are chain coalgebras,  and $\op M$ is a level comonoid.  If $\op M$ is a right $\op P$-module and $\Om C'$ is a $\op P$-coalgebra it is natural to ask when $\ind \theta:\op T(\Om C)\circ \op M\to\op T(\Om C')$ is a morphism of right $\op P$-modules.  Similarly, when $\op M$ is a left $\op P$-module and $\Om C$ is a $\op P$-coalgebra, we can ask when $\ind\theta: \op T(\Om C)\circ \op M\to \op T(\Om C')$ induces a morphism $\widehat{\ind\theta}: \op T(\Om C)\underset {\op P}\circ \op M\to \op T(\Om C')$.   We provide complete answers to these questions in section \ref {sec:enrind}, in terms of Enriched Cobar Duality.

Applying enriched induction, we show in Theorems \ref{thm:diffmod} and \ref{thm:diffmod2} that acyclic models methods can again be used, to prove the existence of right $\op P$-module maps $ \op T(\Om C)\underset {\op P}\circ \op M\to \op T(\Om C')$.  This result is absolutely essential to proving one of the main theorems in \cite {hps2}.

The Enriched Cobar Duality Theorem is most clearly expressed in the language of \emph{bundles of bicategories with connection}, a new categorical concept that is analogous to the geometric notion of fiber bundles with connection.  We present the elements of this theory in section \ref {sec:enrcobdual}.  A complete introduction to the theory will appear in \cite {hess:bicatbndl}.

Evaluating $\diffract$ on the
composition unit $\op{J}$, we obtain $\op{F}=\diffract(\op{J})$, the
\emph{Alexander-Whitney} $\op{A}$-co-ring, which plays a particularly important role.  Let
$\psi_{\op{F}} : \op{F} \rightarrow \op{F}
\underset{\op{A}}{\circ} \op{F}$ denote the composition
comultiplication on $\op{F}$.  As we show in section \ref{sec:aw}, the coproduct $\psi_{\op F}$ admits a counit $\op F \to \op A$, so that there are categories ${}_{(\op{A},\psi_{\op F})}\cat{Mod}$ and $\cat{Mod}_{(\op{A},\psi_{\op F})}$ and not just semicategories.

Note that $\op{F}$ admits the structure of a level comonoid,
since the level comultiplication $\Delta_{\op{J}}$ on $\op{J}$ is
necessarily a morphism of level comonoids.  We can therefore
apply the generalized Milgram transformation to obtain a level
comultiplication $\Delta _{\op{F}}=q\diffract (\Delta _{\op{J}})$
on $\op{F}$, which we prove to be coassociative.

The following theorem summarizes the most important properties of
$\op{F}$, which are proved in section \ref{sec:aw} and section 7.

\begin{AWT}  Let $\op{F} = \diffract{\op{J}}$.
    \begin{enumerate}
    \item $\cat{DASH} = \Alg{A}{\op{F}}$ and
    $\cat{DCSH} = \Coalg{A}{\op{F}}$.
    \medskip
    \item The natural morphism $\op{F}\rightarrow \op{A}$ of
    $\op{A}$-bimodules is a quasi-isomorphism in positive levels, i.e.,
    $\op{F}$ is a free $\op{A}$-bimodule resolution of $\op{A}$.  In fact, $\op F$ is the two-sided Koszul resolution of $\op A$.
     \medskip
    \item Let $X, Y: \cat D\to\alg A$  be functors, where $\cat D$ is a category admitting a set of models $\frak M$ with respect to which  $X$ is free and globally connective and $Y$ is acyclic.  Let $\tau:UX\to UY$ be a natural transformation, where $U $ is the forgetful functor down to $\cat {Ch}$.  Then there is a natural, comultiplicative  transformation
$$\theta:BX\to BY$$
extending $s\tau$,  i.e., the following composite is equal to $\tau$.
$$X\xrightarrow{s}sX\hookrightarrow BX\xrightarrow{\theta}BY\xrightarrow {\text{proj.}}sY\xrightarrow {s^{-1}}Y$$
i.e., for each $d\in \ob \cat D$, the natural chain map $\tau(d):F(d)\to G(d)$ admits a natural DASH-structure.

Furthermore, if $B X$ and $BY$ admit natural, associative products, then there is a comultiplicative transformation
$$\Theta:B^2 X\to B^2Y$$
extending $s\theta$.
     \medskip
    \item Let $X, Y: \cat D\to\coalg A$  be functors, where $\cat D$ is a category admitting a set of models $\frak M$ with respect to which  $X$ is free and globally connective and $Y$ is acyclic.  Let $\tau:UX\to UY$ be a natural transformation, where $U $ is the forgetful functor down to $\cat {Ch}$.  Then there is a natural, multiplicative  transformation
$$\theta:\Om X\to \Om Y$$
extending $s^{-1}\tau$,  i.e., the following composite is equal to $\tau$.
$$X\xrightarrow{s^{-1}}s^{-1}X\hookrightarrow \Om  X\xrightarrow{\theta}\Om  Y\xrightarrow {\text{proj.}}s^{-1}Y\xrightarrow {s}Y$$
i.e., for each $d\in \ob \cat D$, the natural chain map $\tau(d):F(d)\to G(d)$ admits a natural DCSH-structure.

Furthermore, if $\Om X$ and $\Om Y$ admit natural, coassociative coproducts, then there is a multiplicative transformation
$$\Theta:\Om ^2 X\to \Om ^2Y$$
extending $s^{-1}\theta$.
\end{enumerate}
\end{AWT}

The first property asserts that we have attained our goal of
providing a purely operadic description of the categories
$\cat{DASH}$ and $\cat{DCSH}$ and follows immediately from the Bar
and Cobar Duality Theorems, applied to $\op{X}=\op{J}$.

In Example~\ref{ex:counter}, we use $\op{F}$ to provide an example
of a chain coalgebra whose cohomology algebra is realizable as an
algebra over the Steenrod algebra, but which is not of the
homotopy type of the chain complex of any space.

In future articles we will treat numerous possible
generalizations  of the diffracting functor.  In particular, we
will define and apply diffraction over general quadratic operads $\op Q$ and their cofibrant replacements $\op Q_{\infty}$.

The authors would like to express their heartfelt appreciation to
Haynes Miller, for having mentioned the words ``Kleisli category''
at a critical juncture during the course of the research presented
here.  The first author would like to thank the Clay Mathematics
Institute and the organizers of the Workshop on Calculus of
Functors and Operads for having been invited to participate in the
workshop, at which she first learned about the remarkable
suspension Hopf operad of Arone, Bauer, Johnson and Morava, which
proved so critical to developing a deep understanding of
diffraction.  The second author would like to thank the IGAT at
the EPFL for its hospitality during the visits he made during the
course of this research.

\begin{notn}
We denote by $\N$ the set of natural numbers $\{ 0, 1, 2, \ldots
\}$. Let $\vec{m} = (m_{1}, \ldots, m_{k}) \in \N^{k}$. We use the
convention that $m = m_{1} + \cdots + m_{k}$. We set $I_{k,m} = \{
\vec{m} \mid m_{1} + \cdots + m_{k} = m \}$ and $I_{k,m}^{+} = \{
\vec{m} \in I_{k,m} \mid m_{i} > 0\ \forall\ i \}$.

Given objects $A$ and $B$ of a category $\cat C$, respectively of a semicategory $\scat C$, we let $\cat C(A,B)$, respectively $\scat C(A,B)$, denote the set of morphisms with source $A$ and target $B$.

Let $(\cat M,\otimes, I)$ be a monoidal category.  To simplify the language in this paper, we say that a \emph{monoid} in $\cat M$ consists of an object $A$ endowed with an associative multiplication morphism $\mu:A\otimes A\to A$ and denote the category of monoids and their homomorphisms $\cat {Mon}$.  If there is a morphism $\eta: I\to A$ that is a unit with respect to  the multiplication map $\mu$, i.e., $\mu (\eta\otimes Id_{A})=Id_{A}=\mu (Id_{A}\otimes \eta)$, then we call $(A,\mu,\eta)$ a \emph{unital monoid}.  The category of unital monoids is denoted $\cat {Mon}_{*}$.

In the dual case we distinguish similarly between comonoids and counital comonoids, the categories of which are denoted $\cat {Comon}$ and $\cat {Comon}_{*}$.

In this paper we often work in the category $\cat{Ch}$ of connective chain complexes over a commutative ring $R$ and of chain maps between them.  Borrowing terminology from stable homotopy theory, we call a chain complex $C$ \emph{connective} if there is an integer $N$ such that $C_{n}=0$ for all $n<N$.  The tensor product $C\otimes C'$ of connective chain complexes $C$ and $C'$ satisfies
$$(C\otimes C')_{n}= \bigoplus _{l+m=n}C_{l}\otimes C'_{m},$$
which is a finite sum for all $n$ since $C$ and $C'$ are connective. The differential $d''$ on $C\otimes C'$ is given by the usual Leibniz rule, i.e., $d''(c\otimes c')=dc\otimes 1 +(-1)^{|c|}c\otimes d'c'$, where $|c|$ denotes the degree of $c$.  The sign in the definition of $d''$ is there to ensure that $d''\circ d''=0$.

Another convention used consistently throughout this article is the Koszul sign convention for commuting elements  of a graded module or for commuting a morphism of graded modules past an element of the source module.  For example,  if $V$ and $W$ are graded algebras and $v\otimes w, v'\otimes w'\in V\otimes W$, then
$$(v\otimes w)\cdot (v'\otimes w')=(-1)^{|w|\cdot |v'|}vv'\otimes ww'.$$
Futhermore, if $f:V\to V'$ and $g:W\to W'$ are morphisms of graded modules, then for all $v\otimes w\in V\otimes W$,
$$(f\otimes g)(v\otimes w)=(-1)^{|g|\cdot |v|} f(v)\otimes g(w).$$
The source of the Koszul sign convention is the definition of the twisting isomorphism
$$\tau: V\otimes W\longrightarrow W\otimes V: v\otimes w\mapsto (-1)^{|v|\cdot |w|}w\otimes v.$$
The sign in the definition of the twisting isomorphism is chosen so that $\tau$ is a chain map if $V$ and $W$ are chain complexes.

To simplify presentation, all signs in this article are implicit, as they all follow simply from applying the Koszul rule.
\end{notn}

\tableofcontents

\section{Operads and their modules}\label{sec:bimodules}

In this section we first recall the category of symmetric sequences and
review its various monoidal structures. We then recall the definition
of operads, which are monoids with respect to one of the monoidal structures on the category of symmetric sequences.

As a monoid, an operad
has categories of left and right modules, as well as bimodules. We
recall how left modules generalize algebras, and observe that
right modules generalize coalgebras.  In the section 2.3 we discuss how
co-rings, which are bimodules with extra comonoidal structure, can be seen as generalizing morphisms between
(co)algebras.  Modules over the associative operad $\op{A}$ enjoy
certain special properties that are considered separately in section \ref{ssec:modassocop}.  In
particular, an $\op{A}$-bimodule is a cosimplicial object with
cup-pairing in the sense of McClure and
Smith~\cite{mcclure-smith:02}.

For a more thorough treatment
of the various monoidal products of symmetric sequences, and
operads as monoids, we refer the reader to Fresse~\cite{fresse:04}
and Rezk~\cite{rezk:96}.

\subsection{Monoidal products for symmetric sequences}\label{ssec:sigma-modules}

Let $(\cat{M}, \otimes, I)$ be a closed symmetric monoidal
category that is complete and cocomplete. In particular, $\cat{M}$
has an initial object $O$ that satisfies $O \amalg A \cong A
\amalg O \cong A$ for all objects $A$ in $\cat{M}$ (since $A$
satisfies the relevant universal property) and $O \otimes A \cong
A \otimes O \cong O$ (since $\otimes$ commutes with colimits and
$O$ is the colimit of the empty diagram in $\cat{M}$).

A \emph{symmetric sequence} consists of a sequence $\op{X} = \{
\op{X}(n) \mid n \geq 0 \}$ of objects in $\cat{M}$, where each
$\op{X}(n)$ comes equipped with a right action of the symmetric
group $\Sigma_{n}$, that is, an anti-homomorphism $\Sigma_{n}
\rightarrow \mathrm{Aut}_{\cat{M}}(\op{X}(n))$.  We call
$\op{X}(n)$ the \emph{$n$th level} of $\op{X}$.  A morphism of
symmetric sequences $\varphi:\op{X} \rightarrow \op{Y}$ is a
sequence of morphisms $\{ \varphi_{n} : \op{X}(n) \rightarrow
\op{Y}(n) \mid n \geq 0 \}$, where $\varphi_{n}$ is
$\Sigma_{n}$-equivariant. Abusing notation slightly, we denote the category of symmetric sequences and
their morphisms by
$\cat{M}^{\Sigma}$, rather than $\cat{M}^{(\Sigma^{\text{op}})}$.

\begin{defn}\label{defn:level-tensor}
The \emph{level tensor product} of two symmetric sequences
$\op{X}$ and $\op{Y}$ is defined level-wise, that is,
\[
    ( \op{X} \otimes \op{Y} )(n) = \op{X}(n) \otimes \op{Y}(n)
    \quad (n \geq 0)
\]
with the diagonal action of $\Sigma_{n}$.
\end{defn}

The following proposition is easily proved.

\begin{prop}\label{prop:level-tensor}
Let $\op{C} = \{ \op{C}(n) \}_{n \geq 0}$ be the symmetric
sequence with $\op{C}(n) = I$ and trivial $\Sigma_{n}$-action,
for all $n \geq 0$.  Then $(\cat{M}^{\Sigma},\otimes,\op{C})$ is a
closed symmetric monoidal category, called the \emph{level
monoidal structure} on $\cat{M}^{\Sigma}$.
\end{prop}

\begin{rmk}
We use the symbol $\op{C}$ for the unit since it happens to be the
symmetric sequence underlying the commutative operad.
\end{rmk}

A \emph{level monoid} is a monoid with respect to the
level tensor product.  Thus $\op{X}$ is a level monoid if and
only if each $\op{X}(n)$ comes equipped with a
$\Sigma_{n}$-equivariant associative multiplication.  A
\emph{unital level monoid} is a level monoid with unit $\eta : \op{C}
\rightarrow \op{X}$. Dually, we have \emph{level comonoids} and
\emph{counital level comonoids}.

Let $G$ be a discrete group.  The forgetful functor $\cat{M}^{G}
\rightarrow \cat{M}$ from the category of $G$-objects has a left
adjoint, denoted $(-)[G]$.  If $X \in \ob\cat{M}$, then $X[G] :=
\coprod_{g \in G} X_{g}$, where $X_{g}$ is a copy of $X$ labelled
by $g \in G$.  We call $X[G]$, equipped with the obvious right
$G$-action of permutation of summands, the \emph{free right
$G$-object on $X$}.  Note that there is an isomorphism $X[G \times
H] \cong X[G] \otimes X[H]$, natural in $X$, $G$, and $H$.

\begin{ex}\label{ex:ass-as-level-comonoid}
Define the symmetric sequence $\op{A}$ by $\op{A}(n) =
I[\Sigma_{n}]$ for all $n \geq 1$, with $\op{A}(0) = O$.  For $n
\geq 1$, the diagonal $\Sigma_{n} \rightarrow \Sigma_{n} \times
\Sigma_{n}$ determines a diagonal $I[\Sigma_{n}] \rightarrow
I[\Sigma_{n}] \otimes I[\Sigma_{n}]$ that is counital with respect
to the augmentation $I[\Sigma_{n}] \rightarrow I$. It follows that
$\op{A}$ is a counital level comonoid. Furthermore, if $i + j = m$, then
$I[\Sigma_{m}]$ is a left $(\Sigma_{i} \times \Sigma_{j})$-module,
considering $\Sigma_{i} \times \Sigma_{j}$ as a subgroup of
$\Sigma_{m}$ in the obvious way.
\end{ex}

\begin{defn}\label{defn:graded-tensor}
Let $\op{X}$ and $\op{Y}$ be symmetric sequences. We define their
\emph{graded tensor product}, denoted $\op{X} \odot \op{Y}$, by
\[
    ( \op{X} \odot \op{Y} )(m) = \coprod_{i+j=m}
     (\op{X}(i) \otimes \op{Y}(j)) \otimes_{\Sigma_{i} \times \Sigma_{j}}
     I[\Sigma_{m}].
\]
\end{defn}

The proof of the following theorem is straightforward; see, for
example,~\cite{hss:2000}.

\begin{prop}\label{prop:graded-tensor}
Let $\op{U} = \{ \op{U}(n) \}$, where $\op{U}(0) = I$, and
$\op{U}(n)=O$ if $n>0$. Then $(\cat{M}^{\Sigma}, \odot, \op{U})$
is a closed symmetric monoidal category, called the \emph{graded
monoidal structure on} $\cat{M}^{\Sigma}$.
\end{prop}

We call a symmetric sequence $\op{X}$ with an associative
multiplication $\mu : \op{X} \odot \op{X} \rightarrow \op{X}$ a
\emph{graded monoid}.  If furthermore $\mu$ is unital with
respect to a given unit $\nu:\op{U} \rightarrow \op{X}$, then we
call $(\op{X},\mu,\nu)$ a \emph{unital graded monoid}. Dually, we have
\emph{graded comonoids} and \emph{counital graded comonoids}.

\begin{notn}
The symmetric group $\Sigma_{m}$ acts on the left on $I_{m,n}$ via
permutation of entries, namely, $\sigma\cdot\vec{n} :=
(n_{\sigma^{-1}(1)}, \ldots , n_{\sigma^{-1}(m)})$ for $\vec{n} =
(n_{1},\ldots,n_{m}) \in I_{m,n}$.  A permutation $\sigma \in
\Sigma_{m}$ and an $m$-tuple $\vec{n} \in I_{m,n}$ determine the
block permutation $\sigma_{\vec{n}} \in \Sigma_{n}$.

Given a symmetric sequence $\op{X}$, we set $\op{X}[\vec{n}] =
\op{X}(n_{1}) \otimes \cdots \otimes \op{X}(n_{m})$.  The action
of $\sigma \in \Sigma_{m}$ on $I_{m,n}$ defines an isomorphism
$\hat{\sigma} : \op{X}[\vec{n}] \xrightarrow{\cong} \op{X}[\sigma
\vec{n}]$.
\end{notn}

\begin{rmk}
Let $\op{X}$ be a symmetric sequence, and let $m,n \in \N$. We
have the identity
\[
    \op{X}^{\odot m}(n) = \coprod_{m \geq 0}
         \coprod_{\vec{n} \in I_{m,n}}
            \op{X}[\vec{n}] \otimes_{\Sigma_{\vec{n}}}
            I[\Sigma_{n}].
\]
where $\Sigma_{\vec{n}} = \Sigma_{n_{1}} \times \cdots \times
\Sigma_{n_{m}}$. The symmetric group $\Sigma_{m}$ acts on the left
on $\op{X}^{\odot m}(n)$ by $\hat{\sigma}\otimes
\sigma_{\vec{n}}^{-1} : \op{X}[\vec{n}] \otimes I[\Sigma_{n}]
\rightarrow \op{X}[\sigma\vec{n}] \otimes I[\Sigma_{n}]$, while
$\Sigma_{n}$ acts on the right via multiplication.
\end{rmk}

The category of symmetric sequences comes equipped with a third,
non-symmetric, monoidal product, $\boxprod$, called the
\emph{composition product}.

\begin{defn}\label{defn:box}
Let $\op{X}$ and $\op{Y}$ be symmetric sequences.  The \emph{composition product} of $\op X$ and $\op Y$, denoted $\op X\circ \op Y$, is the symmetric sequence with
\[
    (\op{X} \boxprod \op{Y})(n) = \coprod_{m \geq 0} \op{X}(m)
    \otimes_{\Sigma_{m}} (\op{Y}^{\odot m})(n)
\]
for all $n\in \N$.\end{defn}

The following result is proved in~\cite[section
II.1.8]{markl-shnider-stasheff}.

\begin{prop}\label{prop:mon-cat}
Let $\op{J} = \{ \op{J}(n) \}$, where $\op{J}(1) = I$ and
$\op{J}(n) = O$ if $n \neq 1$. Then $(\cat{M}^{\Sigma}, \circ,
\op{J})$ is a monoidal category.
\end{prop}

We will call $(\cat{M}^{\Sigma},\circ,\op{J})$ the
\emph{composition monoidal structure} on $\cat{M}^{\Sigma}$.
Composition monoids, or operads, are discussed in the next
section.

From the definitions, we obtain a useful formula for iterated
composition products in terms of the graded tensor product.

\begin{prop}\label{prop:box-level}
Let $\op{X}$, $\op{Y}$ and $\op{Z}$ be symmetric sequences.  Then
for $n \geq 1$,
\[
    \big(( \op{X} \boxprod \op{Y}) \boxprod \op{Z} \big)(n)
        \cong \coprod_{k \leq m \leq n}
            \left( \op{X}(k) \otimes_{\Sigma_{k}}
            \op{Y}^{\odot k}(m) \right) \otimes_{\Sigma_{m}}
            \op{Z}^{\odot m}(n)
\]
as right $\Sigma_{n}$-objects.
\end{prop}

The natural transformation of the next proposition, which is very easily proved, plays an absolutely crucial role throughout this article.

\begin{prop}\label{prop:level-tensor-box}
Let $\op{X}$, $\op{X}'$, $\op{Y}$, and $\op{Y}'$ be symmetric
sequences.  There is a natural morphism of symmetric sequences
\[
    \iota : ( \op{X} \otimes \op{X}' ) \boxprod ( \op{Y} \otimes \op{Y} ')
    \rightarrow ( \op{X} \boxprod \op{Y} ) \otimes ( \op{X}'
    \boxprod \op{Y}' ).
\]
\end{prop}

We call the morphism $\iota$ the \emph{intertwiner}.

The following proposition seems to be implicit
in the literature.

\begin{prop}\label{prop:extending-to-box}
Let $\op{X}$, $\op{Y}$, and $\op{Z}$ be symmetric sequences.  Let
\[
    \varphi_{\vec{n}} : \op{X}(m) \otimes
    \op{Y}[\vec{n}] \rightarrow \op{Z}(n)
\]
be a family of morphisms, for all $\vec{n} \in I_{m,n}$ and for
all $m,n \geq 0$. Suppose the diagrams
\[
    \xymatrix{
        {\op{X}(m) \otimes \op{Y}[\vec{n}]}
            \ar[r]^{1 \otimes \hat{\sigma}}
            \ar[dd]_{\sigma \otimes Id}
        & {\op{X}(m) \otimes \op{Y}[\sigma\vec{n}]}
            \ar[d]^{\varphi_{\sigma\vec{n}}} \\
        & {\op{Z}(n)}
            \ar[d]^{\sigma_{\vec{n}}}   \\
        {\op{X}(m) \otimes \op{Y}[\vec{n}]}
            \ar[r]_{\varphi_{\vec{n}}}
        & {\op{Z}(n)}
    }
\]
and
\[
    \xymatrix{
        {\op{X}(m) \otimes \op{Y}[\vec{n}]}
            \ar[r]^{\varphi_{\vec{n}}}
            \ar[d]_{1 \otimes \tau_{1} \otimes \cdots \otimes
            \tau_{m}}
        & {\op{Z}(n)}
            \ar[d]^{\tau_{1} \oplus \cdots \oplus \tau_{m}} \\
        {\op{X}(m) \otimes \op{Y}[\vec{n}]}
            \ar[r]_{\varphi_{\vec{n}}}
        & {\op{Z}(n)}
    }
\]
commute for all $\vec{n} \in I_{n,m}$, $\sigma \in \Sigma_{m}$,
and $\tau_{j} \in \Sigma_{n_{j}}$. Then the $\varphi_{\vec{n}}$
extend uniquely to a morphism of symmetric sequences,
\[
    \varphi : \op{X} \boxprod \op{Y} \rightarrow \op{Z}.
\]
\end{prop}

\begin{proof}
The $\Sigma_{n_{1}} \times \cdots \times
\Sigma_{n_{m}}$-equivariance means that for a fixed $n$, we may
extend
\[
    \varphi_{n} := \coprod_{\vec{n} \in I_{m,n}} \varphi_{\vec{n}}
\]
to
\[
    \varphi_{m,n} : \op{X}(m) \otimes \op{Y}^{\odot m}(n) \rightarrow
    \op{Z}(n).
\]
We may then pass to orbits because of $\Sigma_{m}$-equivariance.
\end{proof}

\begin{cor}\label{cor:extend-to-box}
Let $k$, $m_{1}, \ldots, m_{k}$ be natural numbers, and for $i =
1, \ldots, k$, let $\vec{n}_{i} = (n_{i1}, \ldots, n_{im_{i}}) \in
I_{m_{i},n_{i}}$.  Let $n = n_{1} + \cdots + n_{k}$. If a family
of morphisms
\[
    \op{W}(k) \otimes \op{X}(m_{1}) \otimes \cdots \otimes
    \op{X}(m_{k}) \otimes \op{Y}[\vec{n}_{1}] \otimes \cdots
    \otimes
    \op{Y}[\vec{n}_{k}] \rightarrow \op{Z}(n)
\]
is equivariant with respect to the actions of $\Sigma_{k}$,
$\Sigma_{m_{1}} \times \cdots \times \Sigma_{m_{k}}$, and
$\Sigma_{n_{11}} \times \cdots \times \Sigma_{n_{km_{k}}}$, then
the morphisms determine a unique morphism of symmetric sequences,
\[
    (\op{W} \boxprod \op{X}) \boxprod \op{Y} \rightarrow \op{Z}.
\]
\end{cor}

\subsection{Operads and their modules}\label{ssec:bimodules}

An \emph{operad} is a unital monoid in the category of symmetric
sequences, with respect to the composition product.  The structure
morphism in the operad $\op{P}$,
\[
    \gamma : \op{P} \boxprod \op{P} \rightarrow \op{P}
\]
is called the \emph{composition product}, while $\eta : \op{J}
\rightarrow \op{P}$ is the \emph{unit}.

A \emph{left $\op{P}$-module} is a symmetric sequence $\op{M}$
equipped with an associative, unital structure morphism in
$\cat{M}^{\Sigma}$,
\[
    \lambda : \op{P} \boxprod \op{M} \rightarrow \op{M}.
\]
A morphism of left $\op{P}$-modules is a morphism in
$\cat{M}^{\Sigma}$ that commutes with the structure morphisms. We
denote by ${}_{\op{P}}\cat{Mod}$ the category of left
$\op{P}$-modules and $\op{P}$-linear morphisms.

Recall that a $\op{P}$-\emph{algebra} is an object $A$ in $\cat{M}$
equipped with a sequence of structure morphisms
$\varphi_{m}:\op{P}(m) \otimes A^{\otimes m} \rightarrow A$ that
are associative, equivariant, and unital. A morphism of
$\op{P}$-algebras is a morphism in $\cat{M}$, $f : A \rightarrow
B$, that commutes with the structure morphisms of $A$ and $B$.
The category of $\op{P}$-algebras is denoted $\alg{P}$.

\begin{ex}\label{ex:assoc}
Recall the symmetric sequence $\op{A}$ defined in
Example~\ref{ex:ass-as-level-comonoid}.  As is well known, the
permutation of blocks composed with permutations within blocks
defines a composition product on $\op{A}$ that make it into an
operad, called the \emph{associative operad}.  The
$\op{A}$-(co)algebras are precisely the (co)monoids in
$(\cat{M}, \otimes, I)$.

We may form the associative operad in $(\cat{M}^{\Sigma}, \otimes
, \op{C})$ to obtain the operad $\op{A}_{\otimes}$ whose
(co)algebras are precisely the level (co)monoids.  Explicitly,
$\op{A}_{\otimes}(n)$ is the symmetric sequence
$\op{C}[\Sigma_{n}]$. Similarly, forming the associative operad in
$(\cat{M}^{\Sigma} , \odot , \op{U} )$, we obtain the operad
$\op{A}_{\odot}$, whose (co)algebras are the graded (co)monoids.
\end{ex}

There are at least two ways to embed $\op{P}$-algebras as left
$\op{P}$-modules, which we outline below.

Let $c(A)$ be the constant symmetric sequence defined by $c(A)(0)
= O$, $c(A)(n) = A$ for $n \geq 1$, with trivial
$\Sigma_{n}$-action for all $n \geq 0$. Then $c(A)$ is a left
$\op{P}$-module. Indeed, for $\vec{n} \in I_{m,n}$, we have
$c(A)[\vec{n}] = A^{\otimes m}$. Define $\varphi_{\vec{n}}$ as the
composite
\[
    \op{P}(m) \otimes c(A)[\vec{n}] \cong \op{P}(m) \otimes
    A^{\otimes m} \xrightarrow{\varphi_{m}} A = c(A)(n).
\]
One verifies readily that $\varphi_{\vec{n}}$ is equivariant with
respect to the actions of $\Sigma_{m}$ and $\Sigma_{\vec{n}}$.  By Proposition~\ref{prop:extending-to-box}, the family $\{\varphi_{\vec n}\}_{\vec n}$ gives rise to a
morphism of symmetric sequences
\[
    \op{P} \boxprod c(A) \rightarrow c(A).
\]

\begin{prop}\label{prop:c}
The functor $c : \alg{P} \rightarrow {}_{\op{P}}\cat{Mod}$ is
faithful but not necessarily full.
\end{prop}

\begin{proof}
The definition of $c$ on objects is outlined above.  If $f : A
\rightarrow B$ is a morphism of $\op{P}$-algebras, then we set
$c(f)_{n} = f : A \rightarrow B$. Clearly $c$ commutes with the
composition of morphisms.  The fact that $c(f)$ is a morphism of
left $\op{P}$-modules follows directly from the fact that $f$ is a
morphism of $\op{P}$-algebras. Furthermore, if $c(f) = c(g)$, then
evaluation at any level reveals that $f=g$, so $c$ is faithful.

To see that $c$ need not be full, suppose that $\cat{M}$ is
additive, and consider any $\op{P}$-algebra $A$. Then $F : c(A)
\rightarrow c(A)$, $F_{n} = (-1)^{n}Id_{A}$, $n \geq 0$, is a
morphism of left $\op{P}$-modules that is not of the form $c(g)$
for any $g : A \rightarrow A$.
\end{proof}

To a symmetric sequence $\op{X}$ we associate an endofunctor
$\op{X}(-):\cat{M} \rightarrow \cat{M}$, defined for $V \in
\cat{M}$ by
\[
    \op{X}(V) = \coprod_{m \geq 0} \op{X}(m) \otimes_{\Sigma_{m}}
    V^{\otimes m}.
\]
The functor $\op{X} \mapsto \op{X}(-)$ is strongly monoidal, that
is, $\op{X}(\op{Y}(V)) \cong (\op{X}\circ\op{Y})(V)$ for all $V
\in \ob\cat{M}$.  As a consequence, if $\op{P}$ is an operad, then
$\op{P}(-)$ defines a monad.

Let $X \in \ob \cat{M}$.  Define $z(X)$ to be the symmetric sequence
with $z(X)(0) = X$ and $z(X)(n) = O$ for $n > 0$. If $f : X
\rightarrow Y$ is a morphism in $\cat{M}$, set $z(f)_{0} = f$,
while $z(f)_{n}$ is necessarily the identity on the initial object $O$ for $n
> 0$.  Clearly $z$ defines a functor, $\cat{M} \rightarrow
\cat{M}^{\Sigma}$. If $F : z(X) \rightarrow z(Y)$ is a morphism of
symmetric sequences, then $F_{0} : X \rightarrow Y$ is a morphism
in $\cat{M}$, and $F_{n}$ is the identity on $O$.  Thus $F =
z(F_{0})$ and so $z$ is full. Furthermore, if $f,g : X \rightarrow
Y$ and $z(f) = z(g)$, then $z(f)_{0} = z(g)_{0}$, so $f=g$ and $z$
is faithful.

\begin{prop}\label{prop:z}
The functor $z$ restricts to define a full and faithful functor $z
: \alg{P} \rightarrow {}_{\op{P}}\cat{Mod}$.
\end{prop}

\begin{proof}
The structure morphisms on $A$ can be expressed as one morphism
$\varphi : \op{P}(A) \rightarrow A$.  It is an easy exercise to
show that $\op{P} \circ z(A) \cong z(\op{P}(A))$. Therefore $z(A)$
is a left $\op{P}$-module, with structure morphism $\op{P} \circ
z(A) \cong z(\op{P}(A)) \xrightarrow{z(\varphi)} z(A)$.  Since the
isomorphism $\op{P} \circ z(A) \cong z(\op{P}(A))$ is natural, $z$
turns morphisms of $\op{P}$-algebras into morphisms of left
$\op{P}$-modules.  The restriction of a faithful functor is
evidently faithful.

Suppose that $F : z(A) \rightarrow z(B)$ is a morphism of left
$\op{P}$-modules. We need to show that $F_{0} : A \rightarrow B$
is a morphism of $\op{P}$-algebras. This follows from the diagram
expressing the fact that $F$ is a morphism of left
$\op{P}$-modules, in level zero.
\end{proof}

Propositions~\ref{prop:c} and~\ref{prop:z} show that left
$\op{P}$-modules are a generalization of $\op{P}$-algebras.  In
fact, one may consider left $\op{P}$-modules to be \emph{graded}
$\op{P}$-algebras, in the sense that the structure morphisms
define maps
\[
    \op{P}(n) \otimes_{\Sigma_{n}} \op{M}^{\odot n} \rightarrow
    \op{M}.
\]

\begin{rmk}
There exists a functor, $u : \cat{M}^{\Sigma} \rightarrow
(\cat{M}^{\Sigma})^{\Sigma}$, where $u(\op{X})_{n}$ is the
symmetric sequence $z(\op{X}(n))$.  The functor $u$ is strongly monoidal
with respect to the composition product, and hence $u(\op{P})$ is
an operad for any operad $\op{P}$.  One verifies that
$u(\op{P})(\op{X}) = \op{P}\circ\op{X}$; as a result, a left
$\op{P}$-module is exactly a $u(\op{P})$-algebra.
\end{rmk}

A \emph{right $\op{P}$-module} is a symmetric sequence $\op{N}$
along with an associative, unital structure morphism in
$\cat{M}^{\Sigma}$,
\[
    \rho : \op{N} \boxprod \op{P} \rightarrow \op{N}.
\]
A morphism of right $\op{P}$-modules is a morphism in
$\cat{M}^{\Sigma}$ that commutes with the structure morphisms. We
denote by $\cat{Mod}_{\op{P}}$ the category of right
$\op{P}$-modules and $\op{P}$-linear morphisms.

Recall that a $\op{P}$-coalgebra is an object $A$ of $\cat{M}$
along with a sequence of structure morphisms
\[
    \theta_{n} : A \otimes \op{P}(n) \rightarrow A^{\otimes n}
\]
that are associative, equivariant, and unital.  The category of
$\op{P}$-coalgebras and morphisms respecting the structure maps is
denoted $\coalg{P}$.

Let $\op{Q}$ be another operad.  A
\emph{$(\op{P},\op{Q})$-bimodule} is a symmetric sequence $\op{M}$
that is simultaneously a left $\op{P}$-module and a right
$\op{Q}$-module, in such a way that the diagram
\[
    \xymatrix{
        {\op{P} \boxprod \op{M} \boxprod \op{Q}}
            \ar[r]^-{Id \boxprod \rho}
            \ar[d]_{\lambda \boxprod Id}
        & {\op{P} \boxprod \op{M}}
            \ar[d]^{\lambda} \\
        {\op{M} \boxprod \op{Q}}
            \ar[r]_{\rho}
        & {\op{M}}
    }
\]
commutes.  A morphism of $(\op{P},\op{Q})$-bimodules is a morphism
in $\cat{M}^{\Sigma}$ that commutes with both structure morphisms.
We denote by ${}_{\op{P}}\cat{Mod}_{\op{Q}}$ the category of
$(\op{P},\op{Q})$-bimodules and their morphisms.

Let $\op{P}$ be an operad in $\cat{M}$. Recall that in
$\cat{M}^{\Sigma}$, colimits are calculated level-wise. Since
$\cat{M}$ is assumed to be cocomplete, coequalizers exist in
$\cat{M}^{\Sigma}$.  Let $\op{M} \in \ob\cat{Mod}_{\op{P}}$, $\op N \in
{}_{\op{P}}\cat{Mod}$. Define the symmetric sequence $\op{M}
\underset {\op P}\circ  \op{N}$ to be the coequalizer
\[
    \xymatrix{
        {\op{M} \boxprod \op{P} \boxprod \op{N}}
            \ar@<1ex>[r]^-{Id \boxprod \lambda}
            \ar@<-1ex>[r]_-{\rho \boxprod Id}
        & {\op{M} \boxprod \op{N}} \ar[r]
        & {\op{M} \underset {\op P}\circ  \op{N}}.
    }
\]

Composition over $\op{P}$ behaves well with respect to actions of
other operads.

\begin{prop}\label{prop:qo-module-structure}
Let $\op{O}$ and $\op{Q}$ be operads.  If $\op{M} \in
{}_{\op{Q}}\cat{Mod}_{\op{P}}$ and $\op{N} \in
{}_{\op{P}}\cat{Mod}_{\op{O}}$, then $\op{M} \underset {\op P}\circ
\op{N}$ admits a natural $(\op{Q},\op{O})$-bimodule structure.
\end{prop}

Observe that $({}_{\op{P}}\cat{Mod}_{\op{P}},\underset {\op P}\circ , \op P)$ is a monoidal category.

Let $A$ be an object of $\cat{M}$.  Define the \emph{free
graded monoid} on $A$ to be the symmetric sequence $\op{T}(A)$
with $\op{T}(A)(n) = A^{\otimes n}$ and $\op{T}(A)(0) = I$, where
$\Sigma_{n}$ acts by permuting factors. In fact, this construction defines a functor, $\op{T} : \cat{M}
\rightarrow \cat{M}^{\Sigma}$, where $\op{T}(f)_{n} = f^{\otimes
n}$ for a given morphism $f : X \rightarrow Y$.

\begin{prop}\label{prop:right-module}
    \begin{enumerate}
    \item The functor $\op{T} : \cat{M} \rightarrow
    \cat{M}^{\Sigma}$ is faithful but not full.
    \item $\op{T}$ takes values in ${}_{\op{A}}\cat{Mod}$, and is full as
    a functor $\cat{M} \rightarrow {}_{\op{A}}\cat{Mod}$.
    \item Let $\op{P}$ be an operad. An object $A$ in $\cat{M}$ is a
    $\op{P}$-coalgebra if and only if $\op{T}(A)$ is an
    $(\op{A},\op{P})$-bimodule.
    \item
    A morphism $f : A \rightarrow B$ is a morphism of $\op{P}$-coalgebras
    if and only if $\op{T}(f)$ is a morphism of $(\op{A},\op{P})$-bimodules.
    \end{enumerate}
    Thus $\op{T}$ restricts to define a full and faithful functor $\coalg{P}
    \rightarrow {}_{\op{A}}\cat{Mod}_{\op{P}}$.
\end{prop}

\begin{proof}
(1)  Let $f,g:X \rightarrow Y$ be morphisms in $\cat{M}$ such that
$\op{T}(f) = \op{T}(g)$. Thus $\op{T}(f)_{1} = \op{T}(g)_{1}$,
i.e., $f = g$.  Thus $\op{T}$ is faithful.  To see that $\op{T}$ is
not full in general, we suppose that $\cat{M}$ is additive.  Let $X \in
\cat{M}$, and set $F_{n} = -(-1)^{n}Id_{X^{\otimes n}} : X^{\otimes
n} \rightarrow X^{\otimes n}$.  Thus $F : \op{T}(X) \rightarrow
\op{T}(X)$ is a morphism of symmetric sequences, but $F \neq
\op{T}(F_{1})$, since $F_{1} = Id_{X}$.

\medskip
(2) Let $X \in \ob \cat{M}$.  For $n,m \geq 0$ and for each partition
$(n_{1}, \ldots , n_{m}) \in I_{m,n}$, there is a natural
isomorphism $X^{\otimes n_{1}} \otimes \cdots \otimes X^{\otimes
n_{m}} \rightarrow X^{\otimes n}$, where $n = n_{1} + \cdots +
n_{m}$, that is $\Sigma_{n_{1}} \times \cdots \times
\Sigma_{n_{m}}$-equivariant. These isomorphisms add up to define a
morphism $\op{A}(m) \otimes_{\Sigma_{m}}\op{T}(X)^{\odot{m}}(n)
\rightarrow \op{T}(X)(n)$.  It is routine to verify that these
morphisms together form an associative, unital action $\op{A}
\circ \op{T}(X) \rightarrow \op{T}(X)$, such that if $f : X
\rightarrow Y$ is a morphism in $\cat{M}$, then $\op{T}(f)$ is a
morphism of left $\op{A}$-modules.

Let $F : \op{T}(X) \rightarrow \op{T}(Y)$ be a morphism of left
$\op{A}$-modules.  For all $n \geq 1$, the diagram
\[
    \xymatrix{
        (\op{T}(X)(1))^{\otimes n}
            \ar[r]
            \ar[d]_{F_{1}^{\otimes n}}
        & \op{T}(X)(n)
            \ar[d]^{F_{n}}
            \\
        (\op{T}(Y)(1))^{\otimes n}
            \ar[r]
        & \op{T}(Y)(n)
    }
\]
commutes, where the horizontal arrows come from the left
$\op{A}$-module structure morphisms; in this case, they are the
identity morphisms. Thus $F_{n} = F_{1}^{\otimes n}$, so $F =
\op{T}(F_{1})$, and $\op{T} : \cat{M} \rightarrow
{}_{\op{A}}\cat{Mod}$ is full.

\medskip
(3) Suppose that $\op{T}(A)$ is a right $\op{P}$-module with structure
morphism
\[
    \rho : \op{T}(A) \boxprod \op{P} \rightarrow \op{T}(A).
\]
Recall that
\[
    (\op{T}(A) \boxprod \op{P} )(n) = \coprod_{0\leq m \leq n}
    A^{\otimes m} \underset{\Sigma_{m}}{\otimes} \op{P}^{\odot m}(n).
\]
Let $\theta_{n}$ be the component of $\rho_{n}$ with $m=1$:
\[
    \theta_{n} : A \otimes
    \op{P}(n) \rightarrow A^{\otimes n}.
\]
Then $\theta_{n}$ inherits associativity, equivariance and unity
from $\rho$, i.e., $A$ is a $\op P$-coalgebra.

Conversely, suppose that $A$ is a $\op{P}$-coalgebra.  The
structure morphisms $\theta_{n}$ define morphisms in $\cat{M}$,
for $m \geq 1$, and $\vec{n} = (n_{1}, \ldots, n_{m})$,
\[
    A^{\otimes m} \otimes \op{P}[\vec{n}]
    \xrightarrow{\text{shuffle}}
    \bigotimes_{j=1}^{m}
        \left(
            A \otimes \op{P}(n_{j})
        \right)
    \xrightarrow{\theta_{n_{1}} \otimes \cdots \otimes
    \theta_{n_{m}}}
        A^{\otimes n},
\]
where $n = n_{1} + \cdots + n_{m}$.  These morphisms are
equivariant with respect to the actions of $\Sigma_{m}$ and
$\Sigma_{n_{1}} \times \cdots \times \Sigma_{n_{m}}$, and hence
piece together to form a morphism of left $\op A$-modules
\[
    \op{T}(A) \boxprod \op{P} \rightarrow \op{T}(A)
\]
that is associative and unital because of the corresponding
properties of the $\theta_{n}$, i.e., $\op T(A)$ is a right $\op P$-module.

\medskip
(4) Suppose that $f : A \rightarrow B$ is a morphism of
$\op{P}$-coalgebras.  We have already
remarked that $\op{T}(f)$ is a morphism of left $\op{A}$-modules.  Let
$m,n \geq 0$, $\vec{n} \in I_{m,n}$.  In the diagram
\[
    \xymatrix{
        A^{\otimes m} \otimes \op{P}[\vec{n}]
            \ar[rr]^{\rho}
            \ar[dr]_{\text{permute}}
            \ar[ddd]_{f^{\otimes m} \otimes Id}
        & & A^{\otimes n}
            \ar[ddd]^{f^{\otimes n}} \\
        & A \otimes \op{P}(n_{1}) \otimes \cdots \otimes A \otimes
            \op{P}(n_{m})
                \ar[d]^{(f \otimes Id)^{\otimes m}}
                \ar[ur]^{\theta_{\vec{n}}} \\
        & B \otimes \op{P}(n_{1}) \otimes \cdots \otimes B \otimes
            \op{P}(n_{m})
                \ar[dr]^{\theta_{\vec{n}}}
                \\
        B^{\otimes m} \otimes \op{P}[\vec{n}]
            \ar[rr]_{\rho}
            \ar[ur]_{\text{permute}}
        & & B^{\otimes n}
    }
\]
the top and bottom cells commute by definition, the left cell
commutes by naturality of the permutation, and the right cell
commutes because $f$ is a morphism of $\op{P}$-coalgebras.
Therefore the outer diagram commutes, and so $\op{T}(f)$ is a
morphism of right $\op{P}$-modules, i.e., $\op T(f)$ is a morphism of $(\op A, \op P)$-bimodules.

Finally, let $F : \op{T}(A) \rightarrow \op{T}(B)$ be a morphism
of $(\op{A},\op{P})$-bimodules.  Since $F$ is a morphism of left
$\op{A}$-modules it is of the form $F = \op{T}(f)$, where $f =
F_{1} :A \rightarrow B$ in $\cat{M}$. Since $F$ is a morphism of
right $\op{P}$-modules, the diagram
\[
    \xymatrix{
        A \otimes \op{P}(n)
            \ar[r]^{\rho}
            \ar[d]_{F_{1} \otimes Id}
        & A^{\otimes n}
            \ar[d]^{F_{n}} \\
        B \otimes \op{P}(n)
            \ar[r]_{\rho}
        & B^{\otimes n}
    }
\]
commutes.  Since $F_{1} = f$ and $F_{n} = f^{\otimes n}$, we have
shown that $f$ is a morphism of $\op{P}$-coalgebras.
\end{proof}

\begin{ex}\label{ex:hom-rel-op}
The template for bimodules is the \emph{homomorphism bimodule}.
Since our symmetric monoidal category $\cat{M}$ is closed, it has
an internal hom functor $\Hom (X,-)$ that is right
adjoint to $- \otimes X$ for all objects $X$ of $\cat{M}$.  For all $X,Y\in \ob \cat M$, let $ev_{X}:\Hom (X,Y)\otimes X\to Y$ denote the adjoint of the identity $\Hom (X,Y)\xrightarrow =\Hom (X,Y)$.  For all $X,X',Y,Y'\in \ob \cat M$, let
$$t:\Hom (X,Y)\otimes \Hom (X',Y')\to \Hom (X\otimes X', Y\otimes Y')$$
denote the adjoint of the composite
$$ \Hom (X,Y)\otimes \Hom (X',Y')\otimes X\otimes X'\xrightarrow{\cong}\Hom (X,Y)\otimes X \otimes \Hom (X',Y')\otimes X'\xrightarrow{ev_{X}\otimes ev_{X'}}Y\otimes Y'.$$
Finally, for all $X,Y,Z\in \ob \cat M$, let $c:\Hom (Y,Z)\otimes \Hom (X,Y)\to \Hom (X,Z)$ denote the adjoint of the composite
$$\Hom (Y,Z)\otimes \Hom (X,Y)\otimes X\xrightarrow {Id\otimes ev_{X}}\Hom (Y,Z)\otimes Y\xrightarrow{ev_{Y} }Z.$$

Let $X,Y\in \ob\cat{M}$. Recall that the
\emph{endomorphism operad} is defined by $\End_{X}(n) =
\Hom (X^{\otimes n},X)$. The composition product
\[
    \End_{X}(n) \otimes \End_{X}(m_{1}) \otimes \cdots \otimes
    \End_{X}(m_{n})
     \rightarrow \End_{X}(m)
\]
is defined to be the composite
$$\xymatrix{
\End_{X}(n)\otimes\End_{X}(m_{1}) \otimes \cdots \otimes
    \End_{X}(m_{n})\ar[rr]^(0.55){Id\otimes t^{n-1}}&&\Hom (X^{\otimes n}, X)\otimes \Hom (X^{\otimes m}, X^{\otimes n})\ar[d]^{c}\\
    && \Hom (X^{\otimes m}, X)=\End _{X}(m).}$$

Let $\operatorname{Hom}_{X,Y}(m) = \Hom (X^{\otimes m} , Y)$
for $m \geq 1$.  We set $\operatorname{Hom}_{X , Y}(0) = \Hom(I,Y)$.
Define a left $\End_{Y}$-module structure via the composite
$$\xymatrix{
    \End_{Y}(n) \otimes \operatorname{Hom}_{X,Y}(m_{1}) \otimes \cdots \otimes
    \operatorname{Hom}_{X,Y}(m_{n}) \ar[rr]^(0.55){Id\otimes t^{n-1}}&&\Hom(Y^{\otimes n},Y)\otimes \Hom (X^{\otimes m}, Y^{\otimes n})\ar[d]^{p}\\
   &&\operatorname{Hom}_{X,Y}(m).}$$

Similarly, one defines a right $\End_{X}$-module structure.  Since
precomposition commutes with postcomposition, these two structures
are compatible, so that $\operatorname{Hom}_{X,Y}$ is an
$(\End_{Y},\End_{X})$-bimodule.

Similarly, if we set $\CoHom_{X,Y}(n) =
\Hom(X,Y^{\otimes n})$, $\Coend_{X} = \CoHom_{X,X}$, and
$\Coend_{Y} = \CoHom_{Y,Y}$, then $\CoHom_{X,Y}$ is a
$(\Coend_{X},\Coend_{Y})$-bimodule.
\end{ex}

The proof of the following proposition works exactly as in the
classical case of modules over associative rings.

\begin{prop}\label{prop:morphisms}
Let $\theta : \op{Q} \rightarrow \op{Q}'$ and $\varphi : \op{P}
\rightarrow \op{P}'$ be morphisms of operads.  Then $\theta$ and
$\varphi$ induce a functor
\[
    {}_{\op{Q}'}\cat{Mod}_{\op{P}'} \rightarrow
        {}_{\op{Q}}\cat{Mod}_{\op{P}}.
\]
\end{prop}

\begin{ex}
Let $\op{P}$ be an operad.  A $\op{P}$-algebra structure on $X$ is
equivalent, via adjunction, to a morphism of operads, $\op{P}
\rightarrow \End_{X}$.  Let $\op{Q}$ be another operad, let $Y$ be
a $\op{Q}$-algebra with structure morphism $\op{Q} \rightarrow
\End_{Y}$, and let $\op{M}$ be a $(\op{Q},\op{P})$-bimodule.
 By Proposition~\ref{prop:morphisms},
$\operatorname{Hom}_{X,Y}$ is a $(\op{Q},\op{P})$-bimodule.  A morphism of
$(\op{Q},\op{P})$-bimodules, $\op{M} \rightarrow \operatorname{Hom}_{X,Y}$,
corresponds via adjunction to a sequence of morphisms $\op{M}(n)
\otimes_{\Sigma_{n}} X^{\otimes n} \rightarrow Y$, or,
equivalently, a morphism $\op{M}(X) \rightarrow Y$.  We will see
in section~\ref{ssec:comonoidal} how this idea allows us to
generalize the notion of ``algebra morphism''.
\end{ex}

For the next proposition, we remark that the intertwiner $\iota$
of Proposition~\ref{prop:level-tensor-box} yields monoidal
categories $(\cat{Comon}_{\otimes}, \circ, \op{J})$ and
$(\cat{Op}, \otimes, \op{C})$.  Namely, if $(\op{X},\Delta _{\op X})$ and $(\op{Y},\Delta _{\op Y})$
are level comonoids, then the composite
\[
    \op{X} \circ \op{Y} \xrightarrow{\Delta_{\op X}\circ \Delta _{\op Y}} (\op{X} \otimes \op{X}) \circ
    (\op{Y} \otimes \op{Y}) \xrightarrow{\iota} (\op{X} \circ
    \op{Y}) \otimes (\op{X} \circ \op{Y})
\]
defines a coassociative coproduct on $\op X\circ \op Y$, thanks to the coassociativity of $\Delta _{\op X}$ and of $\Delta _{\op Y}$ and to the naturality of the intertwiner.  Furthermore, if $\varepsilon_{\op X}:\op X\to \op C$ and $\varepsilon _{\op Y}:\op Y\to \op C$ are counits, then
\[
    \op{X} \circ \op{Y} \xrightarrow {\varepsilon_{\op X}\circ \varepsilon _{\op Y}}\op{C} \circ \op{C}
    \xrightarrow{\gamma} \op{C}
\]
is a counit for $\op{X} \circ \op{Y}$, where $\gamma$ is the composition multiplication of the operad $\op C$.
Furthermore, the isomorphism $\op{J} \cong \op{J} \otimes \op{J}$
determines a level comonoid structure in $\op{J}$.  Similarly, if
$\op{P}$ and $\op{Q}$ are operads, then $\op{P} \otimes \op{Q}$ is
an operad with composition product
\[
    (\op{P} \otimes \op{Q}) \circ (\op{P} \otimes \op{Q})
    \xrightarrow{\iota} (\op{P} \circ \op{P}) \otimes (\op{Q}
    \circ \op{Q}) \xrightarrow{\gamma \otimes \gamma} \op{P}
    \otimes \op{Q}
\]
and unit $\op{J} \cong \op{J} \otimes \op{J} \rightarrow \op{P}
\otimes \op{Q}$.

Recall that a \emph{Hopf operad} is  a counital comonoid $\op P$ in the
monoidal category $(\cat{Op}, \otimes , \op{C})$ or, equivalently, a unital
monoid in $(\cat{Comon}_{\otimes}, \circ, \op{J})$.

\begin{prop}\label{prop:tensor}
If $\op{P}$ and $\op{Q}$ are Hopf operads, then
$({}_{\op{P}}\cat{Mod}_{\op{Q}},\otimes , \op{C})$ is a
monoidal category.
\end{prop}

\begin{proof}
By definition, the level diagonals on $\op{P}$ and $\op{Q}$
are morphisms of operads, and so by
Proposition~\ref{prop:morphisms} define a functor ${}_{\op{P}
\otimes \op{P}}\cat{Mod}_{\op{Q} \otimes \op{Q}} \rightarrow
{}_{\op{P}}\cat{Mod}_{\op{Q}}$. Let ${}_{\op{P}}\cat{Mod}_{\op{Q}}
\times {}_{\op{P}}\cat{Mod}_{\op{Q}} \rightarrow {}_{\op{P}
\otimes \op{P}}\cat{Mod}_{\op{Q} \otimes \op{Q}}$ be the functor
defined on objects by $(\op{M},\op{N}) \mapsto \op{M} \otimes
\op{N}$. The composite of the above two functors defines a
monoidal product in ${}_{\op{P}}\cat{Mod}_{\op{Q}}$ that is
associative, since the level diagonals in $\op{P}$ and $\op{Q}$
are coassociative.

The counits in $\op{P}$ and $\op{Q}$ are morphisms of operads
$\op{P} \rightarrow \op{C}$ and $\op{Q} \rightarrow \op{C}$ that
define a $(\op{P},\op{Q})$-bimodule structure in $\op{C}$. Let
$\op{M} \in {}_{\op{P}}\cat{Mod}_{\op{Q}}$. A diagram chase using
the counit condition shows that the natural isomorphisms  $\op{M}
\otimes \op{C} \xrightarrow{\cong} \op{M}$ and $\op{C}
\otimes \op{M} \xrightarrow{\cong} \op{M}$ are morphisms of
$(\op{P},\op{Q})$-bimodules.
\end{proof}

\begin{cor}
If $\op{P}$ is a Hopf operad, then
$({}_{\op{P}}\cat{Mod},\otimes,\op{C})$ and
$(\cat{Mod}_{\op{P}},\otimes, \op{C})$ are monoidal
categories.
\end{cor}

\begin{proof}
One observes that $\op{J}$, the unit with respect to the
composition product, is a Hopf operad, ${}_{\op{P}}\cat{Mod} =
{}_{\op{P}}\cat{Mod}_{\op{J}}$, and $\cat{Mod}_{\op{P}} =
{}_{\op{J}}\cat{Mod}_{\op{P}}$.
\end{proof}

\subsection{$\op{P}$-co-rings and their applications}\label{ssec:comonoidal}

In this section we describe an operad-theoretic way to de-couple
the objects and morphisms in categories of (co)algebras over an
operad.  Our method is to use morphisms relative to a
$\op{P}$-co-ring $\op{R}$. This simplifies the established
technique of using coloured operads~\cite{markl:04}, and makes the
composition of morphisms an immediate consequence of structure of
the co-ring.  Note that the constructions in this section work in
any monoidal category with coequalizers.

Recall first from Proposition \ref {prop:qo-module-structure} and its immediate consequences that $({}_{\op{P}}\cat{Mod}_{\op{P}}, \underset {\op P}\circ ,\op{P})$ is a monoidal category.

\begin{defn}\label{def:coring}
Let $\op{P}$ be an operad.  A $\op{P}$-\emph{co-ring} is a pair
$(\op{R},\psi)$, where $\op{R}$ is a $\op{P}$-bimodule and $\psi :
\op{R} \rightarrow \op{R} \underset {\op P}\circ  \op{R}$ is a
coassociative morphism of $\op{P}$-bimodules.  A \emph{counital
$\op{P}$-co-ring} is a triple $(\op{R}, \psi, \varepsilon)$, where
$(\op{R}, \psi)$ is a $\op{P}$ co-ring, and $\varepsilon : \op{R}
\rightarrow \op{P}$ is a morphism of $\op{P}$-bimodules with
respect to which $\psi$ is counital.

A morphism of (counital) $\op{P}$-co-rings is a morphism of
$\op{P}$-bimodules that commutes with the structure morphisms. We
denote by $\cat{CoRing}_{\op P}$ the category of
$\op{P}$-co-rings, and by $\cat{CoRing}_{\op P,*}$ the
category of counital $\op{P}$-co-rings.
\end{defn}

The ``co-ring'' terminology is in analogy to that used in ring
theory; if $R$ is an associative, unital ring, then an $R$-co-ring
is a comonoid in $({}_{R}\cat{Mod}_{R}, \otimes_{R}, R)$. In
section~\ref{sec:delooping}, we construct a large family of
co-rings.

Let $\op{P}$ be an operad, and let $(\op{R}, \psi)$ be a
$\op{P}$-co-ring. We define a semicategory (that is, a category
without identity arrows), $\scat{Mod}_{(\op{P},\op R)}$, as follows.
The objects are all right $\op{P}$-modules.  We let
\[
    \scat{Mod}_{(\op{P},\op R)}(\op{M},\op{N}) =
        \cat{Mod}_{\op{P}}(\op{M} \underset {\op P}\circ
        \op{R},\op{N})
\]
for right $\op{P}$-modules $\op{M}$ and $\op{N}$. Let $\theta :
\op{L} \underset {\op P}\circ  \op{R} \rightarrow \op{M}$ and $\varphi :
\op{M} \underset {\op P}\circ  \op{R} \rightarrow \op{N}$ be morphisms. We
define their composition by
\[
    \op{L} \underset {\op P}\circ  \op{R} \xrightarrow{Id\underset {\op P} \circ \psi}
        \op{L} \underset {\op P}\circ  \op{R} \underset {\op P}\circ  \op{R}
        \xrightarrow{\theta \underset {\op P}\circ Id}
        \op{M} \underset {\op P}\circ  \op{R}
        \xrightarrow{\varphi}
        \op{N}.
\]
An elementary diagram chase using the fact that $\psi : \op{R}
\rightarrow \op{R} \underset {\op P}\circ  \op{R}$ is coassociative shows
that the composition so defined is associative, and so
$\scat{Mod}_{(\op{P},\op R)}$ is indeed a semicategory.

Suppose that $(\op{R}, \psi, \varepsilon)$ is a counital
$\op{P}$-co-ring.  Then $\varepsilon$ is used to define the
identity morphisms needed to make $\scat{Mod}_{(\op{P},\op{R})}$
into a full-fledged category. Indeed, if $\op{M} \in
\cat{Mod}_{\op{P}}$, then $Id\underset {\op P}\circ  \varepsilon : \op{M}
\underset {\op P}\circ  \op{R} \rightarrow \op{M} \underset {\op P}\circ  \op{P}
\cong \op{M}$ is the identity morphism of
$\cat{Mod}_{(\op{P},\op{R})}(\op{M}, \op{M})$.  The fact that it
is neutral with respect to composition as defined above follows
from the fact that $\psi$ is counital with respect to
$\varepsilon$.

Note that a morphism of $\op{P}$-co-rings, $\op{Q} \rightarrow \op{R}$,
defines a functor $\scat{Mod}_{(\op{P},\op{R})} \rightarrow
\scat{Mod}_{(\op{P},\op{Q})}$. Note furthermore that $\op{P}$ may be considered to be a counital
$\op{P}$-co-ring in the obvious manner, so that if $(\op{R},\psi,
\varepsilon)$ is a counital $\op{P}$-co-ring, the counit
$\varepsilon : \op{R} \rightarrow \op{P}$ is a morphism of
counital $\op{P}$-co-rings.  Thus the counit defines a functor
$\cat{Mod}_{\op{P}} \rightarrow \cat{Mod}_{(\op{P},\op{R})}$,
where we have identified $\cat{Mod}_{\op{P}}$ and
$\cat{Mod}_{(\op{P},\op{P})}$.

Let $(\op{R},\psi)$ be a $\op{P}$-co-ring. We denote by
$\sCoalg{P}{\op R}$ the subsemicategory of $\scat{Mod}_{(\op{P},\op{R})}$, the objects of which are
right $\op{P}$-modules of the form $\op{T}(A)$ for $A$ a $\op P$-coalgebra. A morphism in $\sCoalg{P}{\op R}$ is a morphism of
$(\op A,\op{P})$-bimodules,
\[
    \op{T}(A) \underset {\op P}\circ  \op{R} \rightarrow \op{T}(B).
\]
We call $\sCoalg{P}{\op R}$ the semicategory of $\op{P}$-coalgebras and
$\op{R}$-governed morphisms.

Similarly, we define the semicategory
${}_{(\op{P},\op{R})}\scat{Mod}$ of left $\op{P}$-modules and
$\op{R}$-governed morphisms, that becomes a category with a
functor ${}_{\op{P}}\cat{Mod} \rightarrow
{}_{(\op{P},\op{R})}\cat{Mod}$ if $\op{R}$ is counital. We may
furthermore define the full sub(semi)category $\sAlg{P}{\op R}$ of
$\op{P}$-algebras and $\op{R}$-governed morphisms.  A morphism in
$\sAlg{P}{\op R}$ is then a morphism of left $\op{A}$-modules,
\[
    \op{R} \underset {\op P}\circ  c(A) \rightarrow c(B).
\]
We define the composition of $\op{R}$-relative algebra morphisms
similarly, as the composition
\[
    \op{R} \underset {\op P}\circ  c(A) \rightarrow \op{R}
     \underset {\op P}\circ  \op{R}  \underset {\op P}\circ
     c(A) \rightarrow \op{R}  \underset {\op P}\circ  c(B)
     \rightarrow c(C).
\]
We may define in an analogous manner the category of $\op{P}$-algebras and
$\op{R}$-relative morphisms for the $z$ embedding (cf., Proposition \ref {prop:z}).

Note that a $\op{P}$-co-ring $\op{R}$ controls relative morphisms
of both algebras and coalgebras, without dualization. This is one
of the advantages of replacing duality by chirality.

Recall the Kleisli category $\cat{C}^{T}$ associated to a comonad
$T : \cat{C} \rightarrow \cat{C}$, ~\cite{kleisli:65}.  The
objects of $\cat{C}^{T}$ are the same as the objects of $\cat{C}$,
while the arrows are given by $\cat{C}^{T}(a,b) = \cat{C}(Ta,b)$.
We draw the following conclusions from the above discussion.

\begin{prop}\label{prop:kleisli-coalg}
Let $(\op R,\psi, \varepsilon)$ be a counital $\op P$-co-ring. The category $\cat{Mod}_{(\op{P},\op R)}$ is the Kleisli category
associated to the comonad $( - \underset {\op P}\circ  \op{R}, - \underset {\op P}\circ  \psi, -
\underset {\op P}\circ  \varepsilon)$.  The functor
$\op{T}(-)$ exhibits $\Coalg{P}{\op R}$ as a subcategory.
\end{prop}

\begin{prop}\label{prop:kleisli-alg}
Let $(\op R,\psi, \varepsilon)$ be a counital $\op P$-co-ring. The category ${}_{(\op{P},\op R)}\cat{Mod}$ is the Kleisli category
associated to the comonad $(\op{R} \underset {\op P}\circ  -, \psi \underset {\op P}\circ  -, \varepsilon
\underset {\op P}\circ  -)$. The functor $c(-)$
exhibits $\Alg{P}{\op R}$ as a full subcategory.
\end{prop}

Finally, we address the issue of level monoidal structures on
categories of modules.

\begin{prop}
If $\op{P}$ is a Hopf operad and  $\op{R}$ is a counital
$\op{P}$-co-ring in the category of level comonoids, then
$({}_{(\op{P},\psi)}\cat{Mod},\otimes,\op{C})$ and
$(\cat{Mod}_{(\op{P},\psi)},\otimes,\op{C})$ are
monoidal categories.
\end{prop}

The proof is formal.  The compatibility of the structure morphisms
in $\op{R}$ is necessary to make the tensoring of morphisms
compatible with composition.

\subsection{Modules over the associative operad}\label{ssec:modassocop}

In this section, we discuss the cosimplicial structures that arise
in $\op{A}$-bimodules, generalizing the notion of ``operad with
multiplication'' to $\op{A}$-bimodules ``with central morphism''.

Recall that $\op{A}$ is the \emph{reduced} associative operad,
with $\op{A}(n) = I[\Sigma_{n}]$ if $n > 0$ and $\op{A}(0) = O$.
Due to the lack of a unit operation, when we write
``cosimplicial'', we really mean ``precosimplicial'', since we are
missing the necessary codegeneracies.

Gerstenhaber and Voronov~\cite{gerstenhaber-voronov:95} noted that
an \emph{operad with multiplication} $\op{O}$ defines a
cosimplicial object $\op{O}^{\bullet}$; see also McClure and
Smith~\cite[Definition 3.1]{mcclure-smith:02}.  The coface
morphisms arise as a result of the $\op{A}$-bimodule structure
induced from a given morphism of operads $\op{A} \rightarrow
\op{O}$.  The only part of the ``pure'' operad structure in
$\op{O}$ that one needs in order to define the cosimplicial
structure is the unit, $\op{J} \rightarrow \op{O}$, which we replace with the following notion.

Let $\op{M}$ be an $\op{A}$-bimodule. We will call a morphism $v:I
\rightarrow \op{M}(1)$ \emph{central} if the following diagram
commutes.
\[
    \xymatrix{
        {\op{A}}(2)
            \ar[r]^-{\cong}
            \ar[d]^{\cong}
            & \op{A}(2) \otimes I \otimes I
            \ar[r]^-{1 \otimes v \otimes v}
            & \op{A}(2) \otimes \op{M}(1) \otimes \op{M}(1)
                \ar[d]^{\lambda} \\
        I \otimes \op{A}(2)
            \ar[r]_-{v \otimes Id}
            & \op{M}(1) \otimes \op{A}(2)
                \ar[r]_-{\rho}
                & \op{M}(2)
        }
\]
The terminology comes from a classical analogy: if $M$ is a
bimodule over a ring $R$, then $x \in M$ is central if $r \cdot x
= x \cdot r$ for all $r \in R$.

\begin{prop}\label{prop:right-module-cosimp}
Let $\op{M}$ be an $\op{A}$-bimodule. Any choice of central
morphism $v : I \rightarrow \op{M}(1)$ determines a cosimplicial
structure on $\op{M}$.
\end{prop}

\begin{proof}
Define the cosimplicial object $\op{M}^{\bullet}$ as follows. Set
$\op{M}^{n} = \op{M}(n)$. The coface operators $d^{i}$ will be
defined by the right $\op{A}$-module structure for $0 < i < n$,
and by the left $\op{A}$-module structure and the central morphism
$I \rightarrow \op{M}(1)$ for $i=0, n+1$.

Let $u_{q}:I \rightarrow I[\Sigma_{q}]$ be the
inclusion of the summand labelled by $1 \in \Sigma_{q}$, for $q=1,2$.
For $0 < i < n+1$, the coface operator $d^{i}:\op{M}(n)
\rightarrow \op{M}(n+1)$ is defined to be the composite
\[
    \op{M}(n) \cong \op{M}(n) \otimes I^{\otimes n}
        \rightarrow \op{M}(n) \otimes \op{A}(1)^{\otimes(i-1)} \otimes \op{A}(2)
        \otimes \op{A}(1)^{\otimes(n-i)}
        \xrightarrow{\rho} \op{M}(n+1).
\]
where the first arrow is $1 \otimes u_{1}^{\otimes(i-1)} \otimes
u_{2} \otimes u_{1}^{\otimes(n-i)}$.

Define $d^{0}$ to be the composite
\[
    \op{M}(n) \cong I \otimes I \otimes \op{M}(n)
        \xrightarrow{u_{2} \otimes v \otimes Id}
            \op{A}(2) \otimes\op{M}(1) \otimes \op{M}(n)
        \xrightarrow{\lambda}
            \op{M}(n+1)
\]
and $d^{n+1}$ to be the composite
\[
    \op{M}(n) \cong I \otimes \op{M}(n) \otimes I
        \xrightarrow{u_{2} \otimes Id\otimes v}
            \op{A}(2) \otimes\op{M}(n) \otimes \op{M}(1)
        \xrightarrow{\lambda}
            \op{M}(n+1).
\]
The identity $d^{i}d^{j} = d^{j}d^{i-1}$ for $n+2 > i > j > 0$
follows from the associativity of the right $\op{A}$-action.  For
$i=n+2$ and $j=0$, it follows from the associativity of the left
$\op{A}$-action.  For $i = n+2$ and $j=n+1$ or $i = 1$ and $j =
0$, it follows from the compatibility of the left and right
actions and the diagram defining a central morphism.  The remaining cases follow from the compatibility of the left and
right actions.
\end{proof}

\begin{rmk}
Kapranov and Manin~\cite[1.3.6 and 1.3.8]{kapranov-manin:01} have
also defined an ``augmented simplicial'' structure in any right
$\op{A}$-module.
\end{rmk}

Let $X$, $Y$ and $Z$ be (pre)cosimplicial objects in $\cat{M}$.

\begin{defn}\cite{mcclure-smith:02}
A \emph{cup-pairing} $\phi : (X^{\bullet},Y^{\bullet}) \rightarrow
Z^{\bullet}$ is a family of morphisms
\[
    \phi_{p,q} : X^{p} \otimes Y^{q} \rightarrow Z^{p+q}
\]
satisfying
    \begin{enumerate}
    \item\label{cond:d1} $d^{i}\phi_{p,q} =
        \left\{
            \begin{array}{ll}
                \phi_{p+1,q}(d^{i} \otimes Id_{Y})
                    & \text{if }i \leq p \\
                \phi_{p,q+1}(Id_{X} \otimes d^{i-p})
                    & \text{if }i > p
            \end{array}
        \right.$
    \item\label{cond:d2} $\phi_{p+1,q}(d^{p+1} \otimes Id_{Y})
        = \phi_{p,q+1}(Id_{X} \otimes d^{0})$.
    \end{enumerate}
\end{defn}

Evidently, a cup-pairing is a family of maps that yields a pairing
on the level of associated total objects.

\begin{prop}\label{prop:cosimp-mon}
Let $\op{M}$ be an $\op{A}$-bimodule with
central morphism. Then the graded multiplication in $\op{M}$
defines a natural cup-pairing in $\op{M}^{\bullet}$.
\end{prop}

\begin{proof}
We define $\phi_{p,q} : \op{M}(p) \otimes \op{M}(q) \rightarrow
\op{M}(p+q)$ by restricting the graded multiplication in $\op{M}$.
Then condition (\ref{cond:d1}) is satisfied by the compatibility
of the left and right operad actions, and for $d^{0}$ and
$d^{p+q+1}$, the associativity of the graded multiplication.
Similarly, condition (\ref{cond:d2}) follows by associativity.
\end{proof}

\begin{ex}
Let $C$ be a coaugmented coalgebra over the commutative ring $R$.
By Proposition~\ref{prop:right-module}, $\op{T}(C)$ is an
$\op{A}$-bimodule. The coaugmentation $R \rightarrow C$ induces a
unit $R \rightarrow \op{T}(C)(1)$.  If one checks the definitions,
one finds that the cosimplicial $R$-module $\op{T}(C)^{\bullet}$
is just the cosimplicial cobar construction on $C$, and
$\operatorname{Tot}(\op{T}(C)^{\bullet}) \cong \Omega C$.
\end{ex}

\section{Operadic spectroscopy}\label{sec:delooping}

For the remainder of this article, we work in $\cat{Ch}$, the closed, symmetric monoidal
category of connective chain complexes over a commutative ring
$R$, where  the tensor product of chain complexes over $R$ is denoted simply by $\otimes$.  Denote by $\csg{\otimes}$ the category of level
comonoids and of morphisms of symmetric sequences that respect the level coproduct.  Let
$\op{A}$ be the associative operad.  The purpose of this section
is to construct the \emph{diffracting functor}
\[
    \diffract(-) : \csg{\otimes} \rightarrow
    \cat{CoRing}_{\op A}.
\]
In the
following sections we apply the diffracting functor to proving the existence of structured maps between bar constructions or between cobar constructions.

\subsection{The diffracting functor}\label{ssec:delooping-functor}

The following conventions concerning the associative operad hold throughout the rest of this article.

\begin{notn}
Let $\op{A}$ be the associative operad of chain complexes. Let $\delta^{(n)} \in \op{A}(n)$ be the unit of $R[\Sigma_{n}]$ for all $n \geq 1$.  We let $\delta=\delta ^{(2)}$.  Note that
$$\gamma\big(\delta\otimes (\delta ^{(m)}, \delta ^{(n)})\big)=\delta ^{(m+n)}$$
 for all $m,n\geq 1$, where $\gamma:\op A\circ \op A\to \op A$ is the composition multiplication on the operad $\op A$.

Let $\vec{n} \in I_{m,n}$.  We set $\delta^{(\vec{n})} =
\delta^{(n_{1})} \otimes \cdots \otimes \delta^{(n_{m})}$.
\end{notn}

The $\op A$-bimodule underlying the diffraction of a level comonoid $\op X$ is a free $\op{A}$-bimodule generated by an appropriate
suspension of $\op X$. The symmetric sequence that suspends by tensoring is defined as follows.

\begin{defn}\label{defn:sphere}
For $n \geq 1$, let $\op{S}(n) = R\{ s_{n-1} \}$, the free graded
$R$-module generated by an element $s_{n-1}$ of degree $n-1$,
equipped with the sign representation of $\Sigma_{n}$.  Let
$\op{S}(0) = 0$.  We set
$$\Ssigma(n) = \op{S}(n) \otimes
\op{A}^{\sharp}(n),$$
where $\op{A}^{\sharp}(n)$ is the linear dual
of $\op{A}(n)$ for $n \geq 0$.

Let $\alpha_{n}=s_{n-1} \otimes 1 \in \op{S} \otimes \op{A}^{\sharp}(n)$, where $1$ is the unit of the ring $\op A^\sharp (n)$.
\end{defn}

The following proposition was proved in a simplicial version by Bauer, Johnson and
Morava in \cite{bjm:sphere}, based on an idea of Arone.  In its chain complex guise, the proof is quite simple and left to the reader.  It appears as well in \cite {markl-shnider-stasheff}Ê and in \cite {ginzburg-kapranov}.

\begin{prop}\label{prop:sphere}
The symmetric sequence $\op{S}$ is an operad and a cooperad.
\end{prop}

\begin{notn}
Given $m \geq 1$, let $\vec{n} = (n_{1}, \ldots, n_{m}) \in
I_{m,n}$.  We write $\alpha_{\vec{n}}$ for the product
$\alpha_{n_{1}} \otimes \cdots \otimes \alpha_{n_{m}}$.
\end{notn}

We use the following lemma  in defining one part of the differential on $\Phi(\op X)$.

\begin{lem}\label{lem:comonoid}
Let $(\op{X},\Delta)$ be a level comonoid. Then $(\op{S}
\otimes \op{X}) \boxprod \Ssigma$ is naturally a graded
comonoid.
\end{lem}

\begin{proof}
For $x \in \op{X}(m)$, let $\Delta(x) = \sum _{\ell}x_{\ell} \otimes
x^{\ell}$, and write $s_{m-1}x$ rather than $s_{m-1} \otimes x$ for
the corresponding element of $(\op{S} \otimes \op{X})(m)$.  The
desired diagonal $\widetilde \Delta$ on $(\op{S}
\otimes \op{X}) \boxprod \Ssigma$  is the $\Sigma_{n}$-equivariant morphism defined
by the formula
\[
    \widetilde\Delta(s_{m-1}x \otimes \alpha_{\vec{n}})
        = \sum_{\vec{\imath}+\vec{\jmath}=\vec{n}}
            \sum_{\ell}
        \left(
            s_{m-1} x_{\ell} \otimes \alpha_{\vec{\imath}}
        \right) \otimes
        \left(
            s_{m-1} x^{\ell} \otimes \alpha_{\vec{\jmath}}
        \right)
\]
for $\vec{n} \in I_{m,n}$, with implicit signs given by the Koszul rule. Note that the diagonal commutes with
differentials since the diagonal in $\op{X}$ does so. Furthermore,
since morphisms of level comonoids commute with iterated
diagonals, the construction is natural.
\end{proof}

Let $\op{X}$ be a level comonoid.  Let $\op{V}_{\op{X}} =
(\op{S} \otimes \op{X}) \circ \Ssigma$.  Since $\op{V}_{\op{X}}$
is a graded comonoid by Lemma~\ref{lem:comonoid}, it is an
$\op{A}_{\odot}$-coalgebra.  Thus $\op{T}_{\odot}\op{V}_{\op{X}}$
is a $\op{A}_{\odot}$-bimodule.  The zero morphism $\op{U}
\rightarrow (\op{T}_{\odot}\op{V}_{\op{X}})(1)$ is central, and
determines the structure of cosimplicial chain complex with
cup-pairing on $\op{T}_{\odot}\op{V}_{\op{X}}$.  The associated
chain complex is $\op{A} \circ \op{V}_{\op{X}}$ as a symmetric
sequence of graded modules.  The differential $\partial_{c}$ is
the sum of the internal differential induced from the differential $d$ on  $\op{X}$ and
the cobar differential from the diagonal in $\op{V}_{\op{X}}$.  More explicitly, for all $x\in \op X(m)$ with $\Delta (x)=\sum_{\ell}x_{\ell}\otimes x^\ell$ and for all $\vec n\in I_{m,n}$,
$$\del _{c}(s_{m-1}x\otimes \alpha_{\vec n})=s_{m-1}dx \otimes \alpha _{\vec n}+\sum _{\vec \imath +\vec\jmath=\vec n}\sum_{\ell}\delta\otimes\big((s_{m-1}x_{\ell}\otimes \alpha_{\vec\imath})\otimes (s_{m-1}x_{\ell}\otimes \alpha_{\vec\jmath})\big),$$
with implicit signs determined by the Koszul rule.

\begin{defn} Let $\op X$ be a level comonoid of chain complexes.  The \emph{cobar construction} on $\op X$ is the left $\op A$-module $(\op A\circ \op V_{\op X},\del _{c})$ constructed above.
\end{defn}

Now we define a presimplicial structure on $\Ssigma \circ \op{A}(n)$, for each $n$, which gives rise to another piece of the differential of $\Phi(\op X)$ .
Recall that
\[
    (\Ssigma \circ \op{A})(n)
         =  \bigoplus_{m \geq 0} \Ssigma(m)
            \otimes_{\Sigma_{m}}
            \op{A}^{\odot m}(n).
\]
For all $m\geq 0$, let
$$\big(\Ssigma \circ \op{A}(n)\big)_{m}= \Ssigma(m)\otimes_{\Sigma_{m}}\op{A}^{\odot m}(n).$$  For all $m\geq 1$ and and all $0 \leq i
\leq m$, define face maps
$$d_{i}: \big(\Ssigma \circ \op{A}(n)\big)_{m}\to \big(\Ssigma \circ \op{A}(n)\big)_{m-1}$$ by
$d_{i}(\alpha_{m} \otimes \delta^{(\vec{n})})=\alpha_{m-1} \otimes \delta^{(\vec{n}(i))}$, where
$\vec{n}(i) \in I_{m-1,n}$ is given by
\[
    \vec{n}(i)
        = \left\{
            \begin{array}{cl}
                \varnothing   & \text{if }i=0, \\
                (n_{1},\ldots,
                n_{i}+n_{i+1}, \ldots, n_{m})
                    & \text{if }1 \leq i < m, \\
                \varnothing & \text{if }i=m.            \end{array}
        \right.
\]
where $\delta^{(\varnothing)}=1 \in R$. The resulting morphisms are $\Sigma_{n}$-equivariant.

\begin{rmk}
The chain complex associated to $\op{A}^{\bot} \circ \op{A}$ is
precisely the one-sided Koszul resolution of $\op{A}$.  Indeed,
consider the cooperad structure morphism $$\op{A}^{\bot}(m)
\rightarrow \op{A}^{\bot}(m-1) \otimes
(\op{A}^{\bot})^{\odot(m-1)}(m).$$
Since the difference between
$(\op{A}^{\bot})^{\odot(m-1)}(m)$ and $\op{A}^{\odot (m-1)}(m)$ is
a single suspension, desuspension defines a morphism of degree
$-1$, $\op{A}^{\bot}(m) \rightarrow \op{A}^{\bot}(m-1) \otimes
\op{A}^{\odot(m-1)}(m)$, that extends linearly to define a
differential of right $\op A$-modules on $\op{A}^{\bot} \circ \op{A}$, which
is precisely the one defined above. In fact, the simplicial
structure outlined above is just the simplicial bar resolution of
the monoid $\cat{N}$.
\end{rmk}

The levelwise differential on $\op A^\bot\circ \op A$ resulting from the presimplicial structure on each $(\op A^\bot\circ \op A)(n)$ extends to define a differential of $\op A$-bimodules
$\partial_{s}$ on $\op A\circ \op{V}_{\op{X}}\circ \op A$. A straightforward
calculation establishes the following proposition.

\begin{prop}\label{prop:diffl-mho}
The cosimplicial  and simplicial differentials together define a
differential $\partial = \partial_{c} + \partial_{s}$ on
$\op A\circ \op{V}_{\op{X}}\circ \op A$.
\end{prop}

\begin{rmk}\label{rmk:formula-diffl}
The formula for the full differential $\partial$, modulo signs, is
\begin{align*}
    \partial(s_{m-1}x \otimes \alpha_{\vec{n}})
        =& s_{m-1}d x\otimes \alpha_{\vec{n}}\\
        &+ \sum_{\vec{\imath} + \vec{\jmath}=\vec{n}} \sum _{\ell}\delta \otimes\big(
        \left( s_{m-1}x_{\ell} \otimes \alpha_{\vec{\imath}} \right)
        \otimes
        \left( s_{m-1}x^{\ell} \otimes \alpha_{\vec{\jmath}} \right)\big) \\
        &+ \sum_{j=1}^{n-1} s_{m-1}x \otimes \alpha_{\vec{n}-\vec{e}_{j}}
            \otimes (1^{\otimes(j-1)} \otimes \delta \otimes
            1^{\otimes (n-j)})
\end{align*}
where $x\in \op X(m)$, $\Delta(x) = \sum x_{\ell} \otimes x^{\ell}$, and
$\vec{e}_{j}$ is the unit vector that indicates which subinterval
of the partition $\vec{n}$ contains $j$. More precisely, let $r$
satisfy $n_{1} + \cdots + n_{r-1} + 1 \leq j \leq n_{1} + \cdots +
n_{r}$. Then $\vec{e}_{j}$ is the $m$-tuple with zeros everywhere
except for a $1$ in the $r$th position. The signs are determined
by the Koszul convention and the usual rules for differentials
coming from simplicial face maps.
\end{rmk}

\begin{defn}
The \emph{diffracting functor}
\[
    \diffract : \csg{\otimes}\rightarrow
        {}_{\op{A}}\cat{Mod}_{\op{A}}
\]
is defined on objects by $\Phi(\op X)=(\op A\circ \op V_{\op X}\circ \op A, \del)$.  If $\varphi : \op{X} \rightarrow
\op{Y}$ is a morphism of level comonoids, then
 $\diffract(\varphi) = Id_{\op{A}} \circ (
Id_{\op{S}} \otimes \varphi ) \circ Id_{\op{A}^\perp} \circ Id_{\op{A}}$.
\end{defn}

The following proposition states that $\diffract(\op{X})$ admits a natural $\op A$-co-ring structure.

\begin{prop}\label{prop:mho-comonoidal}
For all level comonoids $\op X$, the $\op A$-bimodule $\Phi (\op X)$ admits a natural coassociative coproduct $\psi _{\op X}\in {}_{\op{A}}\cat{Mod}_{\op{A}}\big(\Phi(\op X), \Phi(\op X)\acirc \Phi(\op X)\big)$.
\end{prop}

\begin{proof}
For $r,m \geq 0$, let $\vec{m} \in I_{r,m}$.  Suppose that
$\vec{m}_{i} = (m_{i}^{1} , \ldots , m_{i}^{\ell} ) \in
I_{\ell,m_{i}}$ for some $\ell \geq 1$ and for $1 \leq i \leq r$.
For $1 \leq j \leq \ell$, set
\[
    \vec{m}^{j} = (m_{1}^{j},\ldots, m_{r}^{j}) \in I_{r,m_{j}},
\]
where $m^{j} = \sum_{i=1}^{r}m_{i}^{j}$.  Note that
$\sum_{j=1}^{\ell}m^{j} = m$.  Let $x \in \op{X}(r)$, and denote
by $\Delta^{(\ell)} x = \sum x_{0} \otimes \cdots \otimes
x_{\ell}$ the iterated diagonal. Suppose that $\vec{\ell} \in
I_{r,\ell}$.  We define
\begin{multline*}
    \psi_{\vec{\ell},(\vec{m}_{i})} ( s_{r-1}x \otimes
    \alpha_{\vec{m}} )\\
    = (s_{r-1} x_{0} \otimes \alpha_{\vec{\ell}})
    \otimes \big(( s_{r-1}x_{1} \otimes \alpha_{\vec{m}^{1}} )
    \otimes \cdots
    \otimes ( s_{r-1}x_{\ell} \otimes \alpha_{\vec{m}^{\ell}})\big),
\end{multline*}
with implicit signs arising from application of the Koszul rule.  The map $\psi_{\vec{\ell},(\vec{m}_{i})}$ is of degree zero
because $\sum_{j=1}^{\ell} m^{j} = m$.  These morphisms have the
proper equivariance properties to give rise for all $m$ to a morphism
\[
    \left(( \op{S} \otimes \op{X} ) \circ \Ssigma \right) (m)
    \xrightarrow{\psi_{m}}
    \prod_{\ell \geq 0}
    \left(( \op{S} \otimes \op{X} ) \circ \Ssigma \right)(\ell)
    \otimes_{\Sigma_{\ell}}
    \left(( \op{S} \otimes \op{X} ) \circ \Ssigma \right)^{\odot \ell}(m),
\]
which we gather into $\psi_{\op X}:\op V_{\op X}\to \op V_{\op X}\circ\op V_{\op X}.$
The image falls in the direct sum since the symmetric sequences
vanish in level zero. The composite
\[
   \op V_{\op{X}} \xrightarrow{\psi}  \op V_{\op{X}}\circ  \op V_{\op{X}}
    \xrightarrow{\cong} \op{J} \circ \op V_{\op{X}} \circ \op{J} \circ
    \op V_{\op{X}} \circ \op{J} \xrightarrow{\eta \circ Id \circ \eta \circ
    Id \circ \eta} \op{A} \circ  \op V_{\op{X}}\circ \op{A} \circ  \op V_{\op{X}} \circ \op{A}
\]
we also denote by $\psi_{\op X}$, likewise for the extension of $\psi_{\op X}$ to
the $\op{A}$-bilinear morphism
$$\op A\circ\op V_{\op{X}}\circ\op A\rightarrow
\op A\circ\op V_{\op{X}}\circ\op A\acirc \op A\circ\op V_{\op{X}}\circ\op A.$$
The diagonal $\psi_{\op X}$
so defined commutes with the internal differential by inspection,
with the simplicial differential by direct calculation, and with
the cosimplicial differential because $\Delta$ is coassociative.  Furthermore, the coassociativity of $\psi _{\op X}$ follows from that of the level coproduct on $\op X$.
\end{proof}

 \begin{cor} The diffracting functor $\Phi$ corestricts to a functor
 $$\Phi:\cat {Comon}_{\otimes}\longrightarrow \cat{CoRing}_{\op A}.$$
 \end{cor}

\subsection{The generalized Milgram map}\label{ssec:milgram}

The diffraction of a level comonoid admits a diagonal, via the
\emph{generalized Milgram map} that we define below.  It is
analogous to the Milgram equivalence of cobar constructions, $q :
\Omega( A \otimes B ) \xrightarrow{\simeq} \Omega(A) \otimes
\Omega(B)$, defined in~\cite{milgram:66}.

Let $\op{X}$ and $\op{Y}$ be level comonoids. Then $\op{X}
\otimes \op{Y}$ is again a level comonoid, as remarked in section \ref{ssec:bimodules} .  Let $\ell \geq 1$,
$m \geq 0$, and suppose that $\vec{m} \in I_{\ell,m}$.  If $x \in
\op{X}(\ell)$ and $y \in \op{Y}(\ell)$, then $s_{\ell-1}(x \otimes
y) \otimes \alpha_{\vec{m}} \in \diffract(\op{X} \otimes \op{Y})(m)$.
For $n \geq 1$, we abuse notation slightly and write $\Delta^{(n-1)}y = \sum y_{0} \otimes
\cdots\otimes y_{n-1}$ for the iterated diagonal of $y$.  We define
\[
    q(s_{\ell-1}(x \otimes y) \otimes \alpha_{\vec{m}})
    = \sum (s_{\ell-1}x \otimes \alpha_{\vec{n}} \otimes
    \delta^{(\vec{r})}) \otimes
    (\delta^{(n)} \otimes s_{\ell-1}y_{0} \otimes \cdots
    \otimes s_{\ell-1}y_{n-1} \otimes \alpha_{\vec{r}})
\]
where the sum is over all $n \geq 1$, $\vec{n} \in I_{\ell,n}$,
and $\vec{r} \in I_{n,m}$, and signs are given implicitly via the
usual Koszul convention. We use the notation $\delta^{(\vec{r})} =
\delta^{(r_{1})} \otimes \cdots \otimes \delta^{(r_{n})}$.  Note
that $q$ is constructed via simplicial and cosimplicial operators.
A straightforward, yet quite lengthy, calculation shows that $q$
commutes with the underlying face and coface maps of
$\diffract(\op{X} \otimes \op{Y})$ and $\diffract{\op{X}}\otimes
\diffract{\op{Y}}$.  It follows that $q$ is a morphism in
$\bimod{A}{A}$.

Since the construction of $q$ is clearly natural in both $\op X$ and $\op Y$, we can formulate the following definition.

\begin{defn}
The \emph{generalized Milgram map} is the natural transformation
\[
    q : \diffract(-\otimes-) \rightarrow \diffract(-) \otimes \diffract(-)
\]
of functors
\[
    \csg{\otimes} \times \csg{\otimes} \rightarrow
    {}_{\op{A}}\cat{Mod}_{\op{A}}.
\]
\end{defn}

Let $(\op{X},\Delta)$ be a level comonoid.  If $\Delta$ is a
morphism of level comonoids, then we define a diagonal on
$\diffract{\op{X}}$ to be the composite
\[
    \diffract(\op{X}) \xrightarrow{\diffract{\Delta}}
        \diffract{(\op{X} \otimes \op{X})}
        \xrightarrow{q} \diffract(\op{X}) \otimes
        \diffract(\op{X}).
\]
To determine whether or not this diagonal is coassociative
requires a greater understanding of the kernel of $q$ than exists
at the time of writing.

\section{Duality
theorems}\label{sec:induction-linearization}

In this section, we present the Cobar and Bar Duality Theorems.
These theorems exploit the diffracting functor defined in the
previous section to recognize additional comultiplicative
structure in an algebra morphism $\Omega C \rightarrow \Omega C'$ and multiplicative structure in a coalgebra morphism $BA
\rightarrow BA'$ of bar constructions.

The Cobar and Bar Duality Theorems can be seen as a generalization
of the work of Gugenheim and
Munkholm~\cite{gugenheim-munkholm:74}.  Recall from section \ref{ssec:sigma-modules}  the
symmetric sequence $\op{J}$, which is the unit with respect to the composition product.  It is trivially an operad,  whose coalgebras are all chain complexes.
A chain algebra is thus a monoid in the category of
$\op{J}$-coalgebras. Let $C$ and $C'$ be connected chain
coalgebras and let $\Omega(-)$ be the cobar functor.  Gugenheim
and Munkholm showed that a chain algebra morphism $\Omega C
\rightarrow \Omega C'$ is equivalent to a strongly
homotopy-comultiplicative map $C \rightarrow C'$. The Cobar
Duality Theorem implies that $\Omega C \rightarrow \Omega C'$ is a
morphism of monoids in the category of $\op{J}$-coalgebras if and
only if $C \rightarrow C'$ is a $\Phi(\op{J})$-governed morphism.
As we shall see in section~\ref{sec:aw}, a $\Phi(\op{J})$-governed
morphism is precisely a strongly homotopy-comultiplicative map.

\subsection{Interaction between level monoids and
comonoids}\label{ssection:mon-comon}

In order to state the duality theorems, we need to introduce some
terminology.  The terminology involves morphisms of the form
\[
    \op{M} \otimes \op{X} \rightarrow \op{N}
\]
where $\op{X}$ is a level comonoid while $\op{M}$ and $\op{N}$
are both either level monoids or level comonoids.  The
classical case that serves as our template is the following.

\begin{defn}
Let $(A,\mu_{A})$ and $(B,\mu_{B})$ be monoids in $\cat{M}$, and
let $(C,\Delta_{C})$ be a comonoid.  We call a morphism $\varphi:A
\otimes C \rightarrow B$ \emph{multiplicative} if the diagram
\[
    \xymatrix{
        {A \otimes A \otimes C}
            \ar[rr]^{\mu_{A} \otimes 1}
            \ar[d]_{1 \otimes \Delta_{C}}
            & &     {A \otimes C}
                        \ar[dd]^{\varphi}    \\
        {A \otimes A \otimes C \otimes C}
            \ar[d]_{1 \otimes \tau \otimes 1}
                                    \\
        {A \otimes C \otimes A \otimes C}
            \ar[r]_-{\varphi \otimes \varphi}
            & {B \otimes B}
                \ar[r]_-{\mu_{B}}
            & {B}
        }
\]
commutes.
\end{defn}

The following lemma provides a way of constructing multiplicative
morphisms.

\begin{lem}\label{lem:multiplicative}
Let $X$ be an object in $\cat{M}$, and let $B$ and $C$ be a monoid
and a comonoid, respectively.  Let $\op{A}(X)$ be the free
$\op{A}$-algebra on $X$.  Then a morphism $\hat{\varphi} : X
\otimes C \rightarrow B$ extends to a unique multiplicative
morphism $\varphi : \op{A}(X) \otimes C \rightarrow B$.
\end{lem}

\begin{proof} Let $\Delta ^{(n-1)}:C\to C^{\otimes n}$ and $\mu ^{(n-1)}:B^{\otimes n}\to B$ denote the iterated coproduct on $C$ and the iterated product on $A$, respectively.  Recall that $\op A(X)= \bigoplus _{n\geq 1}X^{\otimes n}$.  Define the restriction of $\varphi$ to  $X^{\otimes n}\otimes C$ to be given by the following composite.
$$X^{\otimes n}\otimes C\xrightarrow{Id\otimes \Delta^{(n-1)}}X^{\otimes n}\otimes C^{\otimes n}\xrightarrow{\text{shuffle}}(X\otimes C)^{\otimes n}\xrightarrow{\hat{\varphi}}B^{\otimes n}\xrightarrow {\mu^{(n-1)}}B$$
Note that both the coassociativity of $\Delta$ and the associativity of $\mu$ are crucial in this definition.
\end{proof}

Let $(\op{X}, \Delta_{\op{X}})$ be a level comonoid.  Let
$\theta : \op{M} \otimes \op{X} \rightarrow \op{N}$ be a morphism
of symmetric sequences.

Suppose that $\op{M}$ and $\op{N}$ are level monoids with
multiplications $\mu_{\op{M}}$ and $\mu_{\op{N}}$, respectively.
We say that $\theta$ is \emph{multiplicative} if the diagram
\[
    \xymatrix{
        {\op{M}^{\otimes 2} \circ \op{X}}
            \ar[d]^{\mu_{\op{M}}\circ Id}
            \ar[rr]^{Id \circ \Delta_{\op{X}}}
            && {\op{M}^{\otimes 2}\circ \op{X}^{\otimes 2}}
                \ar[r]^\iota
            & {(\op{M}\circ \op{X})^{\otimes 2}}
                \ar[r]^-{\theta^{\otimes 2}}
            & {\op{N}^{\otimes 2}}
                \ar[d]^{\mu_{\op{N}}} \\
        {\op{M}\circ\op{X}}
            \ar[rrrr]^{\theta}
            &&&& {\op{N}}
    }
\]
commutes.

The next lemma provides  a process for creating level monoids.

\begin{lem}\label{lem:free-sym-mon}
The free graded monoid on a monoid in $\cat{M}$ is naturally a
level monoid in $\cat{M}^{\Sigma}$, i.e., $\op T$ restricts to a functor
$$\op T:\cat {Mon(M)}\longrightarrow \cat{Mon}_{\otimes}.$$
\end{lem}

\begin{proof}
If $A$ is a monoid in $\cat{M}$, then so too is $A^{\otimes n} =
\op{T}(A)(n)$ for any $n \geq 1$.
\end{proof}

The two notions of ``multiplicative'' above are connected by the
following lemma that we shall need for the Cobar Duality Theorem.

\begin{lem}\label{lem:sym-seq-multi}
Let $A$ and $B$ be monoids in $\cat M$, and let $\op{X}$ be a level comonoid.
For $m \geq 1$, suppose that $\theta_{m} : A \otimes \op{X}(m)
\rightarrow B^{\otimes m}$ is a multiplicative morphism that is
$\Sigma_{m}$-equivariant with respect to the right action on
$\op{X}(m)$ and the right permutation action on $B^{\otimes m}$.
Then the sequence $\{ \theta_{m} \}$ defines a multiplicative
morphism of symmetric sequences $\theta : \op{T}(A) \circ \op{X}
\rightarrow \op{T}(B)$.
\end{lem}

\begin{proof}
The sequence $\{ \theta_{m} \}$ defines a morphism of symmetric
sequences by Proposition~\ref{prop:extending-to-box}.  The
morphism is multiplicative since each $\theta_{m}$ is.
\end{proof}

\subsection{Cobar duality}

We state and prove in this section the first  of the two important duality results that justify our study of the diffracting functor.  Before doing so, we first recall the cobar construction as it applies to noncounital coalgebras, then introduce a bit of helpful notation.

Recall that we have defined the associative operad $\op A$ so that $\op A(0)=0$.  As a result, $\op A$-(co)algebras are non(co)unital.  An $\op A$-coalgebra consists therefore of a chain complex $C$ together with a coassociative chain map $\Delta:C\to C\otimes C$.

Recall further that $\cat {Ch}$ denotes the category of connective (i.e., bounded-below) chain complexes over a commutative ring, endowed with the usual graded tensor product.

The next definition generalizes slightly the usual definition of the cobar construction.

\begin{defn}  The \emph{cobar construction}  is the functor $\Om: \coalg{A}\to \alg A$ defined as follows. For all $C\in \ob \coalg{A}$, the graded $\op A$-algebra underlying $\Om C$ is $\op A(s^{-1} C)$, i.e., $\bigoplus_{n\geq 1}\big(s^{-1}C\big)^{\otimes n}$, where $s^{-1}C$ is the desuspension of $C$, i.e., $(s^{-1}C )_{n}=C_{n+1}$.  The differential $d_{\Om }$ on $\Om C$ is the derivation specified by
$$d_{\Om}s^{-1}=-s^{-1}d +(s^{-1}\otimes s^{-1})\Delta:C\to \op A(s^{-1}C),$$
where $d$  and $\Delta$ are the differential and coproduct of $C$.
\end{defn}

\begin{rmk} Observe that $C_{*}\star \cong R$, where $\star$ is the simplicial one-point set. The normalized chain coalgebra $C_{*}K$ on any simplicial set $K$ admits a counit, given by the chain map $C_{*}K\to C_{*}\star$ induced by the unique simplicial map $K\to \star$.  If $K$ is a pointed simplicial set, i.e., $K$ has a distinguished $0$-simplex, then $C_{*}K$ is  coaugmented by the chain map $C_{*}\star\to C_{*}K$ induced by the simplicial map $\star\to K$ picking out the basepoint.  Thus, if $K$ is a pointed simplicial set, then the short exact sequence $0\to C_{*}\star\to C_{*}K\to \overline C_{*}K\to 0$ of chain coalgebras splits, and we can write $C_{*}K\cong R\oplus \overline C_{*}K$, as chain coalgebras.   Furthermore, any easy computation leads to the conclusion that $\Om R$ is acyclic, so that the projection map $\Om R\otimes \Om \overline C\to \Om \overline C$ is a quasi-isomorphism.  Since the Milgram map $\Om (R\oplus \overline C)\to \Om R\otimes \Om \overline C$ is also a quasi-isomorphism (cf., Appendix in \cite {hps2}), we conclude that the composite $\Om C\to \Om \overline C$ is a quasi-isomorphism.  In other words, we can keep or throw away the counit without changing the homology of the cobar construction.
\end{rmk}

We next introduce an important class of morphisms of symmetric sequences, which we need in expressing our duality theorem for the cobar construction. Let $A,B \in \ob\cat{M}$, $\op{X},\op{Y} \in\ob\cat{M}^{\Sigma}$. Let
\[
    \Theta = \{ \theta_{\vec{n}} : A \otimes \op{X}(m) \otimes
    \op{Y}[\vec{n}] \rightarrow B^{\otimes n_{1}} \otimes \cdots
    \otimes B^{\otimes n_{m}} \}
\]
be a family of equivariant maps in $\cat{M}$, with $m \geq 1$ and $\vec{n} =
(n_{1}, \ldots, n_{m})$.

We define morphisms in $\cat{M}$,
\[
    \theta_{\vec{n_{1}},\cdots,\vec{n_{k}}}:
    \op{T}(A)(k) \otimes \left(
        \left( \op{X}(m_{1}) \otimes \op{Y}[\vec{n}_{1}]
        \right)
        \otimes \cdots \otimes
        \left( \op{X}(m_{k}) \otimes \op{Y}[\vec{n}_{k}]
        \right)
        \right)
        \rightarrow \op{T}(B)(n)
\]
by
\[
    \xymatrix{
        {A^{\otimes k} \otimes_{j=1}^{k}
            \left( \op{X}(m_{j}) \otimes \op{Y}[\vec{n}_{j}]
            \right)}
            \ar[r]^{\text{shuffle}}
            & {\otimes_{j=1}^{k}
            \left(
                A \otimes \op{X}(m_{j}) \otimes
                \op{Y}[\vec{n}_{j}]
            \right)
            }
            \ar[d]_{\otimes_{j}\theta_{\vec{n}_{j}}}
            \\
        &   {\otimes_{j=1}^{k}
            \left(
                B^{\otimes n_{j1}} \otimes \cdots \otimes
                B^{\otimes n_{jm_{j}}}
            \right)}
            \ar[d]_{\text{transpose}}
            \\
        &   {\otimes_{i = 1}^{m}
            \left( B^{\otimes n_{1i}} \otimes \cdots \otimes
            B^{\otimes n_{ki}} \right)}
    }
\]
where $m = \max(m_{1},\ldots,m_{k})$,  $n_{ji}=0$ if $i > m_{j}$,
and $B^{\otimes 0} = I$.

\begin{prop}\label{prop:tt-morph}
The morphisms $\theta_{\vec{n_{1}},\cdots,\vec{n_{k}}}$ are
equivariant with respect to the actions of $\Sigma_{k}$,
$\Sigma_{m_{1}} \times \cdots \times \Sigma_{m_{k}}$, and
$\Sigma_{n_{11}} \times \cdots \times \Sigma_{n_{km_{k}}}$, and so
determine a unique morphism of symmetric sequences,
\[
    \theta : \op{T}(A) \circ (\op{X} \circ \op{Y}) \rightarrow
    \op{T}(B).
\]
\end{prop}

\begin{defn}\label{def:tt-morph}
We will call the morphism $\theta$ of
Proposition~\ref{prop:tt-morph} the \emph{transposed tensor
morphism of symmetric sequences} induced by $\Theta$.

Let $\op X$ be a level comonoid, and let $C$ and $C'$ be chain coalgebras.  In a slight variation on this terminology, a morphism  of right $\op A$-modules $\theta: \op T(C)\acirc \Phi(\op X)\to \op T(C')$ is called a \emph{transposed tensor morphism of right $\op A$-modules} if its restriction  $\theta: \op T(C)\circ (\op S\otimes \op X)\circ \op A^\bot\to \op T(C')$ is a transposed tensor morphism of symmetric sequences of graded $R$-modules.
\end{defn}

We record here for future use a helpful result on sums of transposed tensor maps.

\begin{prop}\label{prop:sum-transptens}  Let $C, C', C''$ be objects in $\coalg A$, and let $\op X$ be a level comonoid. Let $\theta':\op T(C')\acirc \Phi(\op X)\to \op T(C)$ and $\theta'':\op T(C'')\acirc \Phi(\op X)\to \op T(C)$ be transposed tensor morphisms of right $\op A$-modules.  Then there is a unique transposed tensor morphisms of right $\op A$-modules $\theta:\op T(C'\oplus C'')\acirc \Phi(\op X)\to \op T(C)$ restricting to $\theta'$ and $\theta''$.
\end{prop}

\begin{proof}  The transposed tensor morphisms $\theta'$ and $\theta''$ are determined by families of equivariant morphisms in $\cat M$
\[
    \Theta' = \{ \theta'_{\vec{n}} : C' \otimes (\op S\otimes \op{X})(m) \otimes
    \op{A}^\perp[\vec{n}] \rightarrow C^{\otimes n_{1}} \otimes \cdots
    \otimes C^{\otimes n_{m}} \}_{m, \vec n}
\]
 and
 \[
    \Theta'' = \{ \theta''_{\vec{n}} : C'' \otimes (\op S\otimes \op{X})(m) \otimes
    \op{A}^\perp[\vec{n}] \rightarrow C^{\otimes n_{1}} \otimes \cdots
    \otimes C^{\otimes n_{m}} \}_{m, \vec n}.
\]
These  in turn give rise to a family of equivariant morphisms
$$ \Theta = \{ \theta_{\vec{n}} : (C'\oplus C'') \otimes (\op S\otimes \op{X})(m) \otimes
    \op{A}^\perp[\vec{n}] \rightarrow C^{\otimes n_{1}} \otimes \cdots
    \otimes C^{\otimes n_{m}} \}_{m, \vec n},
$$
since the monoidal product in $\cat M$ preserves coproducts.  The family $\Theta$ determines the desired transposed tensor morphism.
\end{proof}

Let $\msym(\op Y\circ \op X, \op Z)$ denote the set of multiplicative morphisms of symmetric sequences, where $\op X$ is a level comonoid, and $\op Y$ and $\op Z$ are level monoids. Let $\cat{Mod}_{\op A}^{\text{tt}}\big(\op T(C)\acirc \Phi (\op X),\op T(C')\big)$ denote the set of transposed tensor morphisms of right $\op A$-modules  from $\op T(C)\acirc \Phi(\op X) $ to $\op T(C')$, where $\op X$ is a level comonoid, and $C$ and $C'$ are chain coalgebras.

Recall that if $\cat C$ is a category and $\cat D$ is a small category, then $\cat {C^D}$ denotes the category of functors from $\cat D$ to $\cat C$ and of natural transformations between them.

\begin{DT} Let $\cat D$ be any small category.
There are mutually inverse, natural isomorphisms
$$\ind: \ttma\big(\op T(-)\acirc \Phi(-), \op T(-)\big)\longrightarrow \msym \big(\op T(\Om-)\circ -,\op T(\Om-)\big)$$
and
$$\lin :\msym \big(\op T(\Om-)\circ -,\op T(\Om-)\big)\longrightarrow\ttma\big(\op T(-)\acirc \Phi(-),\op T(-)\big)$$
of  functors from $(\coalg{A})^{\cat D}\times\cat {Comon}_{\otimes}\times (\coalg{A})^{\cat D}$ to  $\cat{Set}^{\cat D}.$
\end{DT}

The natural isomorphism  $\lin$ is called \emph{linearization}, as it transforms a family of maps between cobar constructions to a family of maps between their generators.  Its inverse, $\ind$, is called \emph {induction}, since it creates a family of maps between cobar constructions that is induced from a family of maps between their generators.

\begin{proof}  We first construct the induction transformation.  Let $F,G\in\ob (\coalg{A})^{\cat D}$, and let $\op X\in\ob \cat {Comon}_{\otimes}$.

Fix $m \geq 1$ and let $\vec{n} = (n_{1}, \ldots, n_{m})$ be an
$m$-tuple of natural numbers.  As usual, $n = n_{1} + \cdots+
n_{m}$.  For brevity, as in section 3, we use the notation $\alpha_{\vec{n}} =
\alpha_{n_{1}} \otimes \cdots \otimes \alpha_{n_{m}}$ for elements of $\op A^\perp[\vec n]$.  Let
$\eta_{\vec{n}} : \op{X}(m) \rightarrow (\op{S}(m) \otimes
\op{X}(m)) \otimes_{\Sigma_{m}} \op{A}^\perp[\vec{n}]$ be the map of
degree $n-1$ defined by $\eta_{\vec{n}}(x) = s_{m-1}x \otimes
\alpha_{\vec{n}}$.

Let $\theta \in\ttma\big(\op T(F)\acirc \Phi(\op X), \op T(G)\big)$.  Define
$\theta_{\vec{n}}^{\bot} : s^{-1}F \otimes \op{X}(m) \rightarrow
(\Omega G)^{\otimes m}$ to be the composite
\[
    \xymatrix{
        {s^{-1}F \otimes \op{X}(m)}
            \ar[r]^-{s \otimes \eta_{\vec{n}}}
            & {F \otimes (\op{S}(m) \otimes \op{X}(m))
            \otimes_{\Sigma_{m}} \op{A}^\perp[\vec{n}]}
                \ar[d]_{\theta_{\vec{n}}} \\
        & G^{\otimes n}
            \ar[d]_{(s^{-1})^{\otimes{n}}} \\
        & (s^{-1}G)^{\otimes n_{1}} \otimes \cdots \otimes
        (s^{-1}G)^{\otimes n_{m}}
            \ar[d]_{\iota_{\vec{n}}}   \\
        & (\Omega G)^{\otimes m},
    }
\]
where $\iota _{\vec n}$ is the obvious factor-by-factor inclusion.  Notice that $\theta_{\vec{n}}^{\perp}$ is homogeneous of degree
$0$.  Let
$$\theta_{m}^{\perp} = \sum_{\vec{n}}
\theta_{\vec{n}}^{\perp} : s^{-1}F \otimes \op{X}(m) \rightarrow
(\Omega G)^{\otimes m}.$$

Denote by $\Delta^{(n-1)} : \op{X} \rightarrow \op{X}^{\otimes n}$
the iterated diagonal in $\op{X}$.  Denote by $\op{A}(-)$ the free
$\op{A}$-algebra functor, as usual.  Extend $\theta_{m}^{\perp}$
multiplicatively to $\op{A}(s^{-1} F) \otimes \op{X}(m)$, as explained in the previous section.
\[
    \xymatrix{
        {(s^{-1}F)^{\otimes n} \otimes \op{X}(m)}
            \ar[r]^{1 \otimes \Delta^{(n-1)}}
            \ar@{.>}[ddrr]
            & {(s^{-1}F)^{\otimes n} \otimes \op{X}(m)^{\otimes n}}
            \ar[r]^{\text{shuffle}}
            & {(s^{-1}F \otimes \op{X}(m))^{\otimes n}}
            \ar[d]_{(\theta_{m}^{\perp})^{\otimes n}} \\
        & & {[(\Omega G)^{\otimes m}]^{\otimes n}}
            \ar[d]_{\mu_{(\Omega G)^{\otimes m}}^{(n-1)}} \\
        & &{(\Omega G)^{\otimes m}}
    }
\]
Observer that the iterated multiplication map $\mu _{(\Om G)^{\otimes m}}^{(n-1)}$ performs a transposition on the factors being multiplied before multiplying them, whence the importance of our starting morphism $\theta$ being a transposed tensor morphism.

We note that $\theta_{m}^{\bot}$ is $\Sigma_{m}$-equivariant
because $\theta_{m} : F \otimes \op{X}(m) \rightarrow G^{\otimes
m}$ is. Let $A = \op{A}(s^{-1}F)$.  We construct a family of
morphisms $\theta_{\vec{m}}^{\bot}: A^{\otimes n} \otimes
\op{X}[\vec{m}] \rightarrow (\Omega G)^{\otimes m}$, which are the
composites
\[
    A^{\otimes n} \otimes \op{X}[\vec{m}]
        \xrightarrow{\text{shuffle}}
    \bigotimes_{i=1}^{m} A \otimes \op{X}(m_{i})
        \xrightarrow{\otimes_{i}\theta_{m_{i}}^{\bot}}
    (\Omega G)^{\otimes m}.
\]
Each morphism $\theta_{\vec{m}}^{\bot}$ is
$\Sigma_{n}$-equivariant and $\Sigma_{m_{1}} \times \cdots \times
\Sigma_{m_{n}}$-equivariant, hence the family determines a
morphism of symmetric sequences
\[
    \Ind(\theta) : \op{T}(A) \circ \op{X} \rightarrow
    \op{T}(\Omega G).
\]
By Lemma \ref{lem:sym-seq-multi}, the morphism $\Ind(\theta)$ is multiplicative, since each
$\theta_{m}^{\perp}$ is.  It remains to show that each
$\Ind(\theta)(n)$ commutes with the differentials.  It suffices to
check this on $s^{-1}F \otimes \op{X}(m)$, since
$\op{T}(\op{A}(F))$ is a free left $\op{A}$-module and $\op{A}(F)$
is a free $\op{A}$-algebra.

We may suppose that the
differentials in $F$, $\op{X}(m)$, and $G$ vanish, since the morphism $\Ind(\theta)$  takes into account
all the internal differentials.  Since $\theta$ is a transposed tensor morphism and $\ind \theta$ is multiplicative, the composite $(\ind \theta)\del _{\Om F}$ is determined by its ``linear part,'' which is the desuspension of $\theta \del _{c}$.  On the other hand, the composite $\del_{\Om G}(\ind \theta)$ is determined by $\theta \del _{s}$ because $\theta$ is a morphism of right $\op A$-modules.  It therefore follows from the equality $\theta(\del_{c}+\del _{s)}=0$ that $\ind \theta$ is a chain map.

Next, we define $\Lin$, the linearization transformation.  Let $\theta \in
\msym \big(\op T(\Om F)\circ \op X,\op T(\Om G)\big)$.  Fix $m,n \geq 1$ and $\vec{n} = (n_{1}, \cdots n_{m})\in I_{m,n}$. Let
$\varepsilon_{\vec{n}} : (\op{S} \otimes \op{X})(m) \otimes
\op A^\perp [\vec{n}] \rightarrow \op{X}(m)$ be the map of degree $n-1$
that sends $s_{m-1}x \otimes \alpha_{\vec{n}}$ to $x$. Define
$\theta_{\vec{n}}^{\bot} : F \otimes (\op{S} \otimes \op{X})(m)
\otimes \op A^\perp[\vec{n}] \rightarrow G^{\otimes n}$ to be the
$\Sigma_{n_{1}} \times \cdots \times \Sigma_{n_{m}}$-equivariant
morphism determined by the composite
\[
    \xymatrix{
        {F \otimes(\op{S}(m) \otimes \op{X}(m)) \otimes
        \op A^\perp[\vec{n}]}
            \ar[r]^-{s^{-1} \otimes \varepsilon_{\vec{n}}}
           \ar@{.>}[ddrr]
            & {s^{-1}F \otimes \op{X}(m)}
                \ar[r]^{\theta_{m}}
                & (\Omega G)^{\otimes m}
                    \ar[d]_{\pi[\vec{n}]} \\
            &   & (s^{-1}G)^{\otimes \vec{n}}
                    \ar[d]_{s^{\otimes n}}      \\
            &   & G^{\otimes \vec{n}}
    }
\]
where we use the notation $G^{\otimes\vec{n}} = G^{\otimes n_{1}}
\otimes \cdots G^{\otimes n_{m}}$, and $\pi[\vec{n}]$ is the
obvious projection.  Note that $\theta_{\vec{n}}^{\bot}$ is
homogeneous of degree $0$.  It inherits $\Sigma_{m}$-equivariance
from $\theta$. Let $\theta_{n}^{\bot} = \sum
\theta_{\vec{n}}^{\bot}$, where the sum is taken over all $m \leq
n$ and all $\vec{n} \in I_{n,m}$. Extend the set $\{
\theta_{n}^{\bot} \}$ to a transposed tensor morphism
\[
    \theta^{\bot} : \op{T}(F) \circ (\op{S} \otimes \op{X})
    \circ \op A^\perp \rightarrow \op{T}(G).
\]
Finally, extend $\theta^{\bot}$ to a morphism of right
$\op{A}$-modules,
\[
    \Lin(\theta) : \op{T}(F) \underset{\op{A}}{\circ} \diffract{\op{X}}
    = \op{T}(F) \circ (\op{S} \otimes \op{X})
    \circ \op A^\perp \circ \op{A} \rightarrow \op{T}(G).
\]

It remains to show that each $\Lin(\theta)(n)$ commutes with the
differentials.  As with $\Ind$, we may suppose that the internal
differentials in $F$, $G$, and $\op{X}$ all vanish.  Because $\lin \theta$  commutes with the right action of $\op A$, the composite $(\lin \theta)\del _{s}$ is determined by suspensions of $\del _{\Om G}\theta$.  The  properties of the coequalizer $\acirc$, together with the fact that $\lin \theta$ is a transposed tensor morphism, imply that $(\lin \theta)\del _{c}$ is determined by suspensions of $\theta \del _{\Om F}$.  It therefore follows from  the equality $ \del _{\Om G}\theta =\theta \del _{\Om F}$ that $\lin \theta$ is a chain map.
\end{proof}

The next two results are immediate consequences of the proof the the Cobar Duality Theorem.  Let $\varepsilon:\Phi (\op J)=\op A\circ \op A^\perp\circ \op A\to \op A$ denote the $\op A$-bimodule morphism specified by $\varepsilon(\alpha _{1})=1$, the unit of the ring $\op A(1)$, and $\varepsilon(\alpha _{n})=0$ for all $n\geq 2$.

\begin{schol} \label{schol:elpropind}Induction satisfies the following elementary properties.
\begin{enumerate}
\item For all chain coalgebras $C$, $\ind (Id_{\op T(C)}\acirc \varepsilon)=Id_{\op T(\Om C)}$.
\item For all $\varphi \in \ttma\big(\op T(C)\acirc \Phi(\op X), \op T(C')\big)$ and for all $\beta\in \cat {Csg}_{\otimes}(\op Y,\op X)$,
$$\ind\Big(\varphi\big(Id_{\op T(C')}\acirc \Phi(\beta)\big)\Big)=\ind \varphi(Id_{\op T(\Om C')}\circ \beta).$$
\end{enumerate}
\end{schol}

\begin{schol}
\label{schol:elproplin} Linearization satisfies the following elementary properties.
\begin{enumerate}
\item For all chain coalgebras $C$, $\lin (Id_{\op T(\Om C)})=Id_{\op T(C)}\circ \varepsilon$.
\item For all $\psi\in \msym(\op T(\Om C)\circ \op X, \op T(\Om C'))$ and for all $\beta\in \cat {Csg}_{\otimes}(\op Y,\op X)$,
$$\lin\Big(\psi\big(Id_{\op T(\Om C')}\circ \beta)\big)\Big)=\lin \psi\big(Id_{\op T(C')}\acirc \Phi(\beta)\big).$$
\end{enumerate}
\end{schol}

\subsection{Bar duality}\label{ssec:bar-duality}

In this section we dualize the Cobar Duality Theorem, essentially by simply
considering left modules rather than right modules.

\begin{defn}\label{defn:comultiplicative}
Suppose that $\op X$, $\op{M}$ and $\op{N}$ are level comonoids, with
diagonals $\Delta_{\op{X}}$,  $\Delta_{\op{M}}$ and $\Delta_{\op{N}}$, respectively.
We say that $\theta:\op X\circ\op M\to \op N$ is \emph{comultiplicative} if the diagram
\[
    \xymatrix{
        {\op{X}\circ \op{M}}
            \ar[r]^{\theta}
            \ar[d]^{\Delta_{\op{X}} \circ \Delta_{\op{M}}}
            & {\op{N}}
                \ar[dd]^{\Delta_{\op{N}}} \\
        {\op{X}^{\otimes 2} \circ \op{M}^{\otimes 2}}
            \ar[d] ^\iota\\
        {(\op{X}\circ \op{M})^{\otimes 2}}
            \ar[r]^-{\theta^{\otimes 2}}
            & {\op{N}^{\otimes 2}}
    }
\]
commutes.
\end{defn}

In particular, if for each $m$ there exists a free graded $R$-module $V_{m}$ such that $\op{N}(m)$ is a cofree coalgebra on $V_{m}$, then a family of equivariant morphisms $\{ \theta_{m} :
(\op{X} \circ \op{M})(m) \rightarrow V_{m} \}$ lifts to a
comultiplicative morphism $\hat{\theta} : \op{X} \circ \op{M}
\rightarrow \op{N}$.

Let $\comsym (\op X\circ\op M,  \op N)$ denote the set of comultiplicative morphisms of symmetric sequence, when $\op X$, $\op M$ and $\op N$ are level comonoids.

Recall that we have defined the associative operad $\op A$ so that $\op A(0)=0$.  As a result, $\op A$-algebras are nonunital, i.e., an $\op A$-algebra consists of a chain complex $A$ together with an associative chain map $\mu:A\otimes A\to A$.

The next definition generalizes slightly the usual definition of the bar construction.

\begin{defn}  The \emph{bar construction}  is the functor $B: \alg{A}\to \coalg A$ defined as follows. For all $A\in \ob \alg{A}$, the graded $\op A$-coalgebra underlying $BA$ is the cofree $\op A$-coalgebra $\op A^{\sharp}(s C)$, i.e., $\prod_{n\geq 1}\big(s A\big)^{\otimes n}$ with the obvious coproduct, where $s A$ is the suspension of $A$, i.e., $(sA )_{n}=A_{n-1}$.  Let $\pi_{k}: \op A^\sharp(sA)\to (sA)^{\otimes k}$ denote the projection.  The differential $d_{B }$ on $BA$ is the coderivation specified by
$$s^{-1}\pi d_{B}=-ds^{-1}\pi_{1} +\mu(s^{-1}\otimes s^{-1})\pi_{2}:\op A^\sharp(sA)\to A,$$
where $d$  and $\mu$ are the differential and product of $A$.
\end{defn}

Analogously to the the case of the cobar construction, if $A$ is actually unital and augmented, then one can either keep or throw away the unit, without changing the homology of the bar construction on $A$.

\begin{BDT}
Let $\cat D$ be any small category.
There are mutually inverse, natural isomorphisms
$$\ind: \am\big(\Phi(-)\acirc c(-),c(-)\big)\longrightarrow \comsym \big(-\circ \,c(B-),c(B-)\big),$$
and
$$\lin :\comsym \big(-\circ\, c(B-),c(B-)\big)\longrightarrow\am\big(\Phi(-)\acirc c(-),c(-)\big)$$
of  functors from $(\alg{A})^{\cat D}\times\cat {Comon}_{\otimes}\times (\alg{A})^{\cat D}$ to  $\cat{Set}^{\cat D}.$
\end{BDT}

 As in the dual case, the natural isomorphisms $\ind$ and $\lin$ are called \emph{induction} and \emph{linearization}, respectively.

\begin{proof}
The proof follows the lines of the Cobar Duality Theorem.

Let $F,G\in \ob (\alg{A})^{\cat D}$, and let $\op X\in \ob\cat {Comon}_{\otimes}$.  First we
define $\Ind$.  Consider $\theta \in
\am\big(\Phi(\op X)\acirc c(F),c(G)\big)$, and fix $m \geq 1$ and $\vec{n}
\in I_{m,n}$.  Define $\theta^{\bot}_{\vec{n}} : \op{X}(m) \otimes
(BF)^{\otimes m} \rightarrow sG$ to be the composite
\[
    \xymatrix{
        {\op{X}(m) \otimes (BF)^{\otimes m}}
            \ar[rr]^{1 \otimes \pi[\vec{n}]}
            \ar@{.>}[d]^{\theta_{\vec{n}}^{\perp}}
        & & {\op{X}(m) \otimes (sF)^{\otimes \vec{n}}}
            \ar[d]^{\eta_{\vec{n}} \otimes (s^{-1})^{\otimes n}}
            \\
        {sG}
        & {G} \ar[l]^{s}
        & {(\op{S} \otimes \op{X})(m) \otimes \op A^\perp[\vec{n}]
        \otimes F^{\otimes n}} \ar[l]^-{\theta_{n}}
    }
\]
where $\pi[\vec{n}]$ is the projection and $\eta_{\vec{n}}$ is the
obvious inclusion.  The morphism $\theta_{\vec{n}}^{\bot}$ is
easily seen to be homogeneous of degree zero, and appropriately
equivariant, since $\theta_{n}$ is.  Let $\theta_{m}^{\bot} =
\sum_{\vec{n} \in I_{m,*}} \theta_{\vec{n}}^{\bot}$ and lift
comultiplicatively to
\[
    \hat{\theta}_{m}^{\bot} : \op{X}(m) \otimes (BF)^{\otimes m}
    \rightarrow BG.
\]
The sequence $\hat{\theta}_{m}^{\bot}$ is equivalent to a
comultiplicative morphism of right $\op{A}$-modules,
\[
    \Ind(\theta) : \op{X} \circ c(BF) \rightarrow c(BG).
\]
To show that $\Ind(\theta)$ commutes with the differentials, we
may assume without loss of generality that all internal
differentials vanish. The result then follows by direct
calculation, using the left $\op{A}$-linearity of $\theta$.

Now we define $\Lin$.  Consider $\theta \in
\comsym \big(\op X\circ\, c(BF),c(BG)\big)$.  Fix $m \geq 1$ and $\vec{n} \in
I_{m,n}$.  Define $\theta_{\vec{n}}^{\bot} : (\op{S} \otimes
\op{X})(m) \otimes \op A^\perp[\vec{n}] \otimes F^{\otimes n}
\rightarrow G$ to be the composite
\[
    \xymatrix{
        {(\op{S} \otimes \op{X})(m) \otimes \op A^\perp[\vec{n}]
        \otimes F^{\otimes n}}
        \ar[rr]^{\varepsilon_{\vec{n}} \otimes s^{\otimes n}}
        \ar@{.>}[dd]_{\theta_{\vec{n}}^{\bot}}
        & & {\op{X}(m) \otimes (sF)^{\otimes \vec{n}}}
            \ar[d] \\
        & & {\op{X}(m) \otimes (BF)^{\otimes m}}
            \ar[d]^{\theta_{m}} \\
        {G}
        & {sG} \ar[l]^{s^{-1}}
        & {BG} \ar@{->>}[l]
    }
\]
where $\varepsilon_{\vec{n}}$ is the evident projection.  It is
clear that $\theta_{\vec{n}}^{\bot}$ is homogeneous of degree
zero.  Let $\theta_{n}^{\bot} = \sum_{\vec{n} \in I_{*,n}}
\theta_{\vec{n}}^{\bot}$.  After verifying equivariance, we may
extend the sequence $\{ \theta_{n}^{\bot} \}$ to a morphism of
left $\op{A}$-modules,
\[
    \lin(\theta) : \op{A} \circ \left( (\op{S} \otimes \op{X}) \circ \op{A}^\perp
    \circ c(F) \right) \rightarrow c(G).
\]
The usual arguments now show that $\lin(\theta)$ commutes with the
differentials.
\end{proof}

\begin{rmk}
Of course, there also exists the evident version of the Bar
Duality Theorem for the embedding $z$ of $\op{A}$-algebras as left
$\op{A}$-modules given in Proposition~\ref{prop:z}.
\end{rmk}

\subsection{Existence of
$\op{M}$-governed morphisms}\label{ssec:acyclic}

The existence results stated here and in section 7 are expressed in terms of acyclic models. We recall the foundations of this method before stating the existence theorems.  Let $\cat D$ be a  category, and let $\mathfrak M$ be a set of objects in $\cat D$.    A functor $X: \cat D\to\cat {Ch}$ is \emph{free} with respect to $\mathfrak M$ if there is a set $\{x_{\mathfrak m}\in X(\mathfrak m)\mid \mathfrak m\in \mathfrak M\}$ such that $\{X(f)(x_{\mathfrak m})\mid f\in \cat D(\mathfrak m,d), \mathfrak m\in \mathfrak M\}$ is an $\mathbb Z$-basis of $X(d)$ for all objects $d$ in $\cat D$. The functor $X$ is \emph{acyclic}  with respect to $\mathfrak M$ if $X(\mathfrak m)$ is acyclic for all $\mathfrak m\in \mathfrak M$.   More generally, if $\cat C$ is a category with a forgetful functor $U$ to $\cat {Ch}$ and $X:\cat D\to \cat C$ is a functor, we say that $X$ is \emph {free}, respectively \emph{acyclic},  with respect to $\mathfrak M$ if $UX$ is.

We now restate Theorems~\ref{thm:exist-1} and~\ref{thm:exist-2} from
the introduction, which are special cases of more highly structured existence results results proved in section 7.  These existence theorems have already played an essential role in \cite {hpst:04}, \cite {hl}, and \cite {hess}.

We work here in the comma category (or under category) $\op J\downarrow \cat {Ch}^\Sigma$.  The objects of $\op J\downarrow \cat {Ch}^\Sigma$ are morphisms of symmetric sequences $\op J\to \op X$ and are called \emph{symmetric sequences under $\op J$}.  The morphisms are commuting triangles.   A level comonoid \emph {under $\op J$} is a morphism of level comonoids $\op J\to \op X$.

Recall that $\op J(n)=0$ for $n=0$ and for all $n\geq 2$ and that $\op J(1)=R$.   The canonical isomorphisms $\op J(1)\cong \op J(1)\otimes \op J(1)$ and $0\otimes 0\cong 0$ define a level comonoidal structure on $\op J$. Given any symmetric sequence $j_{\op X}:\op J\to \op X$ under $\op J$, we let $e_{\op X}=j_{\op X}(1)$.

A functor $X: \cat D\to \cat {Ch}$, where $\cat D$ is a small category, is \emph{globally connective} if there is an integer $N$ such that $X(d)_{n}=0$ for all $n<N$ and for all $d\in \ob \cat D$. The integer $N$ is called a \emph{global lower bound} on the functor $X$.  In particular, a symmetric sequence $\op X$ is globally connective if there is an integer $N$ such that $\op X(m)_{n}=0$ for all $n<N$ and for all $m\geq 0$.

\begin{thm}\label{thm:x-exist-cobar}
Let $X,Y: \cat D\to\coalg A$ be functors, where $\cat D$ is a category admitting a set of models $\frak M$ with respect to which  $X$ is free and globally connective and $Y$ is acyclic.  Let $(\op M, \Delta)$ be a level comonoid under $\op J$.  Let $\tau:UX\to UY$ be a natural transformation, where $U $ is the forgetful functor down to $\cat {Ch}$.  Then there is a natural, multiplicative  transformation
$$\theta:\op T(\Om X)\circ \op M\to \op T(\Om Y)$$
extending $s^{-1}\tau$,  i.e., the following composite is equal to $\tau$.
$$X\xrightarrow{s^{-1}}s^{-1}X\hookrightarrow \op T(\Om X)\circ \op J\to \op T(\Om X)\circ \op M\xrightarrow{\theta}\op T(\Om Y)\xrightarrow {\text{proj.}}s^{-1}Y\xrightarrow {s}Y$$
\end{thm}

Theorem \ref {thm:x-exist-cobar} is the special case $\op P=\op J$ of Theorem \ref{thm:diffmod}, which is proved  in the last section of this article.

The next theorem, dual to Theorem \ref {thm:x-exist-cobar}, has an analogous, dual proof.

\begin{thm}
Let $X,Y: \cat D\to\alg A$ be functors, where $\cat D$ is a category admitting a set of models $\frak M$ with respect to which  $X$ is free and globally connective and $Y$ is acyclic.  Let $(\op M, \Delta)$ be a level comonoid under $\op J$.  Let $\tau:UX\to UY$ be a natural transformation, where $U $ is the forgetful functor down to $\cat {Ch}$.  Then there is a natural, comultiplicative  transformation
$$\theta:\op T(B X)\circ \op M\to \op T(B Y)$$
lifting $s\tau$,  i.e., the following composite is equal to $\tau$.
$$X\xrightarrow{s}sX\hookrightarrow \op c(B X)\circ \op J\to c(BX)\circ \op M\xrightarrow{\theta}c(BY)\xrightarrow {\text{proj.}}sY\xrightarrow {s^{-1}}Y$$
\end{thm}

Our need to establish such existence results is the true raison-d'\^etre of the Cobar and Bar Duality Theorems.  It is not possible to apply acyclic models methods directly to proving the existence of natural, multiplicative maps $\op T(\Om F)\circ \op X\to \op T(\Om G)$, since $\Om F$ is too big to be free on the models with respect to which $F$ is free.  Thanks to  cobar duality, it is enough to show that there is a natural transposed tensor transformation $\op T(F)\acirc \Phi(\op X)\to \op T(G)$, which can be done via acyclic models, then to apply induction to obtain the desired, natural,  multiplicative map $\op T(\Om F)\circ \op X\to \op T(\Om G)$.

\section{The Alexander-Whitney co-ring}\label{sec:aw}

Consider the category
$\cat{DCSH}$, in which the objects are coassociative
chain coalgebras, $C$, and where $\cat{DCSH}(C,C')=\alg{A}(\Om C, \Om C')$.  A chain map $f:C\to C'$ is then called \emph{strongly homotopy comultiplicative} or a \emph{DCSH-map} if there exists $\varphi\in\alg{A}(\Om C, \Om C')$ such that the  composite
$$C\xrightarrow{s^{-1}}s^{-1} C\hookrightarrow \Om C\xrightarrow{\varphi}\Om C'\xrightarrow {\pi}s^{-1}C'\xrightarrow {s}C$$
is exactly $f$, where $\pi$ is the projection.  The category $\cat {DASH}$ and \emph{strongly homotopy multiplicative maps} or \emph{DASH-maps} are defined dually.

Gugenheim and Munkholm also showed that the Alexander-Whitney equivalence
\[
    C_{*}(K \times L) \rightarrow C_{*}(K)\otimes C_{*}(L)
\]
of normalized chains on reduced simplicial sets is naturally a
$DCSH$-map, which implies, as shown in \cite {hpst:04}, that $\Delta_{C(K)}$ is naturally
a $DCSH$-map.

Recall that ${\op{J}}$ is the unit symmetric sequence with respect to the composition product.  In fact,
$\op{J}$ is a non-unital operad, whose algebras are simply chain
complexes.  Note that ${\op{J}}(1)$ is a chain coalgebra $(R\{
u_{0} \},\Delta,\varepsilon)$ with $\Delta u_{0} = u_{0} \otimes
u_{0}$ and $\varepsilon(u_{0})=1$. Let $\op{F} =
\diffract(\op{J})=(\op A\circ \op A^\perp\circ \op A, \del _{\op F})$. Then
\[
    \op{F} = \{ \op{F}(m) \mid m \in \N \}
\]
with generators
\[
    \{ f_{m} = s^{m-1}u_{0} \in \op{F}(m)_{m-1} \mid m \in \N \}
\]
satisfying
\[
    \partial f_{m} = \sum_{i=1}^{m-1} \delta \otimes (f_{i} \otimes
    f_{m-i}) + \sum_{i=0}^{m-2} f_{m-1} \otimes (1^{\otimes i} \otimes
    \delta \otimes 1^{\otimes (m-2-i)}).
\]

By Proposition~\ref{prop:mho-comonoidal}, ${\op{F}}$ is an $\op{A}$-co-ring. Indeed, the formula for the composition comultiplication $\psi_{\op F}$
is particularly simple.  For $n \geq 1$, we have
\[
    \psi_{\op F} f_{n} = \sum_{m \leq n} \sum_{\vec{n} \in I_{m,n}}
    f_{m}
    \otimes ( f_{n_{1}} \otimes \cdots \otimes f_{n_{m}} )
\]
where $n_{i} \geq 1$ for all $i$.  In fact $\op F$ is a counital $\op{A}$-co-ring, with counit $\varepsilon:\op F\to \op A$ specified by $\varepsilon(f_{n})=0$ for all $n>1$ and $\varepsilon (f_{1})=1$.

Since $\Delta : \op{J} \rightarrow \op{J} \otimes \op{J}$ is
cocommutative, it is a morphism of comonoids.  Thus we may compose
\[
    \diffract(\op{J}) \xrightarrow{\diffract(\Delta)}
    \diffract(\op{J}\otimes \op{J}) \xrightarrow{q}
    \diffract(\op{J})\otimes \diffract(\op{J}).
\]
Let $\Delta_{\op{F}}=q \circ \diffract(\Delta_{\op{J}})$.  Then
$(\op{F},\partial_{\op{F}},\Delta_{\op{F}})$ is a level comonoid in the category
$\op{A}$-bimodules. Explicitly,
\begin{eqnarray*}
    \Delta_{\op{F}}(f_{m})
    & = & \sum_{k=1}^{m}\sum_{\vec{\imath} \in I_{k,m}}
    \left(f_{k} \otimes \delta^{(i_{i})}\otimes \cdots \otimes
    \delta^{(i_{k})}\right)
    \otimes
    \left( \delta^{(k)} \otimes f_{i_{1}} \otimes \cdots
    \otimes f_{i_{k}}\right).
\end{eqnarray*}
It is easy to check that $\Delta _{\op F}$ is coassociative.

The Bar Duality Theorem, applied to the comonoid $\op{J}$,
states that a morphism of coalgebras $BA \rightarrow BA'$ is
equivalent to a morphism of left $\op{A}$-modules $\op{F}
\circ_{\op{A}} c(A) \rightarrow c(A')$.

Let $C$ and $C'$ be connected, counital, coassociative chain coalgebras. The Cobar Duality Theorem
tells us that there are mutually inverse isomorphisms, natural in $C$ and $C'$,
$$\ind: \ttma\big(\op T( C)\acirc \op F, \op T(C')\big)\longrightarrow \msym \big(\op T(\Om C),\op T(\Om C')\big)$$
and
$$\lin :\msym \big(\op T(\Om C'),\op T(\Om C')\big)\longrightarrow\ttma\big(\op T(C)\acirc \op F,\op T(C')\big).$$

Similarly, let $A$ and $A'$ be connected, augmented, associative chain algebras. The Bar Duality Theorem says that there are mutually inverse isomorphisms, natural in $A$ and $A'$,
$$\ind: \am\big(\op F\acirc c(A),c(A')\big)\longrightarrow \comsym \big(c(BA),c(BA')\big),$$
and
$$\lin :\comsym \big( c(BA),c(BA')\big)\longrightarrow\am\big(\op F\acirc c(A),c(A')\big).$$

The next theorem is an immediate consequence of these observations.

\begin{thm}\label{thm:dcsh-is-fat} There are isomorphisms of categories
$$\cat{DCSH}\cong \Coalg{A}{\op F}\qquad\text{and}\qquad \cat{DASH}\cong \Alg{A}{\op F}.$$
\end{thm}

Either as a direct consequence of Theorem \ref{thm:dcsh-is-fat} together with an acyclic models argument or as Theorem \ref {thm:x-exist-cobar} in the special case $\op M=\op J$, we obtain the following existence result.

\begin{thm}\label{thm:dcsh-exist}
Let $X, Y: \cat D\to\coalg A$  be functors, where $\cat D$ is a category admitting a set of models $\frak M$ with respect to which  $X$ is free and globally connective and $Y$ is acyclic.  Let $\tau:UX\to UY$ be a natural transformation, where $U $ is the forgetful functor down to $\cat {Ch}$.  Then there is a natural, multiplicative  transformation
$$\theta:\Om X\to \Om Y$$
extending $s^{-1}\tau$,  i.e., the following composite is equal to $\tau$.
$$X\xrightarrow{s^{-1}}s^{-1}X\hookrightarrow \Om  X\xrightarrow{\theta}\Om  Y\xrightarrow {\text{proj.}}s^{-1}Y\xrightarrow {s}Y$$
In other words, for each $d\in \ob \cat D$, the natural chain map $\tau(d):F(d)\to G(d)$ admits a natural DCSH-structure.
\end{thm}

Dually, for algebras we have the existence result below.

\begin{thm}\label{thm:dash-exist}
Let $X, Y: \cat D\to\alg A$  be functors, where $\cat D$ is a category admitting a set of models $\frak M$ with respect to which  $X$ is free and globally connective and $Y$ is acyclic.  Let $\tau:UX\to UY$ be a natural transformation, where $U $ is the forgetful functor down to $\cat {Ch}$.  Then there is a natural, comultiplicative  transformation
$$\theta:BX\to BY$$
extending $s\tau$,  i.e., the following composite is equal to $\tau$.
$$X\xrightarrow{s}sX\hookrightarrow BX\xrightarrow{\theta}BY\xrightarrow {\text{proj.}}sY\xrightarrow {s^{-1}}Y$$
In other words, for each $d\in \ob \cat D$, the natural chain map $\tau(d):F(d)\to G(d)$ admits a natural DASH-structure.
\end{thm}

Since the homology of a strongly homotopy-(co)multiplicative map
is a (co)algebra morphism, it stands to reason that we should
expect $H_{*}(\op{F}) \cong \op{A}$, and this is indeed the case.

\begin{thm}\label{thm:HF}
The counit $\varepsilon : \op{F} \rightarrow \op{A}$ is a
quasi-isomorphism in positive levels.
\end{thm}

\begin{proof}
We consider the first-quadrant spectral sequence  associated to the decreasing filtration
\[
    F_{\ell} = \bigoplus_{m\geq\ell}
        \op{A}(m) \otimes_{\Sigma_{k}}(\op A^\perp \circ \op A)^{\odot m}.
\]

From the definitions, we have $E^{0}_{0,*}=0$, and $E^{0}_{1,*} \cong \Ssigma \circ \op{A}$. Since the
cosimplicial differential raises filtration, the differential
$d^{0}$ in the $E^0$ term is given by the simplicial differential in $\Ssigma \circ
\op{A}$.   As remarked in section \ref{sec:delooping}, $\op A^\perp\circ \op A$ is the bar resolution of $(\mathbf N, +)$, which implies that $E^1_{1,*}=R\{\alpha _{1}\}$.

Let $V$ be a chain complex over $R[G]$, let $G \rightarrow G'$ be a
homomorphism.  Then $H_{*}(V) \otimes_{G} R[G'] \xrightarrow{\cong} H_{*}(V
\otimes_{G} R[G'])$, and so it follows from the K\"unneth theorem that
\begin{eqnarray*}
    E^{1}_{k,*}(n) & \cong &  \op{A}(k) \otimes_{\Sigma_{k}} H_{*}(\Ssigma \circ
    \op{A})^{\odot k}(n) \\
        & \cong & \left\{
            \begin{array}{cl}
                R[\Sigma_{n}]\{ \alpha_{1}^{\otimes n} \}   & k = n \\
                0                                       & k \neq
                n.
            \end{array}
        \right.
\end{eqnarray*}
Therefore the spectral sequence collapses at the $E_{1}$ term, and
so $H(\op{F})(n) = \op{A}(n)$,  for $n
\geq 1$.
\end{proof}

\begin{rmk}  Careful inspection of the definitions leads to the conclusion that the Alexander-Whitney bimodule is exactly the two-sided Koszul resolution of $\op A$.  The calculation above therefore provides another proof that $\op A$ is a Koszul operad \cite{ginzburg-kapranov}.
\end{rmk}

\begin{rmk} Since the bimodule ${\op{F}}$ is a free
${\op{A}}$-bimodule resolution of ${\op{A}}$, we
may use it to do homological algebra. If $\op{M}$ and $\op{N}$ are
right and left ${\op{A}}$-modules, respectively, then
\[
    \mathrm{Tor}^{{\op{A}}}(\op{M},\op{N})
        := H( \op{M} \underset{{\op{A}}}{\circ}
        {\op{F}} \underset{{\op{A}}}{\circ}
        \op{N}).
\]
In particular, we may read off the isomorphisms
\[
    \mathrm{Tor}^{{\op{A}}}(\op{J},\op{J}) \cong \Ssigma
    \quad \text{and} \quad
    \mathrm{Tor}^{{\op{A}}}(\op{J},{\op{A}}) \cong
    {\op{J}}.
\]
\end{rmk}

\begin{ex}\label{ex:counter}
We are now in a position to construct over $R=\Zmod{2}$ an example
of a chain coalgebra $M$ such that
\begin{itemize}
\item  its cohomology
algebra is realizable, i.e., there is a topological space $X$ such that $H^*(X;\mathbb F_{2})\cong H^*M$ as graded algebras, but
\item  $M$ is not
quasi-isomorphic to the chains on any space.
\end{itemize}
This example,
along with~\cite[Example 3.8]{ndombol-thomas:04}, show that the
concepts of ``\emph{shc} algebras'' and ``algebras with cup-$i$
products'' are independent of one another.

Let $M = \Zmod{2}\{ 1_{0}, u_{2}, x_{3}, y_{3}, z_{3}, v_{4},
w_{6} \}$, where subscript indicates degree.  The only non-zero
differential in $M$ is $\partial(v)=x+y$.  All elements other than
$v$ and $w$ are primitive, while $\bar{\psi}(v) = u \otimes u$ and
$\bar{\psi}(w) = x \otimes z + z \otimes y$.  It is readily
verified that $(M,\partial,\psi)$ is a coassociative chain
coalgebra.

Let $W$ be the usual $\Zmod{2}[\Sigma_{2}]$-free resolution of
$\Zmod{2}$.  Specifically, $W_{i}$ is generated by an element
$e_{i}$ with $\partial e_{i} = (1+\tau)e_{i-1}$, where $\tau \in
\Sigma_{2}$ is the transposition.

\begin{prop}
There exists an equivariant morphism
\[
    \theta : W \otimes M \rightarrow M \otimes M
\]
such that $\theta ( e_{0} \otimes - ) = \psi$.
\end{prop}

\begin{proof}
We construct the morphism $\theta$; verification that it is a
chain map is routine and left to the reader.

It suffices to define $\theta$ on generators.  The only non-zero
values that $\theta$ takes on generators are $\theta(e_{1} \otimes
w) = v \otimes z + z \otimes v$, $\theta(e_{3} \otimes v) = v
\otimes x + y \otimes v$, and $\theta(e_{|a|} \otimes a) = a
\otimes a$ for $a \in \{ u, v, w, x, y, z \}$.
\end{proof}

By~\cite{may:70}, $\theta$ defines cup-$i$ products in the $\Zmod{2}$-dual
$M^{\sharp}$, and so $H^{*}(M,\Zmod{2})$ comes equipped with an
action of the mod $2$ Steenrod algebra.  In fact, we have the
following proposition.

\begin{prop}The algebra $H^{*}(M;\Zmod{2})$ admits the structure of an
unstable algebra over the mod $2$ Steenrod algebra, where the only
non-trivial operation is the $\mathrm{Sq}^0$. Moreover, this
algebra is isomorphic to
\[
    H^*(S^2\vee(S^3\times S^3);\Zmod{2})
\]
as unstable algebras over the mod $2$ Steenrod algebra.
\end{prop}

\begin{proof}
An easy exercise in $\Zmod{2}$-linear algebra.
\end{proof}

\begin{prop}\label{prop:example}
The chain coalgebra $M$ is not realizable, i.e., $M$ is not of the
same homotopy type as $C_{*}(X;\Zmod{2})$ for any space $X$.
\end{prop}

\begin{proof}
For the duration of the proof, we suppress the coefficients from
the notation. By~\cite{gugenheim-munkholm:74}, if $X$ is a space,
then the diagonal on $C_{*}(X)$ is strongly
homotopy-comultiplicative, that is, there is a morphism of
symmetric sequences, $\Delta : \op{T}(C_{*}(X)) \circ \op{F}
\rightarrow \op{T}(C_{*}(X) \otimes C_{*}(X))$, such that
$\Delta_{1}$ is the diagonal.  If $C_{*}(X)$ and $M$ are connected
by a sequence of chain coalgebra quasi-isomorphisms, then we may
construct a morphism $\Psi : \op{T}(M) \circ \op{F} \rightarrow
\op{T}(M \otimes M)$ such that $\Psi_{1} = \psi$.  The homotopy
class of such a $\Psi$, compatible with $\Delta$, is unique. We
show that no such $\Psi$ exists.

We attempt to define $\Psi$ on generators $a \otimes f_{k}$, for
$a \in M$ and $k \geq 1$.  Necessarily, $\Psi(a \otimes f_{1}) =
\psi(a)$.  We may define
\[
    \Psi(w \otimes f_{2}) = (1 \otimes v) \otimes (z \otimes 1)
    + (1 \otimes z) \otimes (v \otimes 1)
\]
and $\Psi(a \otimes f_{2}) = 0$ for $a \neq w$.  Any other choice
of morphism $\Psi':M \otimes \op{F}(2) \rightarrow (M \otimes
M)^{\otimes 2}$ is necessarily homotopic to $\Psi$.

Now we try to define $\Psi$ on $M \otimes \op{F}(3)$.  In order to
find a value for $\Psi(w \otimes f_{3})$, we must find an element
that bounds
\[
    (1 \otimes z) \otimes (u \otimes 1) \otimes (u \otimes 1)
    + (1 \otimes u) \otimes (1 \otimes u) \otimes (z \otimes 1),
\]
but no such element exists. Therefore the diagonal on $M$ does not
extend to an $\op{F}$-governed morphism.
\end{proof}

\end{ex}

\section{Enriched cobar duality}\label{sec:enrcobdual}

We prove here an enriched version of the Cobar Duality Theorem, expressing induction and linearization as mutually inverse isomorphisms of \emph {bundles of bicategories with connections}. This seems to be most appropriate language in which to express the high degree of naturality and of compatibility with monoidal struture present in the induction/linearization equivalence.

We begin by defining what it means for a bundle of bicategories to admit a connection and recalling the elements of the theory from \cite {hess:bicatbndl}.   We then introduce two seemingly different bundles with connection over the  monoidal category of level comonoids. The point of Enriched Cobar Duality is to prove that these two bundles are in fact isomorphic.  As a consequence of Enriched Cobar Duality, we obtain in section  \ref {sec:enrind} that the induction functor can itself be enriched.

Recall that a (small) \emph{bicategory} $\bicat A$  consists of
\begin{enumerate}
\item a set $\bicat A_{0}$ of \emph{0-cells};
\item a small category $\bicat A_{1}(a,b)$ for every $a,b\in \bicat A_{0}$,  whose objects are \emph{1-cells} $\varphi: a\to b$  and whose morphisms are \emph{2-cells} $\alpha :\varphi\Rightarrow \varphi '$, where composition in $\bicat A_{1}(a,b)$ of two $2$-cells $\alpha :\varphi\Rightarrow \varphi '$ and $\alpha' :\varphi'\Rightarrow \varphi ''$, called the \emph{vertical composition}, is denoted $\alpha'*\alpha: \varphi\Rightarrow \varphi''$;
\item a composition functor for every $a,b,c\in \bicat A_{0}$
$$\bicat A_{1}(b,c)\times \bicat A_{1}(a,b)\longrightarrow \bicat A_{1}(a,c):(\varphi,\psi)\mapsto \psi\cdot \varphi,$$
 sending a pair of $2$-cells $\alpha :\varphi\Rightarrow \varphi '$ and $\beta :\psi\Rightarrow \psi '$ to their \emph{horizontal composite} $\beta \cdot \alpha:\psi\cdot \varphi\Rightarrow \psi'\cdot\varphi'$;
 \item identity $1$-cells $Id_{a}: a\to a$  and identity $2$-cells $i_{a}:Id_{a}\Rightarrow Id_{a}$ for all $a\in \bicat A_{0}$.
\end{enumerate}
Furthermore, these data must satisfy the following axioms.
\begin{enumerate}
\item For all $a,b,c,d\in \bicat A_{0}$ and for all $\varphi\in \bicat A_{1}(a,b)$, $\psi\in \bicat A_{1}(b,c)$ and $\omega\in \bicat A_{1}(c,d)$, there is a natural choice of $2$-cell isomorphism
$$\alpha _{abcd}:\omega\cdot(\psi\cdot \varphi)\overset \cong{\Longrightarrow}(\omega\cdot \psi)\cdot \varphi.$$
\item For all $a,b,\in A_{0}$ and for all $\varphi\in\bicat A_{1}(a,b)$ there are two natural $2$-cell isomorphisms
$$\lambda_{ab}:\varphi\cdot Id_{a}\overset\cong\Longrightarrow \varphi\quad\text{and}\quad \rho_{ab}:Id_{b}\cdot\varphi\overset\cong\Rightarrow \varphi.$$
\item All appropriate coherence diagrams commute.
\end{enumerate}

The last axiom is undeniably vague.  The interested reader is encouraged to consult, e.g., \cite{borceux}, for precise details.

Given a bicategory $\bicat A$, we let $\bicat A_{1}=\coprod _{a,b\in \bicat A_{0}}\bicat A_{1}(a,b)$, the set of all $1$-cells.  For any $\varphi, \varphi'\in \bicat A_{1}(a,b)$, we let $\bicat A_{2}(\varphi, \varphi')$ denote the set of all $2$-cells from $\varphi$ to $\varphi'$.  Finally, we let
$$\bicat A_{2}=\coprod _{a,b\in \bicat A_{0}}\;\coprod _{\varphi, \varphi' \in \bicat A_{1}(a,b)}\bicat A_{2}(\varphi,\varphi'),$$
the set of all $2$-cells in $\bicat A$.

A \emph{bicategory homomorphism} $\Xi:\bicat A\to \bicat B$, consists of three functions $\Xi_{k}:\bicat A_{k}\to \bicat B_{k}$ respecting all identities and compositions \emph{strictly}.  One can also reasonably define notions of morphisms between bicategories in which identities are strictly respected but in which composition is respected up to a $2$-cell or up to an invertible $2$-cell, but we have no need to work in this generality here.

All proofs concerning the general theory of bicategory bundles with connection can be found in \cite {hess:bicatbndl}.

\subsection{Bicategory bundles with connection}

Let $\bicat B$ be a bicategory with exactly one 0-cell, denoted $*$, and let $\bicat E$ be any bicategory.  Let $\Pi: \bicat E\to \bicat B$ be a bundle in the category of bicategories, i.e., a bicategory homomorphism.  Observe that for all $e,e'\in \bicat E_{0}$ there is a splitting of categories
$$\bicat E_{1}(e,e')=\coprod _{b\in B_{1}}\bicat E_{1}(e,e')_{b},$$
where $\bicat E_{1}(e,e')_{b}=\Pi_{1} ^{-1}(b)\cap \bicat E_{1}(e,e')$, and that for all $\varphi\in \bicat E_{1}(e,e')_{b}$, $\psi\in \bicat E_{1}(e,e')_{b'}$
$$\bicat E_{2}(\varphi, \psi)=\coprod _{\alpha\in \bicat B_{2}(b,b')} \bicat E_{2}(\varphi, \psi)_{\alpha},$$
where $\bicat E_{2}(\varphi, \psi)_{\alpha}=\Pi _{2}^{-1}(\alpha)\cap \bicat E_{2}(\varphi, \psi)$.
In terms of this decomposition,
$$\varphi\in\bicat E_{1}(e,e')_{b}, \psi \in\bicat E_{1}(e',e'')_{b'}\Rightarrow \psi\cdot\varphi\in\bicat E_{1}(e,e'')_{b'\cdot b}.$$

The bundles we study here satisfy the additional criteria specified in the definition below.  Recall that a bicategory $\bicat B$ with exactly one object can also be seen as a monoidal category, which we denote $(\cat B, \otimes)$, where $\ob\cat B=\bicat B_{1}$, $\mor\cat B=\bicat B_{2}$ and $b'\otimes b$ is the same as the composite $b'\cdot b$ of $b$ and $b'$ in $\bicat B_{1}$.

\begin{defn}\label{defi:conn}
An \emph{op-connection} on a bicategory bundle $\Pi:\bicat E\to \bicat B$, where $\bicat B$ has exactly one object, consists of a family of functors
$$\boldsymbol\nabla=\{ \nabla _{e,e'}: \cat B^{op}\to \cat {Set}\mid e,e'\in \bicat E_{0}\}$$
that is natural in $e$ and $e'$ and such that
\begin{enumerate}
\item $\nabla _{e,e'}(b)=\bicat E_{1}(e,e')_{b}$ for all $b\in \ob\cat B$, and therefore, for all $\alpha\in \cat B(b',b)=\cat B^{op}(b,b')$, there is a \emph {parallel transport} function
$$\nabla _{e,e'}(\alpha): \bicat E_{1}(e,e')_{b}\longrightarrow \bicat E_{1}(e,e')_{b'};$$
and
\item parallel transport is compatible with the monoidal structure of $\cat B$, i.e., for all $\alpha\in \cat B^{op}(b,b')$, $\bar\alpha\in \cat B^{op}(\bar b,\bar b')$ and $\varphi\in \bicat E_{1}(e,e')_{b}$,
$\overline\varphi\in \bicat E_{1}(e',e'')_{\bar b}$,
$$\nabla _{e,e''}(\bar\alpha \otimes \alpha)(\overline\varphi\cdot\varphi)=\nabla _{e',e''}(\bar\alpha)(\overline\varphi)\cdot\nabla _{e,e'}(\alpha)(\varphi).$$
\end{enumerate}
\end{defn}

\begin{rmk}  In \cite{hess:bicatbndl}, both op-connections and their duals,  \emph{connections,} are studied in detail.
\end{rmk}

\begin{rmk} The monoidal compatibility condition says that one can parallel transport and then compose or compose and then parallel transport, with both sequences of operations leading to the same result.
\end{rmk}

\begin{defn} The \emph{fiber} of a bicategory bundle with op-connection $\Pi:\bicat E\to \bicat B$ is the bicategory $\bicat F$ with $\bicat F_{0}=\bicat E_{0}$, $\bicat F_{1}(e,e')=\bicat E_{1}(e,e')_{Id_{*}}$ for all $e,e'\in \bicat F_{0}$ and
$$\bicat F_{2}(\varphi, \psi)=\coprod _{\alpha \in \bicat B_{2}(Id_{*}, Id_{*})}\bicat E_{2}(\varphi, \psi)_{\alpha}$$
for all $\varphi, \psi\in \bicat F_{1}(e,e')$.  Here we are using the isomorphism 2-cell $\kappa: Id_{*}\overset {\cong}{\longrightarrow} Id_{*}\cdot Id_{*}$  in $\bicat B$ to define composition  in $\bicat F$ via
$$\bicat E_{1}(e,e')_{Id_{*}}\times \bicat E_{1}(e',e'')_{Id_{*}}\longrightarrow \bicat E_{1}(e,e'')_{Id_{*}\cdot Id_{*}}\overset {\nabla _{e,e''}(\kappa)}{\underset{\cong}{\longrightarrow}}\bicat E_{1}(e,e'')_{Id_{*}}.$$
\end{defn}

We are interested in using bicategory bundles with connection as a tool for comparing complex structures that can be encoded as such bundles.  To carry out comparisons, we need a precise description of bundle morphisms.

\begin{defn}\label{defi:bndlmorph}  Let $\Pi:\bicat E\to \bicat B$ and $\Pi':\bicat E'\to \bicat B'$ be bicategory bundles with op-connections $\boldsymbol \nabla$ and $\boldsymbol \nabla'$.  A morphism from $(\Pi, \boldsymbol\nabla)$ to $(\Pi',\boldsymbol{\nabla'})$ consists of a pair of bicategory homomorphisms $\Gamma: \bicat E\to \bicat E'$ and $\Lambda :\bicat B\to \bicat B'$ such that
$$\xymatrix
{\bicat E\ar[d]_{\Pi}\ar[r]^{\Gamma}&\bicat E'\ar[d]_{\Pi'}\\
\bicat B\ar[r]^{\Lambda}&\bicat B'}
$$
commutes.  Furthermore,  for all $e,e'\in\bicat E_{0}$, $\varphi\in \bicat E(e,e')_{b}$ and $\alpha\in \cat B^{op}(b,b')$,
$$\Gamma_{1}\big(\nabla_{e,e'}(\alpha)(\varphi)\big)=\nabla'_{\Gamma_{0}(e), \Gamma_{0}(e')}\big(L(\alpha)\big)\big(\Gamma_{1}(\varphi)\big),$$
where $L:(\cat B, \otimes)\to (\cat B',\otimes)$ is the strict monoidal functor associated to $\Lambda$.
\end{defn}

\begin{ex}  An elementary example of a bicategory bundle with op-connection is a product bundle.  Let $\bicat F$ be any bicategory, and let $\bicat B$ be a bicategory with exactly one object. The projection bicategory homomorphism $\bicat F\times \bicat B\to \bicat B$ is then a bicategory bundle with connection $\boldsymbol \nabla^\times$, where $\nabla^\times _{e,e'}(b)=\bicat F_{1}(e,e')\times \{b\}$ and $\nabla^\times _{e,e'}(\alpha)(\varphi,b)=(\varphi, b')$ for all $e,e'\in \bicat F_{0}$, $b,b\in \ob\cat B$, $\varphi\in \bicat F_{1}(e,e')_{b}$ and $\alpha \in \cat B^{op}(b,b')$.  Verifying the monoidal compatibility of parallel transport is trivial.  The fiber of the product bundle is, unsurprisingly, $\bicat F$.
\end{ex}

\subsection{Bicategory bundles and coalgebras}\label{ssec:bicatbndl-coalg}

In this section we apply  Cobar Duality to defining two seemingly different bicategory bundles with op-connection over level comonoids in the category $\cat {Ch}$ of connective chain complexes over a commutative ring $R$.  We show in the next section that these two bundles are indeed isomorphic, which leads in section \ref{sec:enrind} to the desired Enriched Induction theorem.

Recall that the category of level comonoids $\cat {Comon}_{\otimes}$ is monoidal with respect to the composition product, as shown at the end of section \ref{ssec:bimodules} .
Since $(\cat {Comon}_{\otimes},\circ,\op J)$ is a monoidal category, we can choose to view it as a bicategory $\bicat{CM}$ with exactly one object $*$, where $\bicat{CM}_{1}=\ob\cat {Comon}_{\otimes}$, with composition defined to be exactly the composition product, and $\bicat {CM}_{2}(\op X, \op Y)=\cat {Comon}_{\otimes}(\op X, \op Y)$.

The first  bicategory bundle over $\bicat {CM}$  expands the category $\cat {DCSH}$.
Consider the bicategory, $\bicat{DC}^\Om$, with the following cell structure.
\medskip
\begin{itemize}
\item $\omdc_{0}=\ob\coalg{A}$;
\medskip
\item $\bicat{DC}^\Om_{1}(C,C') =\coprod _{\op X \in \bicat{CM}_{1}} \cat M_{\text{mult.}}^\Sigma \big(\op T(\Om C)\circ \op X, \op T(\Om C')\big)$ for all $C,C\in \omdc_{0}$; and
\medskip
\item $\omdc_{2}(\varphi, \varphi')=\{\alpha\in \cat {Comon}_{\otimes}(\op Y,\op X)\mid \varphi'=\varphi(Id_{\op T(\Om C)}\circ \alpha)\}$, for all $\varphi\in \cat M_{\text{mult.}}^\Sigma \big(\op T(\Om C)\circ \op X, \op T(\Om C')\big)$ and $\varphi'\in \cat M_{\text{mult.}}^\Sigma \big(\op T(\Om C)\circ \op Y, \op T(\Om C')\big)$.
\end{itemize}
\medskip
Given two 1-cells in $\bicat {DC}^\Om$,  $\varphi\in \cat M_{\text{mult.}}^\Sigma \big(\op T(\Om C)\circ \op X, \op T(\Om C')\big)$ and $\psi\in \cat M_{\text{mult.}}^\Sigma \big(\op T(\Om C')\circ \op Y, \op T(\Om C'')\big)$, their composite $\psi\cdot\varphi$ in $\bicat {DC}_{1}$ is given by the following composition of multiplicative maps of symmetric sequences.
$$\op T(\Om C)\circ (\op X\circ \op Y)\cong\big(\op T(\Om C)\circ \op X\big)\circ \op Y\overset {\varphi\circ Id_{\op Y}}{\longrightarrow} \op T(\Om C')\circ \op Y\overset {\psi}{\longrightarrow} \op T(\Om C'')$$
Composition of $1$-cells is clearly associative up to natural isomorphism, since the composition product of level comonoids is associative up to natural isomorphism. Vertical composition of 2-cells is the same as composition of morphisms in $\cat {Comon}_{\otimes}$, while horizontal composition of $2$-cells is given by the monoidal structure of $\cat {Comon}_{\otimes}$.

Let $\Pi^\Om:\bicat {DC}^\Om\to \bicat{CM}$ denote the bundle specified as follows.
\medskip
\begin{itemize}
\item $\Pi^\Om_{0}(C)=*$ for all $C\in \bicat {DC}^\Om_{0}$;
\medskip
\item $\Pi^\Om _{1}(\varphi)=\op X$ for all $\varphi \in \cat M_{\text{mult.}}^\Sigma \big(\op T(\Om C)\circ \op X, \op T(\Om C')\big)$; and
\medskip
\item $\Pi^\Om _{2}(Id_{\op T(\Om C)}\circ \alpha)=\alpha$.
\end{itemize}
\medskip
The bundle $\Pi^\Om$ admits an op-connection $\boldsymbol{\nabla}^\Om$, which is defined by

$$\nabla_{C,C'}^\Om(-)= \cat M_{\text{mult.}}^\Sigma \big(\op T(\Om C)\circ -, \op T(\Om C')\big):\cat {Comon}_{\otimes}^{op }\longrightarrow \cat{Set}.$$
Note that the underlying category of the fiber of $\Pi^{\Om}$ is exactly $\cat {DCSH}$.

Consider now the bicategory, denoted $\bicat {DC}^\Phi$, with the following cell structure.
\begin{itemize}
\item $\bicat{DC}^\Phi_{0}=\ob\coalg{A}$;
\medskip
\item $\bicat{DC}^\Phi_{1}(C,C'):=\coprod _{\op X \in \bicat{CM}_{1}}\cat {Mod}^{tt}_{\op A}\big(\op T(C)\acirc \Phi(\op X), \op T(C')\big)$; and
\medskip
\item  $\bicat{DC}^\Phi_{2}(\varphi, \overline\varphi)=\{\alpha\in \bicat{CM}_{2}(\op Y, \op X)\mid \overline\varphi=\varphi(Id_{\op T(C)}\acirc \Phi(\alpha)\}$,
for all $\varphi\in \cat {Mod}^{tt}_{\op A}\big(\op T(C)\acirc \Phi(\op X), \op T(C')\big)$ and $\overline\varphi\in \cat {Mod}^{tt}_{\op A}\big(\op T(C)\acirc \Phi(\op Y), \op T(C')\big)$.
\end{itemize}

Given two 1-cells in $\bicat{DC}^\Phi$,  $\varphi\in \cat {Mod}^{tt}_{\op A}\big(\op T(C)\acirc \Phi(\op X), \op T(C')\big)$ and $\varphi'\in \cat {Mod}^{tt}_{\op A}\big(\op T(C')\acirc \Phi(\op X'), \op T(C'')\big)$, their composite $\varphi'\cdot\varphi$ in $\bicat{DC}^\Phi_{1}$ is given by the following formula.
$$\varphi'\cdot\varphi:=\lin\big[\ind(\varphi')\cdot \ind(\varphi)\big]$$
Here, $\ind(\varphi')\cdot \ind(\varphi)$ denotes the composite of $\ind(\varphi):\op T(\Om C)\circ \op X\to \op T(\Om C')$ and of $\ind(\varphi'):\op T(\Om C')\circ \op X'\to \op T(\Om C'')$, seen as $1$-cells of $\omdc$.
Composition of $1$-cells is associative up to isomorphism, since
\begin{align*}
(\varphi''\cdot \varphi')\cdot \varphi&=\lin\Big[\ind\lin[\ind(\varphi'')\cdot\ind(\varphi')]\cdot \ind(\varphi)\Big]\\
&=\lin\big[(\ind(\varphi'')\cdot\ind(\varphi'))\cdot \ind(\varphi)\Big]\\
&\cong\lin\big[\ind(\varphi'')\cdot(\ind(\varphi')\cdot \ind(\varphi))\Big]\\
&=\lin\big[\ind(\varphi'')\cdot\ind\lin[\ind(\varphi')\cdot \ind(\varphi)]\Big]\\
&=\varphi''\cdot (\varphi'\cdot \varphi).
\end{align*}
The second equality holds since $\ind\lin=Id$, while the isomorphism is due to the fact that composition of $1$-cells is associative up to isomorphism in $\omdc$.  Furthermore, if $\varepsilon: \Phi(\op J)\to \op A$ denotes the usual augmentation, then $Id_{\op T(C)}\acirc \varepsilon:\op T(C)\acirc \Phi(\op J)\to \op T(C)$ acts as an identity $1$-cell on $\op T(C)$, with respect to the composition we have defined, since $\ind(Id_{\op T(C)}\acirc \varepsilon)\cong Id_{\op T(\Om C)}$, by Scholium  \ref{schol:elpropind}(1).

Vertical composition of 2-cells is given by composition of morphisms in $\cat {Comon}_{\otimes}$, using the functoriality of $\Phi$, while horizontal composition is given by the  composition product of morphisms in $\cat {Comon}_{\otimes}$.  Horizontal composition of $2$-cells is compatible with composition of $1$-cells, since Scholia \ref{schol:elpropind} and \ref{schol:elproplin} imply that if the following two diagrams commute
$$\xymatrix{
\op T(C)\acirc \Phi(\op X)\ar[r]^(0.6){\varphi}&{\op T(C')}&&{\op T(C')\acirc \Phi(\op X')}\ar[r]^(0.6){\varphi'}&{\op T(C'')}\\
{\op T(C)\acirc \Phi(\op Y)}\ar[r]^(0.6){\psi}\ar[u]^{Id_{\op T(C)}\acirc \Phi(\alpha)}&{\op T(C')}\ar[u]^{=}&&{\op T(C')\acirc \Phi(\op Y')}\ar[r]^(0.6){\psi'}\ar[u]^{Id_{\op T(C')}\acirc \Phi (\alpha')}&{\op T(C'')}\ar[u]^{=}
}$$
then
$$\xymatrix{
{\op T(C)\acirc \Phi(\op X\circ \op X')}\ar[r]^(0.7){\varphi'\cdot \varphi}&{\op T(C'')}\\
{\op T(C)\acirc \Phi(\op Y\circ \op Y')}\ar[u]^{Id_{\op T(C)}\acirc \Phi(\alpha \circ \alpha')}\ar[r]^(0.7){\psi'\cdot \psi}&{\op T(C'')}\ar[u]^{=}
}$$
also commutes.
Since $(\beta\cdot \alpha)\circ (\beta '\cdot \alpha')=(\beta \circ \beta')\cdot (\alpha \circ\alpha')$ for all $\alpha \in \cat {Comon}_{\otimes}(\op X,\op Y)$, $\alpha' \in \cat {Comon}_{\otimes}(\op X',\op Y')$, $\beta \in \cat {Comon}_{\otimes}(\op Y,\op Z)$ and $\beta' \in \cat {Comon}_{\otimes}(\op Y',\op Z')$,  the composition maps
$$\phidc_{1}(C',C'')\times \phidc_{1}(C,C')\to \phidc _{1}(C,C'')$$
are indeed functors.

\begin{prop}\label{prop:invbifun} The natural transformations $\ind$ and $\lin$ induce mutually inverse bicategory homomorphisms
$$\Gamma_{\ind}:\phidc\longrightarrow \omdc\text{  and  }\Gamma_{\lin}:\omdc\longrightarrow\phidc,$$
where
\begin{itemize}
\item $(\Gamma_{\ind})_{0}(C)=C=(\Gamma_{\lin})_{0}(C)$ for all $C\in \phidc_{0}=\omdc_{0}$,
\smallskip
\item $(\Gamma_{\ind})_{1}(\varphi)=\ind \varphi$ for all $\varphi\in \phidc_{1}$ and $(\Gamma_{\lin})_{1}(\psi)=\lin \psi$ for all $\psi\in \omdc_{1}$, and
\smallskip
\item $(\Gamma_{\ind})_{2}(\alpha)=\alpha=(\Gamma_{\lin})_{2}(\alpha)$ for all $\alpha\in \phidc_{2}=\omdc_{2}$.
\end{itemize}
\end{prop}

\begin{proof}  We need only to verify that $(\Gamma_{\ind})_{1}$ and $(\Gamma_{\lin})_{1}$ preserve compositions and identities.  In the case of $(\Gamma_{\ind})_{1}$, we have that
\begin{multline*}
(\Gamma_{\ind})_{1}(\varphi'\cdot\varphi)=\ind\lin[\ind(\varphi')\cdot\ind (\varphi)]\\
=\ind(\varphi')\cdot \ind(\varphi)=(\Gamma_{\ind})_{1}(\varphi')\cdot (\Gamma_{\ind})_{1}(\varphi),
\end{multline*}
for all $\varphi, \varphi'\in \phidc_{1}$.  On the other hand,
\begin{multline*}
(\Gamma_{\lin})_{1}(\psi'\cdot\psi)=\lin[(\ind\lin \psi')(\ind\lin\psi)]\\
=(\lin \psi')\cdot(\lin \psi)=(\Gamma_{\lin})_{1}(\psi')\cdot(\Gamma_{\lin})_{1}(\psi).
\end{multline*}
\end{proof}

There is a bicategory bundle projection $\Pi^\Phi: \bicat{DC}^\Phi\to \bicat{CM}$ analogous to $\Pi^\Om$, which is defined as follows.
\begin{itemize}
\item $\Pi^\Phi _{0}(C)=*$ for all $C\in \bicat{DC}^\Phi_{0}$;
\item $\Pi^\Phi _{1}(\varphi)=\op X$ for all $\varphi \in \cat {Mod}^{tt}_{\op A}\big(\op T(C)\acirc \Phi(\op X), \op T(C')\big)\subset \bicat{DC}^\Phi_{1}(C,C')$;
\item $\Pi^\Phi_{2}(\alpha)=\alpha$.
\end{itemize}
It is easy to verify that $\Pi^\Phi$ is a bicategory homomorphism.  Using the decomposition notation defined above, we have
$$\bicat{DC}^\Phi_{1}(C,C')_{\op X}=\cat {Mod}^{tt}_{\op A}\big(\op T(C)\acirc \Phi(\op X), \op T(C')\big).$$

The bundle $\Pi^\Phi$ admits an op-connection $\boldsymbol{\nabla^\Phi}$, defined similarly to $\boldsymbol{\nabla^\Om}$.

For any $C,C'\in \phidc_{0}$, consider the functor
$$\nabla^\Phi_{C,C'}=\cat {Mod}^{tt}_{\op A}\big(\op T(C)\acirc \Phi(-), \op T(C')\big):\cat {Comon}_{\otimes}^{op }\longrightarrow \cat{Set}.$$
To establish that $\boldsymbol{\nabla^\Phi}$ satisfies monoidal compatibility, we explore its relationship with $\boldsymbol{\nabla^\Om}$, as mediated by $\ind$ and $\lin$.

\begin{lem}\label{lem:indconn} For all $\alpha\in \cat {Comon}_{\otimes}^{op}(\op X, \op Y)$ and $\varphi:\op T(C)\acirc \Phi(\op X)\to \op T(C')$,
$$\ind\big(\nabla^\Phi_{C,C'}(\alpha)(\varphi)\big)=\nabla^\Om_{C,C'}(\alpha)(\ind\varphi).$$
\end{lem}

This is simply Scholium \ref{schol:elpropind}(2) in a different guise, while the next lemma rephrases Scholium \ref{schol:elproplin}(2).

\begin{lem}\label{lem:linconn} For all $\alpha\in \cat {Comon}_{\otimes}^{op}(\op X, \op Y)$ and $\psi:\op T(\Om C)\circ \op X\to \op T(\Om C')$,
$$\lin\big(\nabla^\Om_{C,C'}(\alpha)(\psi)\big)=\nabla^\Phi_{C,C'}(\alpha)(\lin\psi).$$
\end{lem}

\begin{prop}  The op-connection $\boldsymbol{\nabla^\Phi}$ on $\Pi^\Phi:\phidc\to \bicat{CM}$ satisfies monoidal compatibility, i.e.,
$$\nabla^\Phi_{C,C''}(\alpha\circ \beta)(\psi\cdot\varphi)=\nabla^\Phi_{C',C''}(\beta)(\psi)\cdot \nabla^\Phi_{C,C'}(\alpha)(\varphi)$$
for all $\alpha\in \cat {Comon}_{\otimes}^{op}(\op X, \op X')$, $\beta \in \cat {Comon}_{\otimes}^{op}(\op Y, \op Y')$, $\varphi\in \phidc(C,C')_{\op X}$ and $\psi\in \phidc(C',C'')_{\op Y}$.
\end{prop}

\begin{proof} Observe that
\begin{align*}
\nabla^\Phi_{C',C''}(\beta)(\psi)\cdot \nabla^\Phi_{C,C'}(\alpha)(\varphi)&= \lin\Big[\ind\big(\nabla^\Phi_{C',C''}(\beta)(\psi)\big)\cdot \ind\big(\nabla^\Phi_{C,C'}(\alpha)(\varphi)\big)\big]\\
&=\lin\big[\nabla^\Om_{C',C''}(\beta)(\ind\psi)\cdot \nabla^\Om_{C,C'}(\alpha)(\ind\varphi)\big]\\
&\qquad\qquad\qquad\qquad\text{by Lemma \ref{lem:indconn}}\\
&=\lin\big(\nabla^\Om_{C',C''}(\beta)(\ind\psi)\big)\cdot \lin \big(\nabla^\Om_{C,C'}(\alpha)(\ind\varphi)\big)\\
&\qquad\qquad\qquad\qquad\text{by Proposition \ref{prop:invbifun}}\\
&=\lin\big(\nabla^\Om_{C,C''}(\alpha\circ\beta)(\ind\psi\cdot \ind\varphi)\big)\\
&\qquad\qquad\qquad\qquad\text{by monoidal compatibility of $\boldsymbol {\nabla^\Om}$}\\
&=\nabla^\Phi_{C,C''}(\alpha\circ\beta)(\lin\ind\psi\cdot \lin\ind\varphi)\\
&\qquad\qquad\qquad\qquad\text{by Lemma \ref{lem:linconn} and Proposition \ref{prop:invbifun}}\\
&=\nabla^\Phi_{C,C''}(\alpha\circ\beta)(\psi\cdot \varphi).
\end{align*}
\end{proof}

\begin{rmk}  The category underlying the fiber $\bicat F^\Phi$ of $\Pi^\Phi$ is exactly the ``fattened'' category $\Coalg{A}{\op F}$, where $\op F=\Phi (\op J)$ (cf., section 5).  It is a straightforward calculation that for all transposed tensor maps $\varphi:\op T(C)\acirc \op F\to \op T(C')$ and $\psi:\op T(C')\acirc \op F\to \op T(C'')$, the composite $\lin\big[\ind \psi\cdot \ind\varphi\big]$  in $\bicat F^\Phi _{1}$
is the same as the composite in $\Coalg{A}{\op F}$, given by
$$\op T(C)\acirc \op F\xrightarrow{Id_{\op T(C)}\acirc \psi _{\op F}}\op T(C)\acirc \op F\acirc \op F\xrightarrow{\varphi\acirc Id_{\op F}}\op T(C')\acirc \op F\xrightarrow \psi \op T(C''),$$
i.e., composition is defined the same way in both categories.
\end{rmk}

Arranging the results above and their consequences somewhat differently, one obtains the following Enriched Cobar Duality Theorem.

\begin{thm}\label{thm:ecdt} The bicategory homomorphisms $\Gamma_{\ind}$ and $\Gamma_{\lin}$ give rise to mutually inverse morphisms of bicategory bundles with connection
$$\xymatrix{
{\phidc}\ar[d]_{\Pi^\Phi}\ar[r]^{\Gamma_{\ind}} &{\omdc}\ar[d]_{\Pi^\Om}\ar[r]^{\Gamma_{\lin}}\ar[d]_{\Pi^\Om}&{\phidc}\ar[d]_{\Pi^\Phi}\\
{\bicat{CM}}\ar[r]^{=}&{\bicat{CM}}\ar[r]^{=}&{\bicat{CM}.}
}$$
\end{thm}

\begin{proof}  The squares of the diagram above clearly commute.  Furthemore $\ind $ and $\lin$ give rise to families of natural transformations
$$\{\nabla^\Phi_{C,C'}\to\nabla^\Om_{C,C'}\mid C,C'\in \ob \coalg{A}_{+}\}$$
 and
 $$\{\nabla^\Om_{C,C'}\to\nabla^\Phi_{C,C'}\mid C,C'\in \ob \coalg{A}_{+}\},$$
 respectively.  Finally, Lemmas \ref{lem:indconn} and \ref{lem:linconn} imply that the last condition in Definition \ref{defi:bndlmorph}  holds for $\Gamma_{\ind}$ and for $\Gamma_{\lin}$, respectively.
\end{proof}

\section{Enriched induction}\label{sec:enrind}

Let  $\theta: \op T(C)\acirc \Phi(\op X)\to \op T(C')$ be a transposed tensor morphism, where $C$ and $C'$ are chain coalgebras,  and $\op X$ is a level comonoid.  If $\op X$ is a right $\op P$-module and $\Om C'$ is a $\op P$-coalgebra it is natural to ask when $\ind \theta:\op T(\Om C)\circ \op X\to\op T(\Om C')$ is a morphism of right $\op P$-modules.  Similarly, when $\op X$ is a left $\op P$-module and $\Om C$ is a $\op P$-coalgebra, we can ask when $\ind\theta: \op T(\Om C)\circ \op X\to \op T(\Om C')$ induces a morphism $\widehat{\ind\theta}: \op T(\Om C)\underset {\op P}\circ \op X\to \op T(\Om C')$.

In this section, we apply the Enriched Cobar Duality Theorem to answering these questions.

\subsection{Diffracted module maps}

Let $\op P$ be a Hopf operad in the category of chain complexes, i.e., $\op P$ is a level comonoid in the category of operads.  Our most important example of a Hopf operad is the associative operad $\op A$.   The categories of right and left $\op P$-modules, as well as the category of $\op P$-bimodules, are monoidal with respect to the level monoidal product.  For example, given a right $\op P$-action $\rho:\op X\circ\op P\to \op X$, the level tensor product $\op X\otimes \op X$ admits a right $\op P$-action
given by the composite
$$(\op X\otimes\op X)\circ \op P\overset {Id\otimes \Delta_{\op P}}{\longrightarrow}(\op X\otimes\op X)\circ(\op P\otimes \op P)\overset{\iota}{\longrightarrow}(\op X\circ\op P)\otimes (\op X\circ\op P)\overset{\rho\otimes \rho}{\longrightarrow}\op X\otimes \op X.$$

Let $\cat F_{\op P}$ denote the following category.  Its objects are chain coalgebras $C$, together with a choice of $\op P$-coalgebra structure on  $\Om C$  that is compatible with its multiplicative structure.  In other words, there is a fixed  multiplicative right $\op P$-module action $\op T(\Om C)\circ \op P\to \op T(\Om C)$.   Morphisms in $\cat F_{\op P}$ are chain coalgebra morphisms $f:C\to C'$ such that $\Om f$ is a morphism of $\op P$-coalgebras.  When $\op P=\op A$, we write $\cat F=\cat F_{\op A}$, and refer to the objects as \emph{Alexander-Whitney coalgebras} \cite {hpst:04}.

In the definitions below, we specify the additional conditions that must be satisfied by a transposed tensor map in order for the multiplicative map it induces to be more highly structured, in the sense of the questions asked above.

\begin{defn}\label{defi:diffrtmod} Let $C$ be a chain coalgebra, and let $C'$ be an object in $\cat F_{\op P}$ with multiplicative action map $\psi':\op T(\Om C')\circ \op P\to \op T(\Om C')$.  Let $(\op X, \Delta, \rho)$ be a level comonoid in the category of  right $\op P$-modules, i.e.,
$\rho:\op X\circ \op P\to \op X$  is the right action of $\op P$ on $\op X$ and $\Delta:\op X\to \op X\otimes \op X$ is a morphism of right $\op P$-modules.  A transposed tensor morphism $\theta:\op T(C)\acirc \Phi(\op X)\to \op T(C')$ is a \emph {diffracted right $\op P$-module map} if the following diagram commutes.
$$\xymatrix{
{\op T(C)\acirc \Phi(\op X\circ \op P)}\ar[d]_{Id_{\op T(C)}\acirc \Phi (\rho)}\ar[rd]^{\quad\lin\big[\psi'(\ind \theta \circ Id_{\op P})\big]}\\
{\op T(C)\acirc \Phi(\op X)}\ar[r]^{\theta}&{\op T(C')}
}$$
\end{defn}

To help the reader understand this definition, we note that the diagonal arrow in the diagram above is obtained by applying linearization to the composite
$$\op T(\Om C)\circ\op X\circ \op P\xrightarrow{\ind \theta \circ Id_{\op P}}\op T(\Om C')\circ\op P\xrightarrow {\psi'}\op T(\Om C'),$$
which is exactly the composite of $\psi'$ and $\ind \theta$ as $1$-cells of $\omdc$.

\begin{defn}\label{defi:diffbalmod}
Let $C$ be an object of $\cat F_{\op P}$ with multiplicative action map  $\psi:\op T(\Om C)\circ \op P\to \op T(\Om C)$, and let $C'$ be a chain coalgebra .  Let $(\op X, \Delta, \lambda)$ be a level comonoid in the category of  left $\op P$-modules, i.e.,
$\lambda:\op P\circ \op X\to \op X$  is the left action of $\op P$ on $\op X$ and $\Delta:\op X\to \op X\otimes \op X$ is a morphism of left $\op P$-modules.  A transposed tensor morphism $\theta:\op T(C)\acirc \Phi(\op X)\to \op T(C')$ is a \emph {diffracted balanced $\op P$-module map} if the following diagram commutes.
$$\xymatrix{
{\op T(C)\acirc \Phi(\op P\circ \op X)}\ar[d]_{Id_{\op T(C)}\acirc \Phi (\lambda)}\ar[dr]^{\quad\lin\big[\ind \theta(\psi\circ Id_{\op X})\big]}\\
{\op T(C)\acirc \Phi(\op X)}\ar[r]^{\theta}&{\op T(C')}
}$$
\end{defn}

Again with the aim of aiding the reader, we note that the diagonal arrow in the diagram above is obtained by applying linearization to the composite
$$\op T(\Om C)\circ \op P\circ \op X\xrightarrow{\psi\circ Id_{\op X}}\op T(\Om C)\circ \op X\xrightarrow{\ind \theta}\op T(\Om C'),$$
which is exactly the composite of $\psi$ and $\ind\theta$ as $1$-cells of $\omdc$.

The next theorem, which answers the questions asked in the introduction to this section, also justifies the terminology introduced in the two definitions above.

\begin{thm}\label{thm:diff-equiv}
\begin{enumerate}
\item Let $C$ be a chain coalgebra, and let $C'$ be an object in $\cat F_{\op P}$.  Let $(\op X, \Delta, \rho)$ be a level comonoid in the category of  right $\op P$-modules.

A transposed tensor map $\theta:\op T(C)\acirc \Phi(\op X)\to \op T(C')$ is a diffracted right $\op P$-module morphism if and only if $\ind \theta:\op T(\Om C)\circ \op X\to\op T(\Om C')$ is a right $\op P$-module map.

\item Let $C$ be an object in $\cat F_{\op P}$, and let $C'$ be a chain coalgebra.  Let $(\op X, \Delta, \lambda)$ be a level comonoid in the category of  left $\op P$-modules.

A transposed tensor map $\theta:\op T(C)\acirc \Phi(\op X)\to \op T(C')$ is a diffracted balanced $\op P$-module morphism if and only if  $\ind \theta:\op T(\Om C)\circ \op X\to\op T(\Om C')$ induces a morphism of symmetric sequences $\widehat{\ind \theta}:\op T(\Om C)\underset {\op P}\circ \op X\to\op T(\Om C')$.
\end{enumerate}
\end{thm}

\begin{proof}  We use here the notation of Definitions \ref{defi:diffrtmod} and \ref{defi:diffbalmod}.

To prove (1), we need to show that
$$\xymatrix{
{\op T(\Om C)\circ \op X\circ \op P\;}\ar[rr]^{\ind \theta\circ Id_{\op P}}\ar[d]_{Id_{\op T(\Om C)}\circ \rho}&&{\;\op T(\Om C')\circ \op P}\ar[d]_{\psi'}\\
{\op T(\Om C)\circ \op X}\ar[rr]^{\ind\theta}&&{\op T(\Om C')}
}$$
commutes if and only if the diagram in  Definition \ref{defi:diffrtmod} commutes.  This is an immediate consequence, however, of Theorem \ref{thm:ecdt}, which implies that $\ind$ and $\lin$ are mutual inverses preserving compositions and that
$$\ind\big(\theta(Id_{\op T(C)}\acirc \Phi (\rho))\big)=\ind \theta(Id_{\op T(\Om C)}\circ \rho)$$
and
$$\lin\big(\ind\theta(Id_{\op T(\Om C)}\circ \rho)\big)=\theta(Id_{\op T(C)}\acirc \Phi (\rho)).$$

Similarly, to prove (2), we must show that
$$\xymatrix{
{\op T(\Om C)\circ \op P\circ \op X}\ar[rr]^{\psi\circ Id_{\op X}}\ar[d]_{Id_{\op T(\Om C)}\circ \lambda}&&{\op T(\Om C)\circ \op X}\ar[d]_{\ind\theta}\\
{\op T(\Om C)\circ \op X}\ar[rr]^{\ind\theta}&&{\op T(\Om C')}
}$$
commutes if and only if the diagram in Definition  \ref{defi:diffbalmod} commutes.  The proof again follows directly from Theorem \ref{thm:ecdt}.
\end{proof}

\begin{cor}\label{cor:enrind}  Let $(C,\psi)$ and $(C',\psi')$ be objects in $\cat F_{\op P}$, and let $(\op X, \lambda, \rho, \Delta)$ be a level comonoid in the category of $\op P$-bimodules.  A transposed tensor map $\theta: \op T(C)\acirc \Phi(\op X)\to \op T(C')$ that is both a diffracted right $\op P$-module morphism and a diffracted balanced $\op P$-module morphism induces a morphism of right $\op P$-modules $$\widehat{\ind \theta}:\op T(\Om C)\underset{\op P}\circ \op X\to\op T(\Om C').$$
\end{cor}

We record here for future use a helpful result concerning sums of enriched transpose tensor maps.

\begin{cor}\label{cor:sum-enrtransptens} Let $C,C', C''$ be objects in $\cat F_{\op P}$, where $\op P$ is a Hopf operad, and let $\op M$ be a level comonoid in the category of $\op P$-bimodules.  Let $\theta':\op T(\Om C')\underset{\op P}\circ \op M\to \op T(\Om C)$ and $\theta'':\op T(\Om C')\underset{\op P}\circ \op M\to \op T(\Om C)$ be multiplicative morphisms of right $\op P$-modules.   Then there is a unique multiplicative morphism of right $\op P$-modules $\theta:\op T(\Om (C'\oplus C''))\underset{\op P}\circ \op M\to \op T(\Om C)$ restricting to $\theta'$ and $\theta''$.
\end{cor}

\begin{proof}  The hypotheses of the corollary imply that $\lin(\theta'):\op T(C')\acirc \Phi(\op X)\to \op T(C)$ and $\lin (\theta''):\op T(C'')\acirc \Phi (\op X)\to \op T(C)$ are transposed tensor morphisms that are diffracted balanced  and diffracted right $\op P$-module morphisms.  By Proposition \ref{prop:sum-transptens}, there is a unique transposed tensor morphism $\tau:\op T(C'\oplus C'')\acirc \Phi(\op X)\to \op T(C)$ restricting to $\lin (\theta')$  and $\lin (\theta'')$.  The formula for the sum $\tau$ given in the proof of Proposition \ref{prop:sum-transptens} enables us to verify by inspection that $\tau$, too, is a diffracted balanced morphism and a diffracted right $\op P$-morphism.  We may therefore set $\theta=\ind (\tau)$.
\end{proof}

\subsection{Existence of diffracted module maps}

We show here how to apply acyclic models methods in order to prove the existence of diffracted module maps and therefore of highly-structured induced morphisms.  As in the previous section, $\op P$ denotes a Hopf operad in the category of chain complexes.

Throughout this section, as in section \ref{ssec:acyclic}, we work in the comma category (or under category) $\op J\downarrow \cat M^\Sigma$.  The objects of $\op J\downarrow \cat M^\Sigma$ are morphisms of symmetric sequences $\op J\to \op X$ and are called \emph{symmetric sequences under $\op J$}.  The morphisms are commuting triangles.  The usual composition monoidal structure on $\cat M^\Sigma$ induces a monoidal product of symmetric sequences under $\op J$, given by
$$(\op X,j_{\op X})\circ (\op Y,j_{\op Y}):=\big(\op X\circ\op Y, (j_{\op X}\circ j_{\op Y})\kappa_{\circ}\big),$$
where $\kappa_{\circ}:\op J\xrightarrow{\cong}\op J\circ\op J$ is the canonical isomorphism.
The unit is the identity map $\op J\xrightarrow{=}\op J$.  An \emph{operad under $\op J$} is a monoid with respect to this monoidal structure on $\op J\downarrow \cat M^\Sigma$, in which case $j_{\op P}:\op J\to \op P$ is necessarily the unit map.  In fact, since we have defined operads to be unital monoids with respect to the composition product of symmetric sequences, any operad is an operad under $\op J$. Given an operad $(\op P,j_{\op P})$ under $\op J$, we can define in the obvious way categories $\op J\downarrow({}_{\op P}\cat {Mod})$ of left $\op P$-modules under $\op J$, $\op J\downarrow(\cat {Mod}_{\op P})$ of right $\op P$-modules under $\op J$ and $\op J\downarrow({}_{\op P}\cat {Mod}_{\op P})$ of $\op P$-bimodules under $\op J$.

Note that if $(\op X,j_{\op X})$ and $(\op Y, j_{\op Y})$ are symmetric sequences under $\op J$, then there are natural ``inclusion'' morphisms in $\op J\downarrow \cat M^\Sigma(\op X,\op X\circ \op Y)$ and $\op J\downarrow \cat M^\Sigma(\op X, \op Y\circ \op X)$  given by the composites
$$\op X\xrightarrow{\cong}\op X\circ\op J\xrightarrow{Id_{\op X}\circ j_{\op Y}}\op X\circ \op Y$$
and
$$\op X\xrightarrow{\cong}\op J\circ\op X\xrightarrow{j_{\op Y}\circ Id_{\op X}}\op Y\circ \op X.$$
We work in this section with objects under $\op J$ precisely because we need such inclusions to construct diffracted module maps by acyclic models methods.

The level monoidal structure on $\cat M^\Sigma$ induces a symmetric monoidal structure on $\op J\downarrow \cat M^\Sigma$, defined by
$$(\op X, j_{\op X})\otimes (\op Y,j_{\op Y}):=\big(\op X\otimes \op Y, (j_{\op X}\otimes j_{\op Y})\kappa_{\otimes}\big),$$
where $\kappa_{\otimes}:\op J\xrightarrow{\cong}\op J\otimes\op J$ is the canonical isomorphism.  It makes sense therefore to speak of level comonoids under $\op J$.

The discussion above applies to any closed, symmetric monoidal category $\cat M$.  We work henceforth in $\cat {Ch}$, the category of chain complexes over a commutative ring $R$. There are obvious morphisms of  symmetric sequences from $\op J$ into each of $\op A$, $\op S$ and $\op A^\perp$.  If $(\op M, \Delta, j_{\op M})$ is a level comonoid under $\op J$, then there is an induced morphism of left $\op A$-modules $\hat\jmath _{\op M}:\op A\to \Phi (\op M)$ given by the composite
$$\op A\xrightarrow{\cong}\op A\circ \op J\circ \op J\circ\op J\to\op A\circ (\op S\otimes \op M)\circ \op A^\perp\circ \op A,$$
where the second morphism is built from the obvious morphisms $\op J\to \op S\otimes \op M$, $\op J\to \op A^\perp$ and $\op J\to \op A$.  The original diffracting functor thus induces a diffracting functor relative to $\op J$
$$\Phi: \op J\downarrow \cat {Comon}_{\otimes}\longrightarrow \op A\downarrow({}_{\op A}\cat {Mod}_{\op A}).$$

 Given a (Hopf) operad $(\op P, j_{\op P})$ under $\op J$, the image of $j_{\op P}$, denoted $\op P_{\star}$, is a sub (Hopf) operad of $\op P$, concentrated in level $1$ and degree $0$.  In particular, since $\op J(1)=R$, we can write $\op P_{\star}(1)=R\cdot e_{\op P}$, where $e_{\op P}=j(1)$.  For example, $e_{\op A}=\delta ^{(1)}$.  Furthermore,  an object $\op V$ in $\op J\downarrow \cat {Ch}^\Sigma$ is not actually a subsymmetric sequence of the free right $\op P$-module $\op V\circ \op P$.  To be absolutely precise, if slightly pedantic, it is the isomorphic symmetric sequence $\op V\circ \op P_{\star}$ that is  a subsymmetric sequence of the free right $\op P$-module $\op V\circ \op P$.  A similar observation holds for free left $\op P$-modules and free $\op P$-bimodules.

Note that a filtration
$$\cdots \subset F^n\op V\subset F^{n+1}\op V\subset \cdots$$
of a symmetric sequence $\op V$ of graded $R$-modules induces a filtration of the free right $\op P$-bimodule $\op V\circ\op P$, where $F^n(\op V\circ \op P):= F^n\op V\circ \op P$.  There are similar induced filtrations on the free left $\op P$-module and on the free $\op P$-bimodule generated by $\op V$.

\begin{defn}\label{defn:semifree} Let $(\op P, j_{\op P})$ be an operad under $\op J$.  A right $\op P$-module $\op M$ under $\op J$ in $\cat {Ch}$ is \emph{semifree} if there is a symmetric sequence of free graded $R$-modules $\op V$ under $\op J$, endowed with an increasing filtration $$\op J=F^0\op V\subset F^1 \op V\subset \cdots \subset   F^n\op V\subset F^{n+1}\op V\subset \cdots$$ such that
\begin{enumerate}
\item for each $n$, $F^{n+1}\op V$ splits as $F^n\op V\oplus E^{n+1}\op V$, where $E^{n+1}\op V$ is a free graded $R$-module;
\item the  symmetric sequence of graded $R$-modules underlying $\op M$ is $ \op V\circ\op P$; and
\item the differential $\del _{\op M}$ of $\op M$ satisfies $\del _{\op M}(E^{n+1}\op V\circ \op P_{\star})\subset F^{n}\op V\circ \op P$ for all $n$.
\end{enumerate}
\end{defn}

\begin{defn}  A level comonoid of chain complexes $(\op M,\Delta)$ in $\op J\downarrow(\cat {Mod}_{\op P})$  is \emph{semifree} if the underlying right $\op P$-module is semifree on $\op V$ and
$$\Delta (E^{n+1}\op V\circ \op P_{\star})\subset \left(F^n\op M\oplus (E^{n+1}\op V\circ \op P_{\star})\right)^{\otimes 2}.$$
\end{defn}

When $(\op M, \Delta)$ is semifree on $\op V$, the coproduct $\Delta$ is determined by its restriction to $\op V\circ \op P_{\star}$.

\begin{rmk}  There are similar definitions of semifree left $\op P$-modules and of semifree $\op P$-bimodules under $\op J$, as well as of semifree level comonoids in $\op J\downarrow({}_{\op P}\cat {Mod})$ and in $\op J\downarrow({}_{\op P}\cat {Mod}_{\op P})$.  We leave the exact formulation of the definitions to the reader.
\end{rmk}

Recall that a functor $X: \cat D\to \cat {Ch}$, where $\cat D$ is a small category, is globally connective if there is an integer $N$ such that $X(d)_{n}=0$ for all $n<N$ and for all $d\in \ob \cat D$. The integer $N$ is called a global lower bound on the functor $X$.  In particular, a symmetric sequence $\op X$ is globally connective if there is an integer $N$ such that $\op X(m)_{n}=0$ for all $n<N$ and for all $m\geq 0$.  Observe that if $\op X$ is globally connective, then so is $\Phi (\op X)$.  In particular, $\op F=\Phi(\op J)$ is globally connective.

\begin{exs}
\begin{enumerate}
\item Every $\op J$-module is trivially semifree, as is any level comonoid in the category of $\op J$-modules.

\item Given any globally connective, level comonoid $(\op X,\Delta)$ of symmetric sequences  under $\op J$, its diffraction $\Phi(\op X)$ is a semifree $\op A$-bimodule under $\op J$.  Recall from section 4 that the underlying $\op A$-bimodule of $\Phi(\op X)$ is $\op A\circ \op V_{\op X}\circ \op A$, where $\op V_{\op X}=(\op S\otimes \op X)\circ \op A^\bot$.  Let $\op X_{\leq p}$ denote the symmetric subsequence of $\op X$ consisting of elements of degree at most $p$, and let $N$ be the global lower bound on the chain complexes $\op X(m)$.  For any $q\geq 1$, let $\op A_q^\perp$ be the symmetric sequence given by
$$\op A_{q}^\perp(l)= \left\{
        \begin{array}{cl}
           \op A^\perp(l)    & \text{if }l\leq q, \\
            0       & \text{if }l>q.
        \end{array}
        \right.
$$
For all $n\geq 1$, set
$$F^n\op V_{\op X}=\op J+\bigoplus _{0\leq j\leq n}(\op S\otimes \op X_{< N+j})\circ \op A^\perp_{n-j}$$
and set $F^0\op V_{\op X}=\op J$.  It is obvious that $F^{n+1} \op V_{\op X}$ splits, as required by condition (1) in Definition \ref{defn:semifree}. Furthermore, the formula for the differential on $ \Phi (\op X)$ in Remark \ref{rmk:formula-diffl} implies that the bimodule version of condition (3)  is satisfied.

In the case $\op X=\op J$, inspection of the formula for the level coproduct on $\op F$ in section 5 leads easily to the conclusion that $\Phi(\op J)=\op F$ is even a semifree level comonoid under $\op J$.
\end{enumerate}
\end{exs}

The following small, technical result plays a crucial role in the proof of the existence theorem below and is even of interest in and of itself.  It specifies one case in which applying $\op X\circ-$ commutes with coproducts, even though in general $\op X\circ -$ does not preserve colimits.

\begin{lem}\label{lem:coprod}  Let $A$ be an associative chain algebra, and let $\op B$ be a level monoid.  Let $\op M$ and $\op N$ be level comonoids.  Given multiplicative maps of symmetric sequences $\varphi:\op T(A)\circ \op M\to \op B$ and $\theta:\op T(A)\circ \op N\to \op B$ that are determined by the families $\{\varphi_{m}:A\otimes \op M(m)\to \op B(m)\}_{m}$ and  $\{\psi_{m}:A\otimes \op N(m)\to \op B(m)\}_{m}$ of multiplicative equivariant maps, there is a unique multiplicative morphism of symmetric sequences
$\varphi*\psi:\op T(A)\circ (\op M\oplus \op N)\to \op B$ such that
$$\xymatrix{
{\op T(A)\circ\op M}\ar[drr]^{\varphi}\ar[d]_{Id_{\op T(A)}\circ \iota _{\op M}}\\
{\op T(A)\circ(\op M\oplus \op N)}\ar[rr]^{\varphi*\psi}&&{\;\op B}\\
{\op T(A)\circ\op N}\ar[urr]_{\psi}\ar[u]^{Id_{\op T(A)}\circ \iota _{\op N}}
}$$
commutes, where $\iota _{\op M}:\op M\hookrightarrow \op M\oplus \op N$ and $\iota _{\op N}:\op N\hookrightarrow \op M\oplus \op N$ are the inclusions.
\end{lem}

\begin{proof}  Let $\kappa_{m}:A\otimes\big(\op M(m)\oplus \op N(m)\big)\to \big(A\otimes \op M(m)\big)\oplus\big(A\otimes \op N(m)\big)$ denote the canonical isomorphism.  The family
$$\big\{(\varphi_{m}+\psi_{m})\kappa_{m}:A\otimes \big(\op M(m)\oplus \op N(m)\big)\to \op B(m)\}_{m}$$
is multiplicative and equivariant and therefore induces the desired multiplicative morphism of symmetric sequences.  This is obviously the unique possibility.
\end{proof}

We refer the reader to section \ref{ssec:acyclic} for the definition of functors that are \emph{free} or \emph{acyclic} with respect to a set of models.

\begin{thm}\label{thm:diffmod}  Let $X: \cat D\to\coalg {A}$ and $Y:\cat D\to \cat F_{\op P}$ be functors, where $\cat D$ is a category admitting a set of models $\frak M$ with respect to which  $X$ is free and globally connective and $Y$ is acyclic.  Let $(\op M, \Delta, j_{\op M})$ be a globally connective, semifree, level comonoid in the category of right $\op P$-modules under $\op J$.  Let $\tau: UX\to UY$
be a natural transformation, where $U$ is the forgetful functor down to $\cat {Ch}$.  Then there is a
multiplicative natural transformation  of right $\op P$-modules
$$\theta :\op T(\Om X)\circ \op M\longrightarrow \op T(\Om Y)$$
extending $s^{-1}\tau$, i.e., the following composite is equal to $\tau$.
$$X\xrightarrow{s^{-1}}s^{-1}X\hookrightarrow \op T(\Om X)\circ \op J\to \op T(\Om X)\circ \op M\xrightarrow{\theta}\op T(\Om Y)\xrightarrow {\text{proj.}}s^{-1}Y\xrightarrow {\cong}Y$$
\end{thm}

\begin{proof} We prove the existence of a natural, diffracted right $\op P$-module  transformation $\widehat \tau:\op T(X)\acirc \Phi(\op M)\to \op T(Y)$ extending $\tau$, i.e., such that the composite of the following sequence of maps is equal to $\tau$.
$$X\hookrightarrow \op T(X)\acirc\op A \overset{Id_{}\acirc \hat\jmath_{\op M}}{\longrightarrow}\op T(X)\acirc\Phi(\op M)\overset{\widehat\tau}{\longrightarrow} \op T(Y)$$
Theorem \ref{thm:diff-equiv} then implies that $\theta=\ind(\widehat\tau)$ is a right $\op P$-module map, as desired.

The proof proceeds by induction on the degrees of the generating elements $x_{\frak m}$, which is possible since $X$ is globally connective, on filtration degree in $\op M$ and on degree in the ``simplicial'' filtration of $\Phi (\op M)$ already used above, which we begin by recalling.

For any $p\geq 1$, let $\op A_p^\perp$ be the symmetric sequence given by
$$\op A_{p}^\perp(l)= \left\{
        \begin{array}{cl}
           \op A^\perp(l)    & \text{if }l\leq p, \\
            0       & \text{if }l>p.
        \end{array}
        \right.
$$
Given any level comonoid $\op N$, let
$$\widetilde F^p\Phi (\op N)=\op A\circ (\op S\otimes \op N)\circ \op A_{p}^\perp\circ \op A$$
for all $p\geq 1$, defining an increasing filtration of $\Phi (\op N)$, at least as an $\op A$-bimodule of graded $R$-modules. The formula for the full differential $\del$ in $\Phi(\op N)$ (cf., Remark \ref{rmk:formula-diffl}) implies that this is actually a differential filtration and that both the cosimplicial and the simplicial parts of the differential are filtration-lowering on elements of the form $e_{\op A}\otimes s_{k-1}x\otimes \alpha _{\vec l}\otimes \delta^{(\vec m)}$.  To see that the cosimplicial part is filtration-lowering, recall that $\op A^\perp(0)=0$, which implies that $\psi(\alpha_{n})=\sum _{0< k<n} \alpha_{k}\otimes \alpha _{n-k}$.

Before stating the inductive hypotheses, we fix further notation.  Let $I_{n}^p=\{(l_{1},...,l_{n})\in \mathbb N^n\mid 1\leq l_{i}\leq p \;\forall i\}$.  Let $\op V$ denote the filtered symmetric sequence of free graded $R$-modules underlying the semifree level comonoid $\op M$, so that $\op M\cong \op V\circ \op P$ as symmetric sequences of graded modules.  Recall that for each $n$, we assume that there is a splitting $F^{n}\op V=F^{n-1} \op V\oplus E^{n}\op V$ as free graded $R$-modules.  Furthermore, $F^{n-1}\op M\oplus E^{n}\op V\circ \op P_{\star}$ is a sub level comonoid of $\op M$, for all $n$.   For each $k$, we fix a basis $B_{n,k}$ of $E^{N}\op V(k)$.  Given $v\in B_{n,k}$, we let $\tilde v=v\otimes e_{\op P}^{\otimes k}\in E^{n}\op V\circ \op P_{\star}$.

For all $d\in \ob \cat D$, let $\del_{d}$ denote the differential on $\op T\big(X(d)\big)\acirc \Phi(\op M)$, which is formed from the differential on $\op T\big(X(d)\big)$ naturally induced by that on $X(d)$ and from the differential on $\Phi(\op M)$.  The differential on $\op T\big(Y(d)\big)$ is denoted $\del_{d} '$.

For all $m,n\in \mathbb N$ and $p\geq 1$ and for all $d\in \ob\cat D$, let
$$\op K_{m,n,p}(d)=\op T\big(X(d)_{<m}\big)\acirc \Phi(\op M)+\op T\big(X(d)_{m}\big)\acirc \widetilde F^p\Phi (F^{n-1}\op M\oplus E^{n}\op V\circ \op P_{\star}),$$
$$\op K_{m,n}(d)=\op T\big(X(d)_{<m}\big)\acirc \Phi(\op M)+\op T\big(X(d)_{m}\big)\acirc \Phi (F^n\op M),$$
and
$$\op K_{m}(d)=\op T\big(X(d)_{\leq m}\big)\acirc \Phi(\op M),$$
which are right $\op A$-submodules of $\op T\big(X(d)\big)\acirc\Phi(\op M)$.

Our inductive hypotheses can be formulated as follows.
\medskip
\begin{description}
\item[$\text{H}_{m,n,p}$] For all objects $d$ in $\cat D$ and for all $j\leq m, k\leq n, l\leq p$, there are natural transposed tensor maps of right $\op A$-modules
$$\widehat \tau_{j,k,l}(d):\op K_{j,k,l}(d)\longrightarrow \op T\big(Y(d)\big)$$
that restrict to diffracted right $\op P$-module maps on $\op T\big(X(d)_{<m}\big)\acirc \Phi(\op M)$ and that extend $\tau(d):X(d)_{\leq m}\to Y(d)$.   Furthermore, if $j\leq j'$, $k\leq k'$ and $l\leq l'$, then the restriction of $\widehat \tau_{j',k',l'}(d)$ to $\op K_{j,k,l}(d)$ is equal to $\widehat \tau_{j,k,l}(d)$ for all $d\in \ob\cat D$.
\medskip
\item[$\text{H}'_{m,n}$] For all objects $d$ in $\cat D$ and for all $j\leq m, k\leq n$, there are natural transposed tensor maps
$$\widehat \tau_{j,k}(d):\op K_{j,k}(d)\longrightarrow \op T\big(Y(d)\big)$$
that are diffracted right $\op P$-module maps and that extend $\tau(d):X(d)_{\leq m}\to Y(d)$.   Furthermore, if $j\leq j'$ and $k\leq k'$, then the restriction of $\widehat \tau_{j',k'}(d)$ to $\op K_{j,k}(d)$ is equal to $\widehat \tau _{j,k}(d)$ for all $d\in \ob\cat D$.
\medskip
\item[$\text{H}''_{m}$] For all objects $d$ in $\cat D$ and for all $j\leq m$, there are natural transposed tensor maps
$$\widehat \tau_{j}(d):\op K_{j}(d)\longrightarrow \op T\big(Y(d)\big)$$
that are diffracted right $\op P$-module maps and that extend $\tau(d):X(d)_{\leq m}\to Y(d)$.   Furthermore, if $j\leq j'$, then the restriction of $\widehat \tau_{j'}(d)$ to $\op K_{j}(d)$ is equal to $\widehat \tau _{j}(d)$ for all $d\in \ob\cat D$.
\end{description}

If $\text{H}'_{m,n}$ is satisfied for all $n$, then the colimit $\widehat \tau_{m}$ of the $\widehat \tau_{m,n}$'s satisfies $\text{H}''_{m}$. Finally, if $\text{H}''_{m}$ is satisfied for all $m$, then the colimit $\widehat \tau$ of the $\widehat \tau_{m}$'s satisfies the conditions of the theorem.

Let $N$ denote the global lower bound on $X$. To complete the proof, we need to show that
\begin{enumerate}
\item $\text{H}_{N,0,1}$ holds.
\item $\text{H}_{m,n,p}$ implies $\text{H}_{m,n,p+1}$ for all $m\geq N$, $n\in \mathbb N$ and all $p\geq 1$;
\item if $\text{H}_{m,n,p}$ holds for all $p$, then $\text{H}'_{m,n}$ holds;
\item $\text{H}'_{m,n}$ implies $\text{H}_{m,n+1,1}$ for all $m\geq N$, $n\in \mathbb N$; and
\item $\text{H}''_{m}$ implies $\text{H}'_{m+1,0,1}$ for all $m\geq N$.
\end{enumerate}

\emph{Proof of (1).}  Observe that $E^0\op V\circ \op P_{\star}\cong\op P_{\star}$, so that for all $d\in \ob \cat D$,
$$\op K_{N,0,1}(d)=\op T\big(X(d)_{N}\big)\acirc \widetilde F^1\Phi(\op P_{\star}).$$
Define $\widehat \tau'(d):X(d)_{N}\otimes (\op S\otimes \op P_{\star})\otimes \op A^\perp_{1}\to Y$ by
$$\widehat \tau'(d)(x\otimes s_{0}e_{\op P}\otimes \alpha _{1}):=\tau(d)(x),$$
which is differential, since $\tau(d)$ is. Extend $\widehat\tau'(d)$ to a transposed tensor morphism of right $\op A$-modules
$$\widehat\tau _{N,0,1}:\op T\big(X(d)_{N}\big)\acirc \widetilde F^1\Phi (E^n\op V\circ \op P_{\star})\to \op T\big(Y(d)\big),$$
which clearly satisfies $\text{H}_{0,0,1}$.

\medskip

\emph{Proof of (2).} Suppose that $\text{H}_{m,n,p}$ holds, for some $m,n\in \mathbb N$ and $p\geq 1$.  If there is no $\mathfrak m\in \mathfrak M$ such that $\deg x_{\mathfrak m}=m$, then $\op K_{m,n,p}(d)=\op K_{m-1,n,p}(d)$ for all $d\in \ob \cat D$ and for all $n,p$, which implies that $\text{H}_{m,n,p+1}$ is equivalent to $\text{H}_{m-1,n,p+1}$ for all $n,p$. Since $\text{H}_{m-1,n,p+1}$ follows from $\text{H}_{m,n,p}$, we therefore obtain for free in this case that $\text{H}_{m,n,p+1}$ holds.

Suppose therefore that there exists $\mathfrak m\in \mathfrak M$ such that $\deg x_{\mathfrak m}=m$.   Let $v\in B_{n,k}$, and let $\vec l\in I_{k}^{p+1}$, with $l=\sum _{i}l_{i}$.  Consider
$$w_{\mathfrak m,\vec l}\,(v)=x_{\mathfrak m}\otimes e_{\op A}\otimes \big(s_{k-1}\tilde v\otimes \alpha _{\vec l}\big)\otimes e_{\op A}^{\otimes l}\in X(\mathfrak m)\otimes \widetilde F^{p+1}\Phi(\op M).$$
Since $(\op M, \Delta)$ is semifree and both the cosimplicial and the the simplicial parts of the differential on $\Phi(\op M)$ are filtration-lowering, $\del _{\mathfrak m}w_{\mathfrak m,\vec l}\,(v)\in \op K_{m,n,p}(\mathfrak m)$.  Its image under $\widehat\tau_{m,n,p}(\mathfrak m)$,
$$\widehat\tau_{m,n,p}(\mathfrak m)\big(\del _{\mathfrak m}w_{\mathfrak m,\vec l}\,(v)\big)\in \op T\big(Y(\mathfrak m)\big),$$
is thus already defined.  On the other hand, since $Y(\mathfrak m)$ is acyclic, $\op T\big (Y(\mathfrak m)\big)$ is levelwise acyclic.  The cycle $\widehat\tau_{m,n,p}(\mathfrak m)\big(\del _{\mathfrak m}w_{\mathfrak m,\vec l}\,(v)\big)$ must therefore be a boundary, i.e., there exists $\omega _{\mathfrak m,\vec l}\,(v)\in \op T\big(Y(\mathfrak m)\big)(l)$ such that $\del_{\mathfrak m}'\omega _{\mathfrak m,\vec l}\,(v)=\widehat\tau_{m,n,p}(\mathfrak m)\big(\del _{\mathfrak m}w_{\mathfrak m,\vec l}\,(v)\big)$.

For any $d\in \ob\cat D$, consider the following morphism of symmetric sequences of chain complexes.
$$\widehat\tau_{m,n,p}(d)+\omega(d):\op K_{m,n,p}(d)\oplus \big[X(d)_{m}\otimes (\op S\otimes E^{n}\op V\circ \op P_{\star})\circ \op A_{p+1}^\perp\big]\longrightarrow \op T(Y),$$
where $\omega (d)$ is defined as follows. If $x\in X(d)_{m}$ can be expressed in the basis
$$B_{m}=\{X(f)(x_{\mathfrak m})\mid \mathfrak m\in \mathfrak M, \deg x_{\mathfrak m}=m, f\in \cat D(\mathfrak m,d)\}$$
as $x=\sum _{B_{m}}a_{\mathfrak m,f}X(f)(x_{\mathfrak m})$, where $a_{\mathfrak m,f}\in R$ for all $\mathfrak m$ and $f$, and $v\in E^{n}\op V$ is an element of the fixed basis, then for all $\vec l\in I_{k}^{p+1}$,
$$\omega(d)\Big(x\otimes  \big(s_{k-1}\tilde v\otimes \alpha _{\vec l}\otimes e_{\op A}^{\otimes p+1}\big)\Big)=\sum _{B_{m}}a_{\mathfrak m,f}\cdot Y(f)\big(\omega _{\mathfrak m,\vec l}\,(v)\big).$$
By Proposition \ref{prop:tt-morph}, we see that $\widehat\tau_{m,n,p}(d)+\omega(d)$ naturally induces a transposed tensor map of right $\op A$-modules
$$\widehat\tau_{m,n,p+1}(d):\op K_{m,n,p+1}(d)\longrightarrow \op T\big(Y(d)\big)$$
such that the restriction of $\widehat\tau_{m,n,p+1}(d)$ to $\op K_{m,n,p}(d)$ agrees with $\widehat\tau_{m,n,p}(d)$.

\medskip

\emph{Proof of (3).}  Let
$$\widehat\tau _{m,n}'(d):\op T\big(X(d)_{<m}\big)\acirc \Phi(\op M)+\op T\big(X(d)_{m}\big)\acirc \Phi (F^{n-1}\op M\oplus E^{n}\op V\circ \op P_{\star})\to \op T\big(Y(d)\big)$$
denote the colimit of the $\widehat\tau _{m,n,p}$'s.
Applying the induction functor, we obtain a natural, mutiplicative and differential morphism of symmetric sequences
$$\ind \widehat\tau_{m,n}'(d):\op T\big(\Om X(d)_{\leq m}\big)\circ (F^{n-1}\op M\oplus E^{n}\op V\circ\op P_{\star})\longrightarrow \op T\big(\Om Y(d)\big)$$
such that the restriction of $\ind\widehat\tau_{m,n}'(d)$ to $\op T\big(\Om X(d)_{\leq m}\big)\circ (F^{n-1}\op M)$ agrees with $\ind\widehat\tau _{m,n-1}(d)$.  Furthermore, $\ind \widehat\tau_{m,n}'(d)$ is determined by the family of multiplicative, equivariant maps of chain complexes
$$\{(\ind \widehat\tau_{m,n}'(d))_{k}:\Om X(d)_{\leq m}\otimes \big(F^{n-1}\op M(k)\oplus E^{n}\op V(k)\big)\to \big(\Om Y(d)\big)^{\otimes k}\}_{k}.$$
From this family we derive a new family of multiplicative, equivariant maps of chain complexes
$$\psi'(d)\big((\ind\widehat\tau_{m,n}'(d))_{k}\otimes Id_{\op P[\vec l]}\big): \Om X(d)_{\leq m}\otimes E^{n}\op V(k)\otimes \op P[\vec l]\to \big(\Om Y(d)\big)^{\otimes l}$$
for all  $k, l\in \mathbb N$, $\vec l\in I_{k,l}$,
which gives rise to a natural multiplicative morphism of symmetric sequences of graded abelian groups
$$\varphi(d):\op T\big(\Om X(d)_{\leq m}\big)\circ (E^{n}\op V\circ \op P)\to \op T(\Om Y(d)).$$
Applying Lemma \ref{lem:coprod} to $\varphi(d)$ and $\ind\widehat\tau_{m,n}(d)$, we obtain a multiplicative morphism of symmetric sequences,
$$\varphi(d)*\ind\widehat\tau_{m,n}'(d):\op T\big(\Om X(d)_{\leq m}\big)\circ (F^{n-1}\op M\oplus E^{n}\op V\circ \op P)\to \op T(\Om Y(d)),$$
which agrees with $\ind\widehat\tau_{m,n}'(d)$ on $\op T\big(\Om X(d)_{\leq m}\big)\circ (F^{n-1}\op M\oplus E^{n}\op V\circ\op P_{\star})$ and is therefore differential, since $\ind\widehat\tau_{m,n}'(d)$ and the action of $\op P$ on $\op T\big(\Om Y(d)\big)$ are differential maps.  Recall that $F^{n}\op M=F^{n-1}\op M\oplus E^{n}\op V\circ \op P$.

To conclude the proof that $\text{H}'_{m,n}$ holds, let $\widehat\tau_{m,n}(d)$ be equal to
\begin{multline*}
\widehat\tau_{m}(d)+\lin\big(\varphi(d)*\ind\widehat\tau'_{m,n}(d)\big):\\
\op T\big(X(d)_{<m}\big)\acirc \Phi(\op M)+\op T\big(X(d)_{\leq m}\big)\acirc \Phi (F^{n}\op M)\longrightarrow \op T\big(Y(d)\big).
\end{multline*}
\medskip

\emph{Proof of (4).}  As in the proof of (2),  if there is no $\mathfrak m\in \mathfrak M$ such that $\deg x_{\mathfrak m}=m$, there is nothing to do. Suppose therefore that there exists $\mathfrak m\in \mathfrak M$ such that $\deg x_{\mathfrak m}=m$.

Let $v\in B_{n+1,k}$.  Consider
$$w_{\mathfrak m}\,(v)=x_{\mathfrak m}\otimes e_{\op A}\otimes \big(s_{k-1}\tilde v\otimes \alpha _{1}^{\otimes k}\big)\otimes e_{\op A}^{\otimes k}\in X(\mathfrak m)\otimes \widetilde F^{1}\Phi( E^{n+1}\op V\circ \op P_{\star}).$$

Since $(\op M, \Delta)$ is semifree and both the cosimplicial and the the simplicial parts of the differential on $\Phi(\op M)$ are filtration-lowering, $\del _{\mathfrak m}w_{\mathfrak m}\,(v)\in \op K_{m,n}$.  Its image under $\widehat\tau_{m,n}(\mathfrak m)$,
$$\widehat\tau_{m,n}(\mathfrak m)\big(\del _{\mathfrak m}w_{\mathfrak m,\vec l}\,(v)\big)\in \op T\big(Y(\mathfrak m)\big),$$
is thus already defined.  On the other hand, since $Y(\mathfrak m)$ is acyclic, $\op T\big (Y(\mathfrak m)\big)$ is levelwise acyclic.  The cycle $\widehat\tau_{m,n}(\mathfrak m)\big(\del _{\mathfrak m}w_{\mathfrak m}\,(v)\big)$ must therefore be a boundary, i.e., there exists $\omega _{\mathfrak m}\,(v)\in \op T\big(Y(\mathfrak m)\big)(l)$ such that $\del_{\mathfrak m}'\omega _{\mathfrak m}\,(v)=\widehat\tau_{m,n}(\mathfrak m)\big(\del _{\mathfrak m}w_{\mathfrak m}\,(v)\big)$.  In particular, when $k=1$, our induction hypotheses imply that we can choose $\omega _{\mathfrak m}(v)=\tau (x_{\mathfrak m})$, which we do.

For any $d\in \ob\cat D$, consider the following morphism of symmetric sequences of chain complexes.
$$\widehat\tau_{m,n}(d)+\omega(d):\op K_{m,n}(d)\oplus \big[X(d)_{m}\otimes (\op S\otimes E^{n+1}\op V\circ \op P_{\star})\circ \op A_{1}^\perp\big]\longrightarrow \op T\big(Y(d)\big),$$
where $\omega (d)$ is defined as follows. If $x\in X(d)_{m}$ can be expressed in the basis
$$B_{m}=\{X(f)(x_{\mathfrak m})\mid \mathfrak m\in \mathfrak M, \deg x_{\mathfrak m}=m, f\in \cat D(\mathfrak m,d)\}$$
as $x=\sum _{B_{m}}a_{\mathfrak m,f}X(f)(x_{\mathfrak m})$, where $a_{\mathfrak m,f}\in R$ for all $\mathfrak m$ and $f$, and $v\in E^{n+1}\op V$ is an element of the fixed basis, then
$$\omega(d)\Big(x\otimes  \big(s_{k-1}\tilde v\otimes \alpha _{1}^{\otimes k}\otimes e_{\op A}^{\otimes k}\big)\Big)=\sum _{B_{m}}a_{\mathfrak m,f}\cdot Y(f)\big(\omega _{\mathfrak m}\,(v)\big).$$
By Proposition \ref{prop:tt-morph} , we see that $\widehat\tau_{m,n}(d)+\omega(d)$ naturally induces a transposed tensor map of right $\op A$-modules
$$\widehat\tau'(d):\op T\big(X(d)_{\leq m}\big)\acirc \widetilde F^{1}\Phi(F^n\op M\oplus E^{n+1}\op V\circ \op P_{\star})\longrightarrow \op T\big(Y(d)\big)$$
such that the restriction of $\widehat\tau'(d)$ to $\op T\big(X(d)_{\leq m}\big)\acirc \Phi(F^n\op M)$ agrees with $\widehat\tau_{m,n}(d)$.

To conclude the proof that $\text{H}_{m,n+1,1}$ holds, let $\widehat\tau_{m,n+1,1}(d)$ be equal to
$$
\widehat\tau_{m,n}(d)+\widehat\tau'(d):\op K_{m,n+1,1}\longrightarrow \op T\big(Y(d)\big).
$$
\medskip
\emph{Proof of (5).}  The argument in this case is almost identical to that in the proof of (4).

\end{proof}

The proof of the next theorem is very similar to that of Theorem \ref{thm:diffmod} and therefore left to the reader.

\begin{thm}\label{thm:diffmod2}  Let $X: \cat D\to\cat {F}_{\op P}$ and $Y:\cat D\to \coalg{A}$ be functors, where $\cat D$ is a category admitting a set of models $\frak M$ with respect to which  $X$ is free and globally connective and $Y$ is acyclic.  Let $(\op M, \Delta)$ be a semifree level comonoid in the category of left $\op
P$-modules under $\op J$.  Let $\tau: UX\to UY$
be a natural transformation, where $U$ is the forgetful functor down to $\cat {Ch}$.  Then there is a natural, multiplicative  transformation
$$\theta:\op T(\Om X)\underset {\op P}\circ \op M\to \op T(\Om Y)$$
extending $s^{-1}\tau$,  i.e., the following composite is equal to $\tau$.
$$X\xrightarrow{s^{-1}}s^{-1}X\hookrightarrow \op T(\Om X)\circ \op J\to \op T(\Om X)\underset {\op P}\circ \op M\xrightarrow{\theta}\op T(\Om Y)\xrightarrow {\text{proj.}}s^{-1}Y\xrightarrow {s}Y$$
\end{thm}

\begin{cor}\label{cor:acycmod} Let $X,Y: \cat D\to\cat F_{\op P}$ be functors, where $\cat D$ is a category admitting a set of models $\frak M$ with respect to which  $X$ is free and globally connective and $Y$ is acyclic.  Let $(\op M, \Delta)$ be a semifree level comonoid in the category of $\op P$-bimodules under $\op J$.  Let $\tau:UX\to UY$ be a natural transformation, where $U $ is the forgetful functor down to $\cat {Ch}$.  Then there is a natural, multiplicative  transformation of right $\op P$-modules
$$\theta:\op T(\Om X)\underset {\op P}\circ \op M\to \op T(\Om Y)$$
extending $s^{-1}\tau$,  i.e., the following composite is equal to $\tau$.
$$X\xrightarrow{s^{-1}}s^{-1}X\hookrightarrow \op T(\Om X)\circ \op J\to \op T(\Om X)\underset {\op P}\circ \op M\xrightarrow{\theta}\op T(\Om Y)\xrightarrow {\text{proj.}}s^{-1}Y\xrightarrow {s}Y$$
\end{cor}

Applying this corollary to the case $\op M=\op F=\Phi(\op J)$, we obtain a significant special case.

\begin{cor}\label{cor:last} Let $X, Y: \cat D\to\cat F_{\op A}$  be functors, where $\cat D$ is a category admitting a set of models $\frak M$ with respect to which  $X$ is free and globally connective and $Y$ is acyclic.  Let $\tau:UX\to UY$ be a natural transformation, where $U $ is the forgetful functor down to $\cat {Ch}$.  Then there is a natural, multiplicative  transformation
$$\Theta:\Om^2 X\to \Om^2 Y$$
extending $s^{-1}\tau$,  i.e., the following composite is equal to $\tau$.
$$X\xrightarrow{s^{-1}s^{-1}}s^{-1}(s^{-1}X)\hookrightarrow \Om ^2 X\xrightarrow{\Theta}\Om ^2 Y\xrightarrow {\text{proj.}}s^{-1}(s^{-1}Y)\xrightarrow {ss}Y$$
\end{cor}

\begin{proof}  Note that for all objects $d$ in $\cat D$, the cobar constructions $\Om X(d)$ and $\Om Y(d)$ admit natural, coassociative coproducts that are morphisms of algebras, since $X(d)$ and $Y(d)$ are objects in $\cat F_{\op A}$.  It therefore makes sense to apply the cobar construction again, to $\Om X(d)$ and $\Om Y(d)$.

By Corollary \ref{cor:acycmod}, there is a natural, multiplicative  transformation of right $\op A$-modules
$$\theta:\op T(\Om X)\underset {\op A}\circ \op F\to \op T(\Om Y)$$
extending $s^{-1}\tau$.  Applying induction to $\theta$, we obtain a multiplicative transformation
$$\ind(\theta):\op T(\Om^2 X)\circ \op J\to \op T(\Om ^2 Y).$$
We can then take $\Theta$ to be the restriction of $\ind(\theta)$ to level $1$.
\end{proof}

\begin{rmk}  One can certainly dualize everything done here to develop Enriched Bar Duality and Enriched Induction for bar duality.
\end{rmk}

\providecommand{\bysame}{\leavevmode\hbox
to3em{\hrulefill}\thinspace}
\providecommand{\MR}{\relax\ifhmode\unskip\space\fi MR }
\providecommand{\MRhref}[2]{%
  \href{http://www.ams.org/mathscinet-getitem?mr=#1}{#2}
} \providecommand{\href}[2]{#2}

\end{document}